\documentclass[a4paper,11pt,
 leqno
]{article}

\usepackage{comment}
\usepackage[utf8]{inputenc}
\usepackage[T1]{fontenc}
\usepackage{amsthm}
\usepackage{needspace}
\usepackage{amssymb}
\usepackage{bbm}
\usepackage{mathrsfs}
\usepackage{amsmath}
\usepackage{amsfonts} 
\usepackage{mathtools}
\usepackage{graphicx} 
\usepackage{float}
\usepackage{dsfont}
\usepackage{enumerate}
\usepackage{enumitem}
\usepackage[a4paper]{geometry}
\usepackage{csquotes}
\usepackage{IEEEtrantools}
\usepackage{appendix}
\usepackage{constants}
\usepackage[obeyDraft]{todonotes}
\usepackage{tikz-cd}
\usetikzlibrary{decorations.pathreplacing,decorations.pathmorphing}
\usetikzlibrary{patterns}
\usepackage{pgfplots}
\pgfplotsset{compat=1.18}
\usepackage{enumerate}
\usepackage{mathabx}
\usepackage{booktabs}
\usepackage
{hyperref}
\hypersetup{linktocpage,
  colorlinks   = true, 
  urlcolor     = blue, 
  linkcolor    = blue,
  citecolor   = blue 
}\usepackage[all]{hypcap}  
\usepackage{breakurl}

\usepackage{titlesec}

\titleformat{\subsection}[runin]   % run-in: no line break after the heading
  {\normalfont\normalsize\bfseries}
  {\thesubsection}
  {1em}
  {}
  [.]                               % period after the run-in title

\titlespacing*{\subsection}
  {0pt}      % left indentation
  {1.5ex}    % space before
  {1em}      % horizontal space after the run-in title

\numberwithin{equation}{section}

\usepackage{caption}
\title{
\normalsize
{\textbf{{THE FATE OF UNIVERSAL TIGHTNESS FOR  CRITICAL \\ GAUSSIAN FREE FIELD CLUSTERS 
}}}}
\author{}
\date{}
\makeatletter
\newcommand{\1}{\mathbbm 1}
\newcommand{\Z}{\mathbb Z}
\newcommand{\R}{\mathbb R}
\renewcommand{\P}{\mathbb{P}}
\newcommand{\E}{\mathbb{E}}

\renewcommand{\epsilon}{\varepsilon}

\newcommand{\capacity}{\text{cap}}
\newcommand{\variance}{\text{var}}
\newcommand{\cov}{\text{cov}}

\newcommand{\Fi}{\Theta}

\newtheorem{theorem}{Theorem}[section]
\newtheorem{lemma}[theorem]{Lemma}
\newtheorem{proposition}[theorem]{Proposition}
\newtheorem{corollary}[theorem]{Corollary}
\theoremstyle{definition}

\newtheorem{remark}[theorem]{Remark}

\newcommand{\floor}[1]{\left\lfloor{#1}\right\rfloor}

\newcommand{\cref}[1]{\text{\upshape\ref{#1}}}

\newcommand{\Slab}{\mathbb{S}_N}

\newcommand{\G}{\mathcal{G}}
\newcommand{\tG}{\tilde{\mathcal{G}}}
\definecolor{InkBlue}{RGB}{0.00,0.00,0.00} % fountain-pen blue, marks changes of notation
\newcommand{\nn}[1]{{\color{InkBlue}#1}}
\newcommand{\K}{\mathscr{C}}

\newcommand{\tSlab}{\tilde{\mathbb{S}}_N}
\newcommand{\Fbox}{F}

\newcommand{\LU}{\text{LocUniq}}

\newcommand{\diam}{\text{diam}}

\newcommand{\tB}{\tilde{B}}
\newcommand{\ttB}{\tilde{B}^\dagger}

\newcommand{\rad}{\text{rad}}

\newcommand{\LV}{M}
\newcommand{\Dcross}{h_N^\star}

\newcommand{\I}{{\cal I}}

\renewcommand{\tilde}{\widetilde}
\renewcommand{\hat}{\widehat}
\renewcommand{\epsilon}{\varepsilon}
\newcommand{\sv}{\underline{v}}
\newcommand{\bv}{\bar{v}}

\newconstantfamily{c}{symbol=c}

\newcommand{\defeq}{\stackrel{\text{def.}}{=}}

\usepackage{color}
\definecolor{Red}{rgb}{1,0,0}
\definecolor{Blue}{RGB}{90,145,120} % darkened Cambridge blue (original #A3C1AD = 163,193,173), used for marking changes via \blue
\definecolor{Olive}{rgb}{0.41,0.55,0.13}
\definecolor{Yarok}{rgb}{0,0.5,0}
\definecolor{Green}{rgb}{0,1,0}
\definecolor{MGreen}{rgb}{0,0.8,0}
\definecolor{DGreen}{rgb}{0,0.55,0}
\definecolor{Yellow}{rgb}{1,1,0}
\definecolor{Cyan}{rgb}{0,1,1}
\definecolor{Magenta}{rgb}{1,0,1}
\definecolor{Orange}{rgb}{1,.5,0}
\definecolor{Violet}{rgb}{.5,0,.5}
\definecolor{Purple}{rgb}{.75,0,.25}
\definecolor{Brown}{rgb}{.75,.5,.25}
\definecolor{Grey}{rgb}{.7,.7,.7}
\definecolor{Black}{rgb}{0,0,0}

\def\black{\color{Black}}

\makeatletter
\newcommand{\namedlabel}[2]{%
  \begingroup
  \phantomsection
  \def\@currentlabel{#1}%
  \label{#2}%
  \endgroup
}
\makeatother

\usepackage{titletoc}

\dottedcontents{section}[4em]{}{2.9em}{0.7pc}
\dottedcontents{subsection}[0em]{}{3.3em}{1pc}

\usepackage{multicol}
\usepackage{parcolumns}

\usepackage{tabu}
\usepackage{caption}
\begin{document}
\maketitle
\vspace{0.1cm}
\begin{center}
\vspace{-1.9cm}
Wen Zhang

\end{center}
\vspace{0.1cm}
\begin{abstract}
\centering
\begin{minipage}{0.90\textwidth}
{
We study critical level-set percolation of the Gaussian free field on the cable graph of two-dimensional slabs, a family of geometries ranging from the planar case to three dimensions.
For Bernoulli percolation on any graph, the volume of the largest cluster in a ball is tight on the scale of its mean at every value of the parameter $p$ (including at $p=p_c$).
We show that this universal tightness fails for the free field: at local scales where the slab looks flat, the volume of the largest cluster is of the order of its mean with probability $\frac12$, and of strictly smaller order otherwise.
We identify the limiting law explicitly: for instance, for the planar field, the fraction of a mesoscopic ball occupied by the largest cluster converges to $\frac12\delta_0+\frac12\operatorname{Unif}[0,1]$.
When the geometry is effectively three-dimensional, tightness is restored. We identify a sharp criterion that governs this dichotomy. As it turns out, flatness already sets in at scales much smaller than the slab thickness.
Along the way, we derive up-to-constant cluster volume estimates, which exhibit (iterated-)logarithmic corrections, and show that, despite this unconventional scaling behaviour, the critical scaling relations hold uniformly across all slab heights and scales.
}
\end{minipage}
\end{abstract}
\thispagestyle{empty}

\vspace{8.5cm}
\begin{flushright}

\noindent\rule{5cm}{0.4pt} \hfill 27 November 2026 \\
\bigskip
\end{flushright}

\begin{multicols}{2}
\noindent Imperial College London\\
 Department of Mathematics\\
 %180 Queen's Gate\\
 London SW7 2AZ, UK \\
  \url{kz2623@ic.ac.uk} 

\columnbreak
\begin{flushright}

\end{flushright}
\end{multicols}

\newpage

\section{Introduction}
A remarkable result of Hutchcroft \cite[Theorem~2.2]{hutchcroft_2021_power} asserts that for Bernoulli
percolation on an arbitrary countable weighted graph, the volume $M_r$ of the
largest cluster in a box of radius $r$ is tight on the scale of its mean at every
value of the density $p$. This property is known as \emph{universal tightness}.
Results of this type were previously known in specific settings, notably for
critical bond percolation on $\Z^2$ \cite{borgs_1999_unifo,borgs_2001_the} (see also \cite[Theorem~8]{kesten_1986_the} for
the corresponding statement for the incipient infinite cluster), on
high-dimensional tori \cite{heydenreich_2007_rando} and on the complete graph \cite{bollobas_1984_the,aldous_1997_brown}. For the cable-graph
Gaussian free field $\varphi$, the same tightness result for critical clusters
was established in \cite{drewitz_2024_clust,cai_2025_incip} on a class of transient graphs with
polynomial volume growth and polynomial decay of the Green function, provided the
critical one-arm probability decays like the square root of the Green function,
which includes $\Z^d$, $3\leq d\leq 5$.

\smallskip

We show that, in contrast, tightness fails near the lower critical dimension. On the cable graph of $\Z^2$ with killing rate $N^{-2}$, and more generally on slabs $\Z^2\times(\Z/h_N\Z)$ ranging from the planar case to three dimensions, $M_r$ can be of the order of its mean with asymptotic probability exactly $\frac12$, and of strictly smaller order otherwise. This happens precisely when $\Dcross\ll r\ll N$, where the crossover scale $\Dcross\approx h_N/\log(N/h_N)$ lies well below the slab height (Theorem~\ref{thm:antiCon}). In this regime, we identify the limiting law of $M_r$ explicitly (Theorem~\ref{thm:limit}). For instance, if $h_N=1$ and $r=\sqrt N$, the density of the largest cluster in $B_r$ converges in law to $\frac12\delta_0+\frac12\operatorname{Unif}[0,1]$. In particular, for small $\varepsilon$, both $\P_N(M_r\geq\varepsilon\E_N[M_r])$ and $\P_N(M_r\geq3\varepsilon\E_N[M_r])$ are close to $\frac12$, whereas for Bernoulli percolation the BK inequality bounds the latter by the square of the former. This multiplicative decay underlies the universal tightness in~\cite{hutchcroft_2021_power}.

The mechanism at play here is a random shift of the level, cf.~Section~\ref{sec:proofIdea} and Figure~\ref{fig:tools}. Inside $B_r$, the field splits into a field killed outside a larger ball and an independent harmonic part generated by the field further away, which is nearly constant on $B_r$. Disjoint large pieces of cluster in $B_r$ are thus positively correlated through this common shift, which defeats the BK argument. The shift matters only when it exceeds the critical window at scale $r$, which happens exactly when $\Dcross\ll r\ll N$; the field then looks locally supercritical or subcritical with probability $\frac12$ each. Notably, the variance of the field does not decide the outcome: the variance diverges when $h_N=1$ and $r\asymp N$, where tightness holds; and the variance stays bounded when $h_N=\sqrt N$, where tightness can fail. The previous heuristic accounts for the asymptotic $(\tfrac12, \tfrac12)$ split (see~\eqref{eq:anticonBig}-\eqref{eq:anticonSmall}), while identifying the full limiting law requires substantially more work, which we describe below.

Our results also yield up-to-constant bounds for the expected size of the largest cluster and for the volume tail of the cluster of the origin (see \eqref{eq:thmM0rExp} and Theorem~\ref{thm:volTailMain}), which are not power laws. Combined with our previous work \cite{rodriguez_2025_on} and the two-point function estimate of \cite{lupu_2016_from}, they show that, despite the tightness anomaly and this unconventional scaling, the critical (hyper)scaling relations persist when formulated directly in terms of the observables, as in Kesten's work \cite{kesten_1987_scali}; see Table~\ref{t:critExp} and Remark~\ref{rem:classicalvgeneralscaling}.
\medskip

We now describe our setup and results more precisely. 
We study percolation of level sets of the Gaussian free field on two-dimensional slabs with spatial scale $N\geq 1$ and height $1\leq h_N\leq N$. This setting was introduced in our previous work~\cite{rodriguez_2025_on} as a rigorous analogue of what physicists refer to as the ``$\epsilon$-expansion'' (cf.~\cite{wilson_1974_the}) around the lower critical dimension, corresponding here to $\epsilon=d-2$. 
To make the notion of slabs precise,
let ${h}\!:\![1,\infty)\!\mapsto\!\R$ be a non-decreasing function satisfying $1\leq {h}(t)\leq t$ for all $t\!>\!0$. For $N\!\geq\!1$, let
\begin{equation}\label{def:slab}
    h_N\defeq \floor{{h}(N)},
    \qquad\text{and}\qquad
    \mathbb{S}_N\defeq \mathbb{Z}^2\times(\Z/h_N\Z).
\end{equation}
We refer to $\mathbb{S}_N$ as a (cylindrical) slab. The periodic boundary condition in the vertical direction is imposed only for convenience; one may instead consider free boundary conditions, { for which one might expect similar results to hold}. For $x=(x_1,x_2,x_3)\in\mathbb S_N$, we let
$\|x\|\coloneqq |(x_1,x_2,\widehat x_3)|,$
where $|\cdot|$ is the Euclidean norm in $\R^3$ and $\widehat x_3\in \Z$ is the representative of $x_3\in\mathbb Z/h_N\mathbb Z$ satisfying
$|\widehat x_3|=d_{\mathbb Z/h_N\mathbb Z}(0,x_3)\leq \frac{h_N}{2}.$
For $r\geq1$, we set $B_r\coloneqq\{x\in\mathbb S_N:\|x\|< r\}$, { $B(x,r)\coloneqq x+B_r$, $x\in\Slab$}.

We consider $\varphi=(\varphi_x)_{x \in \Slab}$, the centered Gaussian field with law $\P_N^h$, $N \geq 1$, whose covariance is given by the Green's function $g_N(\cdot,\cdot)$ of the continuous-time simple random walk $X$ on $\Slab$ with unit jump rate to nearest neighbors and killing rate $N^{-2}$; see Section~\ref{sec:facts} for the precise setup. We write $\E_N[\cdot]$ for the expectation corresponding to $\P_N$. The law $\P_N$ is invariant under translations of $\Slab$. By \cite[Theorem 2.1]{rodriguez_2025_on},
\begin{equation}\label{eq:gN0}
    g_N \defeq g_N(0,0)=\E_N[\varphi_0^2]
    \sim \frac{3}{\pi}\frac{\log N}{h_N}+g_{\Z^3},
    \qquad \text{as } N\to\infty,
\end{equation}
where $\sim$ means that the ratio of both sides tends to one in the limit and
$ g_{\Z^3}$ denotes the on-diagonal Green's function of the random walk on $\Z^3$.

The field $\varphi$ extends to a continuous field on the corresponding metric graph, or cable system $\tSlab$; see Section~\ref{sec:facts} for the precise setup. Our results concern the connectivity of the random set
\[
\{\varphi\geq 0\}\defeq \{x\in \tSlab:\varphi_x\geq 0\},
\]
which, as shown in \cite[Theorem 1.1]{rodriguez_2025_on}, manifests critical behavior. For $U,V\subset\tSlab$, we write $\{U\leftrightarrow V\}$ for the event that $U$ and $V$ are connected by a continuous path in $\{x\in \tSlab:\varphi_x\geq 0\}$. For $U\subset\tSlab$, we write $|U|$ for the cardinality of $U\cap\Slab$, which we refer to as the volume of $U$, and set
$\operatorname{rad}(U)\defeq \sup_{x\in U}\|x\|$.
The quantity $M_r$ alluded to above, which will be of primary interest, is defined as follows. For $r\geq 1$, possibly depending on $N$, let 
\begin{equation}\label{def:df}
  M_r \defeq \sup\{|\K|: \text{$\K$ is a connected component of $\{\varphi\geq 0\}\cap \tB_r$}\}, \qquad \overline{M}_r\defeq\frac{M_r}{\E_N[M_r]},
\end{equation}
where $\tB_r\subset\tSlab$ denotes the cable graph analogue of $B_r$, obtained by including all cables whose endpoints lie in $B_r${, and $\tB(x,r)$ is defined similarly}. We refer to~\eqref{def:tBr} for its precise definition.

\subsection{\texorpdfstring{{ Anticoncentration and mixture limits}}{Anticoncentration and mixture limits}}

As we will see, the dichotomy between concentration and anticoncentration is governed by a ratio of capacities, which we now introduce. Informally, capacity measures how visible a set is to the underlying random walk $X$ introduced above \eqref{eq:gN0}; see Section~\ref{sec:facts} for the precise definition. 
For $A\subset\Slab$, we write $\capacity(A)$ for the capacity of $A$ associated with $X$, and $\capacity^r(A)$ for the corresponding capacity when $X$ is additionally killed upon hitting $\partial B_r$. Note that $\capacity(\cdot)$ and $\capacity^r(\cdot)$ implicitly depend on $N$ through $X$.
We will soon be considering the limiting behaviour of $M_r$ as $N\to\infty$, with $r$ possibly tending to infinity with $N$. 
For reasons that will become clear shortly, we consider the condition
\begin{equation*}\label{def:AntiCon}
  \tag{AntiCon}
  \lim_{N\to\infty}\frac{\capacity^{2r}(B_r)}{\capacity(B_r)}=\infty.
\end{equation*}
%We denote by $(\mathrm{Con})$\namedlabel{Con}{def:Con} the complementary regime, in which the ratio $\capacity^{2r}(B_r)/\capacity(B_r)$ remains bounded as $N\to\infty$. 
{ We denote by $(\mathrm{Con})$\namedlabel{Con}{def:Con} the regime in which $\limsup_{N}\capacity^{2r}(B_r)/\capacity(B_r)<\infty$, and we then write $K\defeq\sup_N\capacity^{2r}(B_r)/\capacity(B_r)$. If the ratio $\capacity^{2r}(B_r)/\capacity(B_r)$ has a limit in $[0,\infty]$ as $N\to\infty$, then \eqref{def:Con} and \eqref{def:AntiCon} form a perfect dichotomy. In general, every sequence $r=r_N$ has a subsequence along which exactly one of the two conditions holds;}
{ see \eqref{eq:introAntiCon2} for an explicit description of these regimes in terms of $r$, $h_N$ and $N$.}

In the sequel, for non-negative quantities $a$ and $b$, we write $a\asymp b$ if there exist constants $0<c\leq c'<\infty$ such that $cb\leq a\leq c'b$ uniformly in all varying parameters.  We write $a_N\ll b_N$ to mean $\lim_N\frac{a_N}{b_N}=0$ and $a \lesssim b$ to mean $a\leq c\,b$, where $c>0$ is a universal constant independent of all relevant parameters. We will also use the notation $a \vee b =\max\{ a,b \}$ and $a \wedge b =\min \{a,b\}$. Let $\Cl[c]{range}\in[1,\infty)$ be an arbitrary fixed constant.

\begin{theorem}\label{thm:antiCon} 
 For all $N\geq 1$ and $1\leq r\leq \Cr{range}N$, with $\overline{M}_r$ as in \eqref{def:df}, the following hold:
\begin{enumerate}[label=(\roman*)]
  \item { Under \eqref{def:Con},} for all $1\leq t\leq c\,\E_N[M_r]$,
  \begin{equation}\label{eq:tailEstlargestVolgena=0}
  \P_N(t^{-1}\leq \overline{M}_r \leq t )\geq 1-ct^{-c'},
  \end{equation}
  { where $c,c'\in(0,\infty)$ depend only on $K$ and $\Cr{range}$.}
  \item { If $r=r_N$ satisfies} \eqref{def:AntiCon},
  \begin{align}
  \label{eq:anticonBig}
  \lim_{t\to\infty}\liminf_{N}
  \P_N(t^{-1}\leq \overline{M}_r\leq t)
  &=
  \lim_{t\to\infty}\limsup_{N}
  \P_N(t^{-1}\leq \overline{M}_r\leq t)
  =\frac12,\\
  \label{eq:anticonSmall}
  \lim_{t\to\infty}\liminf_{N}
  \P_N(\overline{M}_r < t^{-1})
  &=
  \lim_{t\to\infty}\limsup_{N}
  \P_N(\overline{M}_r < t^{-1})
  =\frac12.
  \end{align}
\end{enumerate}
\end{theorem}                               
Thus, under \eqref{def:AntiCon}, asymptotically half of the distribution of $M_r$ lies on the scale of its expectation, while the other half is supported on strictly smaller scales.

\medskip

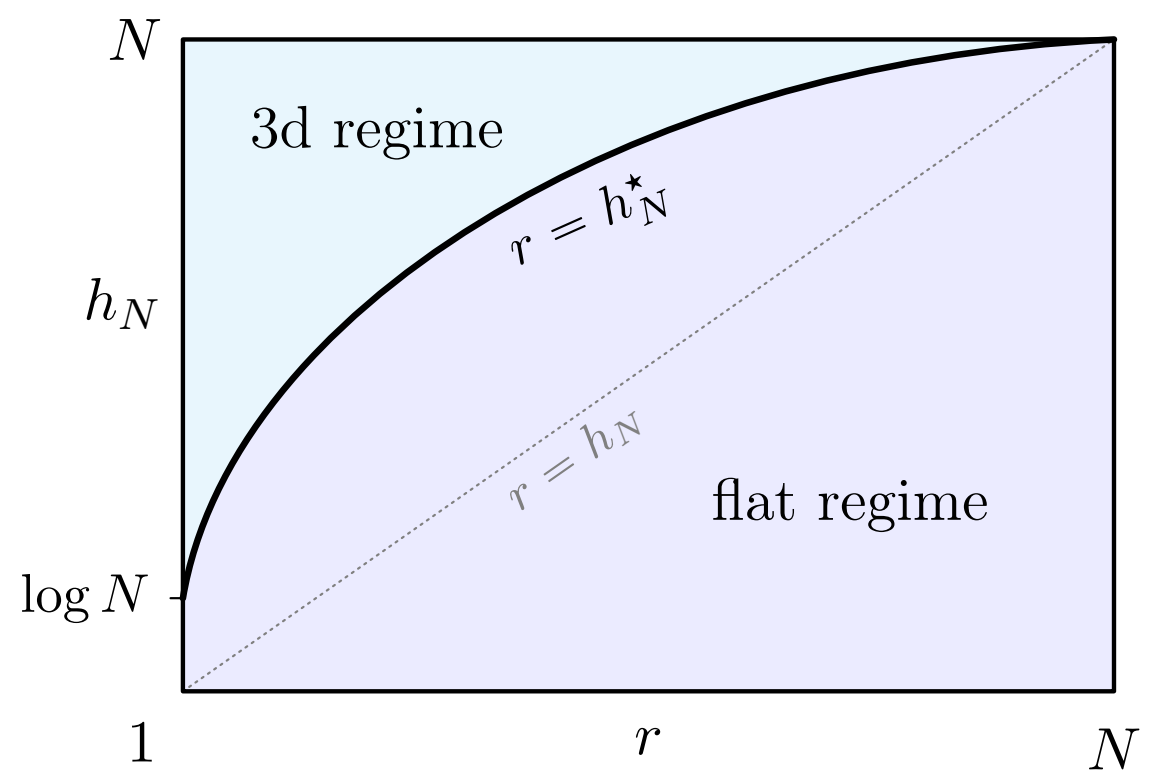
\begin{figure}[t]
\centering
\resizebox{0.48\textwidth}{!}{%
\begin{tikzpicture}[line cap=round, line join=round]

\def\xmax{6}
\def\ymax{4.2}
% Schematic, not to scale. The crossover curve r = h_N^star meets the vertical axis at
% h_N ~ log N and bends concavely towards the upper-right corner.
\def\crosscurve{(0,0.6) .. controls (0.3,2.4) and (2.8,4.05) .. (\xmax,\ymax)}

% Flat regime
\path[fill=blue!8] (0,0) -- (0,0.6)
  -- \crosscurve -- (\xmax,0) -- cycle;
% 3d regime
\path[fill=cyan!8] (0,\ymax) -- (0,0.6)
  -- \crosscurve -- cycle;

% diagonal r = h_N below the merging point, for reference
\draw[densely dotted, gray] (0,0) -- (\xmax,\ymax);
\node[gray, rotate=35, anchor=west, font=\footnotesize] at (2,1.15) {$r=h_N$};

% crossover curve
\draw[very thick] \crosscurve;
\node[rotate=24, anchor=west, font=\small] at (2,2.75) {$r=h_N^\star$};

% bounding box
\draw[thick] (0,0) rectangle (\xmax,\ymax);

% axis labels and ticks
\node[anchor=north, yshift=-3pt] at (0.5*\xmax,0) {$r$};
\node[left] at (0,0.5*\ymax+0.4) {$h_N$};
\node[anchor=north east, xshift=-1pt, yshift=-2pt] at (0,0) {$1$};
\node[anchor=north, yshift=-3pt] at (\xmax,0) {$N$};
\node[left] at (0,\ymax) {$N$};
\draw (0,0.6) -- (-0.08,0.6) node[left, font=\small] {$\log N$};

% regime labels
\node at (1.25,3.6) {\normalsize 3d regime};
\node at (4.3,1.2) {\normalsize flat regime};

\end{tikzpicture}%
%[Previous drawing of Figure 1:]
%\begin{tikzpicture}
%
%\def\xmax{6}
%\def\ymax{4.2}
%
%% Flat regime
%\path[fill=blue!8]
%    (0,0)
%    .. controls (0.45,1.15) and (1.25,1.95) .. (2.4,2.55)
%    .. controls (3.55,3.15) and (5.00,3.85) .. (\xmax,\ymax)
%    -- (\xmax,0)
%    -- cycle;
%
%% 3d regime
%\path[fill=cyan!8]
%    (0,0)
%    -- (0,\ymax)
%    -- (\xmax,\ymax)
%    .. controls (5.00,3.85) and (3.55,3.15) .. (2.4,2.55)
%    .. controls (1.25,1.95) and (0.45,1.15) .. (0,0)
%    -- cycle;
%
%% Crossover curve
%\draw[very thick]
%    (0,0)
%    .. controls (0.45,1.15) and (1.25,1.95) .. (2.4,2.55)
%    .. controls (3.55,3.15) and (5.00,3.85) .. (\xmax,\ymax);
%
%% Bounding box
%\draw[thick] (0,0) rectangle (\xmax,\ymax);
%
%% Axis labels
%\node[anchor=north, yshift=-3pt] at (0.5*\xmax,0) {$r$};
%\node[left] at (0,0.5*\ymax) {$h_N$};
%
%% Endpoint labels
%\node[anchor=north east, xshift=-1pt, yshift=-2pt] at (0,0) {$1$};
%\node[anchor=north, yshift=-3pt] at (\xmax,0) {$N$};
%\node[left] at (0,\ymax) {$N$};
%
%% Regime labels
%\node at (1.45,3.35) {\normalsize 3d regime};
%\node at (4.45,1.45) {\normalsize flat regime};
%
%% Boundary label
%\node[rotate=26] at (3.75,3.55) {$r=h_N^\star$};
%
%\end{tikzpicture}%
}
\caption{The crossover between the three-dimensional and flat regimes { (not to scale)}. { For $h_N\lesssim \log N$, every scale $r\geq1$ lies in the flat regime.} When $r=N$, the model is always in the flat regime, whereas when $h_N=N$, the geometry is effectively three-dimensional. The anticoncentration appears exactly in the flat regime where $r$ additionally satisfies $\Dcross\ll r\ll N$.}
\label{fig:crossover}
\end{figure}

Our analysis leading to Theorem~\ref{thm:antiCon} also yields the precise order of $\E_N[M_r]$. To state this, we write $\log^\star(x)\defeq \max\{\log x,1\}$ and introduce the crossover scale
\begin{equation}\label{def:Dcross}
    \Dcross\defeq \frac{h_N}{\log^\star(N/h_N)}.
\end{equation}
Then, for all $1\leq r\leq \Cr{range}N$, we have
\begin{equation}\label{eq:thmM0rExp}
	\E_N[M_r] \asymp \bv(r)\defeq 
    \begin{cases}
        r^{5/2},
        & r\leq \Dcross,
        \\[1.5ex]
        \displaystyle
        r^2(r\wedge h_N)\Big(
        \frac{g_N h_N}{\log^\star\!(N/r)}
        \Big)^{-1/2},
        & r>\Dcross.
    \end{cases}.
  \end{equation}
%{ By the one-arm estimate of \cite[Theorem~1.1]{rodriguez_2025_on} (see \eqref{eq:one-arm} below) and since $|B_r|\asymp r^2(r\wedge h_N)$, one has $\bar v(r)\asymp |B_r|\,\P_N(\operatorname{rad}(\K)\geq r)$ in both regimes.}
{(Here and below, the letter $v$ stands for volume.)} The role of $\Dcross$ is particularly transparent in \eqref{eq:thmM0rExp}. Below this scale, that is, when $r$ lies a factor $\log^\star\!(N/h_N)$ below the slab height $h_N$, $M_r$ matches the { three-dimensional scaling of} \cite[Theorem 1.1]{drewitz_2024_clust} and is insensitive to the slab height. Once $r$ exceeds $\Dcross$, the cluster begins to feel the finite extent of the third direction, and $h_N$ enters explicitly. We refer to this latter regime as the \emph{flat} regime. See Figure~\ref{fig:crossover} for a visualisation.

{ For instance, take $h_N=1$ and $r=\sqrt N$, as in the example at the beginning of this section. Then \eqref{eq:thmM0rExp} gives $\E_N[M_r]\asymp N\asymp|B_r|$, while the capacity estimates of \cite{rodriguez_2025_on} give $\capacity^{2r}(B_r)/\capacity(B_r)\asymp\log N$, so that \eqref{def:AntiCon} holds. Similarly, for $h_N=1$ and $r=N$, \eqref{def:Con} holds and Theorem~\ref{thm:antiCon} together with \eqref{eq:thmM0rExp} shows that $M_N$ is tight on the scale $N^2/\sqrt{\log N}$, in line with the prediction of \cite{cao_2021_level} that critical clusters of the two-dimensional Gaussian free field are logarithmic fractals, and with the continuum results of \cite{aru_2020_first,aru_2020_the}, where first passage sets of the two-dimensional continuum Gaussian free field and Brownian loop soup clusters are shown to have non-trivial Minkowski content in the gauge $r^2|\log r|^{1/2}$.}

\smallskip
The~\eqref{def:AntiCon} and~\eqref{def:Con} dichotomy is closely tied to the distinction between the three-dimensional and flat regimes. Indeed, a short calculation (see Lemma~\ref{lem:equivFrh}) shows that
\begin{equation}\label{eq:introAntiCon2}
    \eqref{def:AntiCon}\text{ holds if and only if } \Dcross\ll r\ll N.
\end{equation}
Thus, \eqref{def:AntiCon} essentially corresponds to the flat regime, with the additional proviso that $r\ll N$. In other words, anticoncentration occurs only when clusters are observed at a \emph{local} scale within the flat regime. This is a manifestation of a phenomenon we refer to as \emph{killing sensitivity}, whereby restricting a cluster to a smaller scale $r$, or imposing boundary conditions outside $\tB_r$, can substantially alter its connectivity properties; cf.~the right picture in Figure~\ref{fig:tools}. Handling this distinction is a principal difficulty in our analysis of $M_r$ (more on this below). { By contrast, killing sensitivity is absent on $\Z^d$, $d\geq3$: there $\capacity_{\Z^d}(B_r)\asymp\capacity^{2r}_{\Z^d}(B_r)\asymp r^{d-2}$, so the analogue of \eqref{def:Con} always holds, in line with the tightness of $M_r$ established in \cite{drewitz_2024_clust}. The same holds on $\Slab$ when $h_N=N$, since then $\capacity(B_r)\asymp\capacity^{2r}(B_r)\asymp r$.

\medskip

A closer look at~\eqref{eq:anticonSmall} reveals that, under~\eqref{def:AntiCon}, every subsequential limiting law of $\overline{M}_r$ has an atom of mass $\frac12$ at $0$.
It is then natural to wonder if the full limiting law of $M_r$ can be identified.
To state our next result, let $g_N^r(\cdot,\cdot)$ denote the Green's function when $X$ is additionally killed upon hitting $\partial B_r$. Write $e_1=(1,0,0)\in\Slab$ for the first horizontal coordinate vector, and set
\[
\sigma_{N,r}\defeq\sqrt{\frac{g_N(0,re_1)}{g_N^r(0,0)}}.
\]
Let $Z$ be a standard normal random variable, $\Phi$ its distribution function, and $Z_+\defeq \max\{Z,0\}$. 

\begin{theorem}\label{thm:limit}
  If $1\leq r=r_N\leq \Cr{range}N$ satisfies \eqref{def:AntiCon}, $\lim_Nr=\infty$  and $\lim_N\sigma_{N,r}=\sigma$, then
  \begin{equation}\label{eq:limit_density}
  \begin{aligned}
    &\frac{M_r}{|B_r|}\stackrel{d}{\longrightarrow}\big(2\Phi(\sigma Z)-1\big)_+, &&\text{if }\sigma\in(0,\infty],\\[1ex]
    &\frac{M_r}{\sigma_{N,r}|B_r|}\stackrel{d}{\longrightarrow}\sqrt{\tfrac{2}{\pi}}\,Z_+, &&\text{if }\sigma=0,
  \end{aligned}
  \end{equation}
  where $\Phi(\infty\cdot z)\defeq\mathbf 1_{\{z>0\}}$ for $z\neq0$.
\end{theorem}
We refer to Theorem~\ref{thm:limit_precise} for a stronger version of Theorem~\ref{thm:limit} which does not assume the existence of $\sigma$. It also yields the sharp asymptotics
\begin{equation}\label{eq:EcvgClaim}
  \E_N[M_r]\sim |B_r|\frac{1}{\pi}\arctan(\sigma_{N,r}),
\end{equation}
under \eqref{def:AntiCon} and $\lim_Nr=\infty$ only, which refines \eqref{eq:thmM0rExp} in this regime.
Heuristically, $\sigma Z$ corresponds to the harmonic shift $\eta_0^{2r}$ from the outside of $\tB_{2r}$, normalised by the standard deviation $\sqrt{g_N^{2r}(0,0)}$ of the killed field $\psi^{2r}$, and $2\Phi(\sigma Z)-1$ is the probability that $|\psi^{2r}_x|\leq\eta_0^{2r}$, which is the asymptotic density of the largest cluster of $\{\psi^{2r}\geq-\eta_0^{2r}\}\cap\tB_r$; see Section~\ref{sec:proofIdea} and \eqref{eq:intro_supcrit_lim} below.

For critical Bernoulli site percolation on the triangular lattice, \cite[Theorem 5]{camia_2019_confo} establishes a limiting law for $\overline{M}_r$ that is supported on $(0,\infty)$, has no atoms, and has a stretched exponential upper tail.
Theorem~\ref{thm:limit}, which is valid beyond the planar ($h_N=1$) setting throughout the \eqref{def:AntiCon} regime when $r\ll N$, shows that the limiting law in our setting retains an atom
of mass $\frac12$ at $0$ and admits an explicit description in
terms of a single Gaussian variable. 

To unpack Theorem~\ref{thm:limit},
let $g^{\Z^3}(0)$ denote the on-diagonal Green's function of the simple random walk on $\Z^3$. Adapting the arguments of \cite[Theorem 2.1 and Proposition 5.1]{rodriguez_2025_on}, one obtains that, under the assumption of Theorem~\ref{thm:limit},
\begin{equation}\label{eq:sigma}
  \sigma_{N,r}^2\sim \frac{\frac{3}{\pi}\frac{1}{h_N}\log(N/r)}{\frac{3}{\pi}\frac{1}{h_N}\log(r)+g_{\Z^3}(0)} .
\end{equation}
We do not use \eqref{eq:sigma} in our proofs; it only serves to make the examples below explicit. The parameter $\sigma_{N,r}$ governs the density of $M_r$ through~\eqref{eq:EcvgClaim}. As one might expect, $\arctan(\sigma_{N,r})$ is comparable to the one-arm probability to scale $r$, given below in~\eqref{eq:one-arm}.
As $\sigma_{N,r}$ varies from $0$ to $\infty$, $M_r$ passes
from a regime of vanishing density to one where it occupies almost
all of $B_r$ with probability $\frac12$.
We now illustrate these regimes through several special cases. For a random variable $X$, we denote by $\mathscr{L}(X)$ the law of $X$.
\begin{itemize}
  \item \emph{Low density.}
  If $h_N\gg\log N$, then $\sigma_{N,r}\to0$, and
  Theorem~\ref{thm:limit} together with~\eqref{eq:sigma} yields
  \begin{equation}\label{eq:low_density_lim}
    \E_N[M_r]\sim
    |B_r|\sqrt{\frac{3}{\pi^3}}
    \sqrt{\frac{\log(N/r)}{g_{\Z^3}(0)h_N}}
    \quad\text{and}\quad
    \overline{M}_r\stackrel{d}{\longrightarrow}\sqrt{2\pi}Z_+.
  \end{equation}
  The constant $\sqrt{3/\pi^3}$ appearing in \eqref{eq:low_density_lim} is no stranger to the present model. It also governs the one-arm asymptotics in the ``thin'' regime; see~\cite[Corollary~1.4]{rodriguez_2025_on}. The thinness condition imposed in~\cite[Corollary~1.4]{rodriguez_2025_on} is stronger than
\eqref{def:AntiCon}, and~\eqref{eq:low_density_lim} suggests that it could perhaps be relaxed to~\eqref{def:AntiCon}.
  \item \emph{Intermediate density.}
  If $r=\sqrt N$ and $h_N\ll\log N$, then
  $\sigma_{N,r}\to1$, and
  \begin{equation}\label{eq:intermediate_density_lim}
    \E_N[M_r]\sim\frac14|B_r|
    \quad\text{and}\quad
    \frac{M_r}{|B_r|}\stackrel{d}{\longrightarrow}X,
    \qquad
    \mathscr{L}(X)=\frac12\delta_0+\frac12\operatorname{Unif}[0,1].
  \end{equation}

  \item \emph{Full density with probability $\frac12$.}
  If $h_N\ll\log N$ and $\log r\ll\log N$, then
  $\sigma_{N,r}\to\infty$, and
  \begin{equation}\label{eq:full_density_lim}
    \E_N[M_r]\sim\frac12|B_r|
    \quad\text{and}\quad
    \frac{M_r}{|B_r|}\stackrel{d}{\longrightarrow}Y,
    \qquad
    \mathscr{L}(Y)=\frac12\delta_0+\frac12\delta_1.
  \end{equation}
\end{itemize}

\medskip

\subsection{\texorpdfstring{{A double-edged sword: killing sensitivity}}{A double-edged sword: killing sensitivity}}\label{sec:proofIdea}
As one might expect, a central difference between working with $\Slab$ and higher dimensional settings lies in the slow decay of correlation. 
One manifestation of this is the killing sensitivity phenomenon (and hence the~\eqref{def:AntiCon} and~\eqref{def:Con} dichotomy in Theorem~\ref{thm:antiCon}), already discussed
around~\eqref{eq:introAntiCon2}.
It turns out that this sensitivity under~\eqref{def:AntiCon} is a double-edged sword. On the one hand, the direct effect is that local or killed objects cannot in general be replaced by their unrestricted counterparts, which rules out several strategies available in higher dimensions. This is one of the main difficulties in the derivation of Theorem~\ref{thm:antiCon} and~\eqref{eq:thmM0rExp}.
On the other hand, the discrepancy between killed and unkilled objects is so pronounced that the appropriate analogues of unkilled critical objects turn out to be killed objects at an off-critical level.
Working away from criticality provides us with the additional control needed to identify the precise limiting laws in Theorem~\ref{thm:limit}.

\smallskip

We now explain this in more detail, and refer to Figure~\ref{fig:tools} for an accompanying illustration.
Let us first introduce, for any continuous function $\phi:\tSlab\to\R$, $a\in\R$ and $r\geq1$,
\begin{equation}\label{def:Mr_general}
M_r^a(\phi)\defeq
\sup\big\{|\K|:\K\text{ is a connected component of }
\{\phi\geq a\}\cap\tB_r\big\}.
\end{equation}
We omit the superscript when $a=0$ and the argument when
$\phi=\varphi$. Thus $M_r^a=M_r^a(\varphi)$ and $M_r=M_r^0(\varphi)$. We let $\psi^r$ denote the Gaussian free field killed outside $\tB_r$, see~\eqref{eq:field_decomp} for its precise definition.

\medskip
\noindent\textbf{The problem.}\quad
We start by illustrating the complications arising from killing
sensitivity through the proofs of Theorem~\ref{thm:antiCon}
and~\eqref{eq:thmM0rExp}.
The upper bounds in~\eqref{eq:thmM0rExp} follow from the two-point function estimate supplied by \cite[Proposition 5.2]{lupu_2016_from} and \cite[Theorem 2.1]{rodriguez_2025_on}. Upper tightness of $\overline{M}_r$ is immediate from
Markov's inequality.
The main work therefore concerns the lower bounds.
Following naively the strategy developed in \cite{drewitz_2024_clust}, the starting point is to show that
\begin{equation}\label{eq:proofstep1}
  \P_N(M_r\geq c\,\bv(r))\geq c',
\end{equation}
with $\bv(\cdot)$ as in \eqref{eq:thmM0rExp}. We refer to Proposition~\ref{pro:manyMeso} for the precise near-critical statement underlying \eqref{eq:proofstep1}.
The main difficulty in a statement like \eqref{eq:proofstep1} is the mismatch between the ambient scale governing the field (which is $N$) and the scale of the  clusters under consideration (which is $r$). As we explain in Section~\ref{sec:mainprop}, if the goal were instead to show
\begin{equation}\label{eq:proofstep1killed}
  \P_N(M_r(\psi^r)\geq c\sv(r))\geq c,
\end{equation}
and $\sv(r) \asymp \E_N[M_r(\psi^r)]$ is the killed analogue of $\bv(r)$, this mismatch would disappear, and the proof would essentially be identical to the argument of \cite{drewitz_2024_clust}.
Killing sensitivity is precisely the obstruction to passing
from~\eqref{eq:proofstep1killed} to~\eqref{eq:proofstep1}.
In higher dimensions, and more generally in the setting of~\cite{drewitz_2024_clust}, the additional killing preserves the connectivity estimates relevant to $M_r$, allowing $\varphi$ and $\psi^r$ to be used interchangeably for this purpose.
On $\Slab$, this fails under \eqref{def:AntiCon}, so that passing from~\eqref{eq:proofstep1killed} to~\eqref{eq:proofstep1} requires a new approach.

\begingroup
\colorlet{sensitivityRed}{red!65!black}
\colorlet{sensitivityBlue}{blue!75!black}
\colorlet{sensitivityGray}{gray!85}
\colorlet{sensitivityShade}{cyan!9}
\tikzset{sensitivity stroke/.style={line width=0.75pt}}
\tikzset{killing sensitivity/.style={
  line cap=round,
  line join=round,
  loop trace/.style={
    sensitivity stroke,
    decorate,
    decoration={random steps, segment length=1.45pt, amplitude=0.67pt}
  },
  small loop trace/.style={
    loop trace,
    decoration={random steps, segment length=0.70pt, amplitude=0.24pt}
  },
  largest cluster/.style={sensitivityRed},
  smaller cluster/.style={sensitivityBlue},
  largest local piece/.style={}
}}

% The same drawing is used before and after clipping to B_r.
\def\drawambientloops{%
  \pgfmathsetseed{4291}%
  \foreach \trace/\tracestyle in {%
    {(-2.05,0.34) (-1.96,0.72) (-1.77,0.84) (-1.81,1.12)
     (-1.56,1.42) (-1.27,1.38) (-1.12,1.65) (-0.89,1.52)
     (-1.00,1.26) (-0.88,1.06) (-1.09,0.86) (-0.96,0.69)
     (-0.80,0.72) (-0.71,0.61) (-0.53,0.58) (-0.51,0.44)
     (-0.37,0.42) (-0.31,0.25) (-0.41,0.13) (-0.59,0.19)
     (-0.68,0.08) (-0.80,0.17) (-0.99,0.18) (-1.10,0.31)
     (-1.28,0.45) (-1.52,0.35) (-1.62,0.14) (-1.93,0.17)}/,
    {(-1.98,0.45) (-2.13,0.34) (-2.09,0.08) (-2.23,-0.11)
     (-2.05,-0.28) (-2.08,-0.49) (-1.94,-0.64) (-1.76,-0.59)
     (-1.70,-0.35) (-1.82,-0.22) (-1.61,-0.04) (-1.57,0.17)
     (-1.74,0.33)}/,
    {(-1.31,1.65) (-1.12,1.84) (-0.93,1.82) (-0.79,2.02)
     (-0.49,1.98) (-0.28,1.76) (-0.11,1.82) (0.08,1.68)
     (0.19,1.51) (0.43,1.56) (0.60,1.33) (0.40,1.17)
     (0.13,1.22) (-0.05,1.39) (-0.32,1.21) (-0.54,1.38)
     (-0.58,1.62) (-0.80,1.64) (-1.03,1.40) (-1.23,1.45)}/,
    {(0.12,1.42) (0.05,1.69) (0.35,1.89) (0.60,1.78)
     (0.79,1.88) (1.04,1.79) (1.02,1.56) (1.29,1.53)
     (1.34,1.34) (1.63,1.23) (1.66,1.02) (1.40,0.87)
     (1.23,1.05) (0.95,1.03) (0.75,1.29) (0.54,1.22)
     (0.44,1.43) (0.26,1.31)}/,
    {(1.29,1.23) (1.62,1.30) (1.81,1.10) (1.74,0.89)
     (1.94,0.69) (1.87,0.48) (2.10,0.33) (2.07,0.06)
     (1.84,-0.02) (1.73,-0.27) (1.50,-0.34) (1.37,-0.12)
     (1.16,-0.22) (0.96,-0.30) (0.75,-0.22) (0.67,-0.05)
     (0.46,-0.15) (0.30,-0.08) (0.28,0.07) (0.41,0.15)
     (0.34,0.33) (0.49,0.42) (0.66,0.33) (0.72,0.42)
     (0.93,0.37) (1.11,0.33) (1.23,0.19) (1.45,0.12)
     (1.41,0.33) (1.24,0.48) (1.28,0.68) (1.48,0.76) (1.34,0.98)}/,
    {(1.77,0.02) (1.96,-0.14) (1.84,-0.39) (1.98,-0.55)
     (1.83,-0.73) (1.80,-0.97) (1.57,-1.00) (1.54,-1.20)
     (1.32,-1.39) (1.12,-1.48) (0.88,-1.63) (0.74,-1.41)
     (0.90,-1.25) (0.87,-1.08) (1.04,-0.91) (1.31,-0.88)
     (1.38,-0.65) (1.60,-0.56) (1.54,-0.33) (1.60,-0.19)}/,
    {(1.12,-1.44) (0.98,-1.82) (0.68,-1.75) (0.52,-2.02)
     (0.26,-2.06) (0.13,-1.89) (-0.09,-1.95) (-0.32,-1.80)
     (-0.42,-1.56) (-0.27,-1.41) (-0.36,-1.22) (-0.29,-1.01)
     (-0.39,-0.82) (-0.49,-0.63) (-0.44,-0.47) (-0.29,-0.42)
     (-0.24,-0.52) (-0.08,-0.51) (0.02,-0.38) (0.21,-0.44)
     (0.36,-0.39) (0.49,-0.53) (0.62,-0.52) (0.69,-0.77)
     (0.91,-0.84) (0.89,-1.08) (0.77,-1.23) (0.89,-1.38)}/largest local piece,
    {(-0.19,-1.69) (-0.38,-1.94) (-0.62,-1.88) (-0.77,-1.96)
     (-0.99,-1.80) (-0.93,-1.61) (-1.14,-1.47) (-1.06,-1.30)
     (-0.87,-1.36) (-0.73,-1.26) (-0.58,-1.42) (-0.36,-1.44)}/,
    % Smaller overlapping loops reinforce each of the three local pieces.
    {(-0.62,0.50) (-0.53,0.64) (-0.40,0.60) (-0.32,0.73)
     (-0.20,0.68) (-0.16,0.55) (-0.24,0.49) (-0.25,0.36)
     (-0.39,0.31) (-0.49,0.39) (-0.61,0.38)}/,
    {(-0.51,0.19) (-0.41,0.27) (-0.36,0.38) (-0.25,0.34)
     (-0.23,0.23) (-0.14,0.16) (-0.22,0.08) (-0.35,0.12)
     (-0.45,0.05) (-0.54,0.11)}/,
    {(0.38,0.24) (0.20,0.24) (0.13,0.35) (0.19,0.47)
     (0.31,0.43) (0.37,0.55) (0.53,0.57) (0.57,0.46)
     (0.50,0.35) (0.56,0.27) (0.48,0.18)}/,
    {(0.52,0.02) (0.61,0.13) (0.74,0.07) (0.81,0.14)
     (0.88,0.04) (0.85,-0.11) (0.78,-0.18) (0.66,-0.14)
     (0.59,-0.19) (0.48,-0.11)}/,
    {(-0.48,-0.57) (-0.57,-0.45) (-0.51,-0.33) (-0.35,-0.36)
     (-0.28,-0.27) (-0.17,-0.31) (-0.15,-0.46) (-0.25,-0.55)
     (-0.26,-0.67) (-0.40,-0.72)}/largest local piece,
    {(-0.16,-0.57) (-0.03,-0.49) (0.04,-0.35) (0.17,-0.30)
     (0.26,-0.39) (0.32,-0.49) (0.30,-0.63) (0.19,-0.68) (0.20,-0.79)
     (0.04,-0.84) (-0.02,-0.72) (-0.14,-0.71)}/largest local piece%
  }
    \draw[loop trace,\tracestyle] plot coordinates \trace -- cycle;
}

% Ambient soup with the components inside the observation ball highlighted.
\newcommand{\drawsensitivityambient}{%
\begin{scope}
  \fill[sensitivityShade] (0,0) circle (1.08);
  \begin{scope}[sensitivityGray]
    \drawambientloops
    % Other ambient clusters lying outside B_r.
    \pgfmathsetseed{6123}
    \foreach \trace in {%
      {(-1.80,-1.23) (-1.76,-1.10) (-1.68,-1.12) (-1.65,-1.03)
       (-1.53,-1.02) (-1.50,-1.13) (-1.39,-1.15) (-1.42,-1.24)
       (-1.52,-1.27) (-1.57,-1.38) (-1.68,-1.34) (-1.72,-1.25)},
      {(-1.55,-1.24) (-1.49,-1.18) (-1.40,-1.21) (-1.36,-1.31)
       (-1.43,-1.34) (-1.45,-1.43) (-1.57,-1.45) (-1.63,-1.38)
       (-1.58,-1.32)},
      {(-1.48,-0.33) (-1.49,-0.26) (-1.42,-0.23) (-1.39,-0.27)
       (-1.31,-0.29) (-1.33,-0.36) (-1.38,-0.40) (-1.44,-0.37)},
      {(-1.56,1.04) (-1.50,1.12) (-1.43,1.10) (-1.39,1.16)
       (-1.28,1.11) (-1.27,1.01) (-1.36,0.98) (-1.39,0.90)
       (-1.48,0.91) (-1.47,0.98)},
      {(-1.60,0.92) (-1.56,0.99) (-1.50,0.98) (-1.45,1.02)
       (-1.40,0.95) (-1.45,0.89) (-1.49,0.85) (-1.57,0.86)},
      {(-0.39,1.53) (-0.33,1.58) (-0.25,1.59) (-0.22,1.53)
       (-0.28,1.47) (-0.36,1.48)},
      {(0.10,2.13) (0.14,2.18) (0.20,2.17) (0.25,2.19)
       (0.27,2.13) (0.23,2.09) (0.16,2.08)},
      {(0.79,1.51) (0.83,1.58) (0.91,1.59) (0.93,1.54)
       (0.99,1.48) (0.94,1.42) (0.88,1.44) (0.82,1.43)},
      {(1.48,0.48) (1.49,0.57) (1.55,0.61) (1.61,0.55)
       (1.70,0.56) (1.73,0.47) (1.66,0.44) (1.61,0.38)
       (1.56,0.43)},
      {(1.16,-1.17) (1.19,-1.09) (1.25,-1.12) (1.32,-1.09)
       (1.36,-1.15) (1.33,-1.22) (1.25,-1.24) (1.20,-1.29)
       (1.15,-1.25)},
      {(0.15,-1.59) (0.21,-1.51) (0.19,-1.45) (0.25,-1.38)
       (0.31,-1.42) (0.39,-1.39) (0.46,-1.46) (0.40,-1.52)
       (0.43,-1.59) (0.34,-1.66) (0.26,-1.62) (0.21,-1.66)}%
    }
      \draw[small loop trace] plot coordinates \trace -- cycle;
  \end{scope}
  \begin{scope}[smaller cluster, largest local piece/.style={largest cluster}]
    \clip (0,0) circle (1.08);
    \drawambientloops
  \end{scope}
  % Small ambient clusters contained in B_r.
  \pgfmathsetseed{6149}
  \foreach \trace in {%
    {(-0.87,-0.19) (-0.84,-0.14) (-0.78,-0.16) (-0.72,-0.13)
     (-0.69,-0.20) (-0.73,-0.24) (-0.80,-0.26) (-0.85,-0.23)},
    {(0.06,0.78) (0.09,0.84) (0.16,0.85) (0.19,0.79)
     (0.15,0.74) (0.10,0.72)}%
  }
    \draw[small loop trace, smaller cluster] plot coordinates \trace -- cycle;
  \draw[semithick] (0,0) circle (2.34);
  \draw[densely dashed] (0,0) circle (1.08);
  \node[above left, font=\small] at (145:2.34) {$B_N$};
  \node[below left, font=\small] at (220:1.08) {$B_r$};
\end{scope}
}

% The killed-field cluster is formed by overlapping loops entirely in B_r.
\newcommand{\drawsensitivitykilled}{%
  \fill[sensitivityShade] (0,0) circle (1.08);
  \pgfmathsetseed{1387}
  \foreach \trace in {%
    {(-0.82,-0.29) (-0.71,-0.44) (-0.50,-0.33) (-0.41,-0.42)
     (-0.29,-0.29) (-0.11,-0.31) (-0.05,-0.17) (0.10,-0.21)
     (0.20,-0.05) (0.40,-0.01) (0.41,0.15) (0.57,0.26)
     (0.52,0.42) (0.32,0.36) (0.26,0.22) (0.07,0.27)
     (-0.01,0.11) (-0.19,0.16) (-0.22,0.04) (-0.37,0.07)
     (-0.42,-0.09) (-0.59,-0.01) (-0.66,-0.18) (-0.84,-0.12)},
    {(-0.44,-0.04) (-0.57,0.08) (-0.76,0.11) (-0.77,0.28)
     (-0.58,0.36) (-0.65,0.50) (-0.50,0.62) (-0.40,0.80)
     (-0.23,0.75) (-0.12,0.62) (-0.25,0.53) (-0.19,0.40)
     (-0.37,0.31) (-0.32,0.17) (-0.24,0.08) (-0.32,-0.02)},
    {(0.23,0.23) (0.19,0.39) (0.06,0.52) (0.13,0.68)
     (0.32,0.66) (0.41,0.82) (0.56,0.77) (0.64,0.63)
     (0.79,0.61) (0.76,0.42) (0.59,0.37) (0.60,0.23)
     (0.46,0.19) (0.40,0.31) (0.31,0.28)},
    {(-0.47,-0.25) (-0.62,-0.25) (-0.70,-0.43) (-0.61,-0.57)
     (-0.52,-0.56) (-0.48,-0.72) (-0.28,-0.77) (-0.22,-0.91)
     (-0.04,-0.93) (0.04,-0.76) (-0.10,-0.64) (-0.06,-0.45)
     (-0.22,-0.44) (-0.26,-0.30) (-0.37,-0.34)},
    {(0.32,0.10) (0.51,0.13) (0.63,0.02) (0.84,-0.01)
     (0.90,-0.18) (0.82,-0.33) (0.72,-0.30) (0.67,-0.48)
     (0.47,-0.55) (0.38,-0.69) (0.23,-0.55) (0.32,-0.41)
     (0.23,-0.26) (0.34,-0.09) (0.26,-0.02)},
    {(-0.10,-0.45) (0.05,-0.33) (0.24,-0.36) (0.28,-0.50)
     (0.43,-0.59) (0.38,-0.74) (0.22,-0.73) (0.14,-0.85)
     (0.01,-0.73) (0.03,-0.58) (-0.14,-0.58)},
    {(-0.33,0.66) (-0.28,0.83) (-0.13,0.83) (-0.09,0.91)
     (0.02,0.88) (0.06,0.77) (-0.03,0.68) (-0.15,0.70) (-0.22,0.57)},
    {(0.60,-0.41) (0.77,-0.43) (0.75,-0.55) (0.63,-0.57)
     (0.57,-0.72) (0.47,-0.64) (0.50,-0.48)}%
  }
    \draw[loop trace, largest cluster]
      plot coordinates \trace -- cycle;

  % Small clusters separated from the largest one.
  \pgfmathsetseed{2753}
  \foreach \trace in {%
    {(-0.90,0.46) (-0.89,0.53) (-0.84,0.56) (-0.79,0.53)
     (-0.80,0.46) (-0.86,0.43)},
    {(0.10,0.88) (0.13,0.93) (0.19,0.96) (0.23,0.93)
     (0.22,0.86) (0.17,0.84) (0.15,0.87)},
    {(0.86,0.22) (0.85,0.28) (0.90,0.31) (0.94,0.28)
     (0.95,0.22) (0.89,0.20)},
    {(0.84,0.30) (0.81,0.34) (0.84,0.39) (0.89,0.38)
     (0.91,0.33) (0.89,0.28)},
    {(-0.67,-0.67) (-0.61,-0.65) (-0.57,-0.69) (-0.60,-0.74)
     (-0.56,-0.77) (-0.61,-0.81) (-0.67,-0.77) (-0.69,-0.73)}%
  }
    \draw[small loop trace, smaller cluster] plot coordinates \trace -- cycle;
  \draw[semithick] (0,0) circle (1.08);
}

\begin{figure}[!htbp]
\centering
\resizebox{0.95\textwidth}{!}{%
\begin{tikzpicture}[x=0.85cm, y=1.14cm, baseline=0pt, line cap=round, line join=round]
\begin{scope}
  \fill[sensitivityShade] (-1,-2) rectangle (1,2);
  \draw[->] (-3.6,-2) -- (-3.6,2) node[above, font=\small] {$\eta^{\lambda r}$};
  \draw[->] (-3.6,0) -- (3.5,0) node[right,font=\small] {$x$};
  \node[font=\scriptsize, gray, anchor=east] at (-3.62,0) {$0$};
  \foreach \t/\l in {-1/-r, 1/r, -2.2/-\lambda r, 2.2/\lambda r}
    \draw (\t,0.07) -- (\t,-0.07) node[below, xshift=6pt, font=\scriptsize] {$\l$};
  \def\zlA{(-3.4,0.9) (-3.25,1.5) (-3.1,0.6) (-2.95,1.8) (-2.8,0.7) (-2.65,1.6) (-2.5,0.4) (-2.35,1.3) (-2.2,1.05)}
  \def\hmA{(-2.2,1.05) (-1.6,1.25) (-1,1.12) (-0.5,1.08) (0,1.1) (0.5,1.07) (1,1.1) (1.6,0.85) (2.2,1.0)}
  \def\zrA{(2.2,1.0) (2.35,1.7) (2.5,0.5) (2.65,1.4) (2.8,0.3) (2.95,1.5) (3.1,0.8) (3.25,1.6)}
  \def\zlB{(-3.4,-1.05) (-3.25,-0.65) (-3.1,-0.2) (-2.95,0.55) (-2.8,-0.35) (-2.65,0.25) (-2.5,-0.7) (-2.35,-0.1) (-2.2,-0.75)}
  \def\hmB{(-2.2,-0.75) (-1.6,-0.95) (-1,-0.82) (-0.5,-0.78) (0,-0.8) (0.5,-0.77) (1,-0.8) (1.6,-0.6) (2.2,-0.7)}
  \def\zrB{(2.2,-0.7) (2.35,0.05) (2.5,-0.6) (2.65,0.4) (2.8,-0.25) (2.95,0.2) (3.1,-0.35) (3.25,-0.95)}
  \begin{scope}[sensitivity stroke, sensitivityGray, densely dashed]
    \draw plot coordinates \zlB;
    \draw plot[smooth, tension=0.7] coordinates \hmB;
    \draw plot coordinates \zrB;
  \end{scope}
  \begin{scope}[sensitivity stroke, sensitivityBlue]
    \draw plot coordinates \zlA;
    \draw plot[smooth, tension=0.7] coordinates \hmA;
    \draw plot coordinates \zrA;
  \end{scope}
  \node[sensitivityBlue, font=\small] at (0,1.55) {$\eta^{\lambda r}>0$};
  \node[sensitivityGray, font=\small] at (0,-1.15) {$\eta^{\lambda r}<0$};
\end{scope}
\end{tikzpicture}%
\hspace{0.72cm}%
\begin{tikzpicture}[killing sensitivity, baseline=0pt]
  \drawsensitivityambient
  \draw[<-, >=stealth, sensitivityRed, sensitivity stroke]
    (0.43,-0.60) .. controls (0.62,-1.08) and (1.23,-1.49) .. (2.13,-1.47)
    node[right, font=\small] {$\K_r^{\max}(\varphi)$};
  \begin{scope}[shift={(3.68,0)}]
    \drawsensitivitykilled
    \node[below right, font=\small] at (330:1.08) {$B_r$};
    \draw[<-, >=stealth, sensitivityRed, sensitivity stroke]
      (0.33,0.84) .. controls (0.40,1.02) and (0.44,1.24) .. (0.36,1.49)
      node[above, font=\small] {$\K_r^{\max}(\psi^r)$};
  \end{scope}
\end{tikzpicture}%
}
\caption{
Left: the comparison tool. Inside $\tB_{\lambda r}$, the harmonic part $\eta^{\lambda r}$ generated by the field outside $\tB_{\lambda r}$ is almost constant on $\tB_r$ (shaded), so that it acts on $\tB_r$ as a random shift of the level; under \eqref{def:AntiCon} this shift is much larger than its oscillation over $\tB_r$.  Since $\eta^{\lambda r}$ is centred Gaussian, two cases occur, each with probability $\frac12$: in the constructive case (solid), the shift is positive and pushes the field in $\tB_r$ into a supercritical state; in the destructive case (grey), 
%the field far away looks alike, but 
the negative shift pulls the field in $\tB_r$ below level $0$ (dashed), killing the large local clusters. This is the dichotomy behind Theorem~\ref{thm:antiCon}.
Right: killing sensitivity. For a field $\phi$, we use $\K_r^{\max}(\phi)$ to denote the largest connected component in $\{\phi\geq 0\}\cap\tB_r$, so $|\K_r^{\max}(\phi)|=M_r(\phi)$. 
The right figure depicts what happens in the case when $\eta^{\lambda r}<0$ under~\eqref{def:AntiCon}. This shift
raises the level of the killed field beyond its critical window.
The cluster $\K_r^{\max}(\varphi)$ then exhibits subcritical
connectivity at scale $r$, while $\K_r^{\max}(\psi^r)$ is a critical cluster of the killed field. Under~\eqref{def:Con}, as in higher dimensions, the harmonic shift is at most of the order of the
critical window. Hence such a picture cannot appear, and both fields will exhibit critical connectivity at scale $r$.
}
\label{fig:tools}
\end{figure}
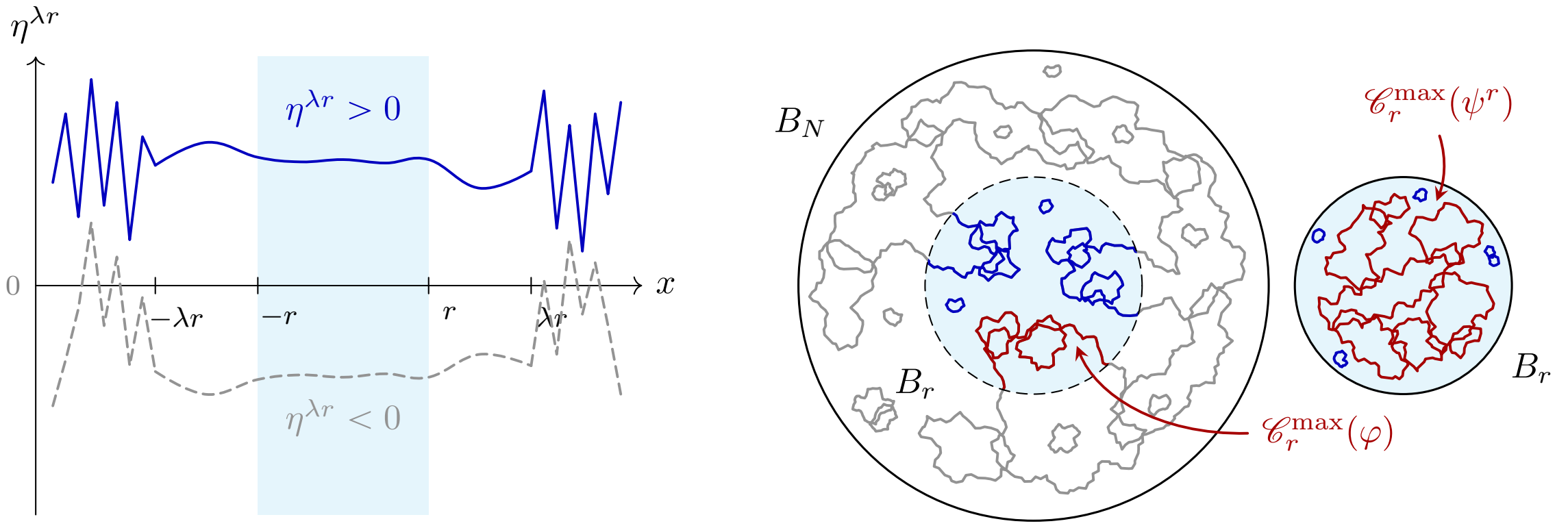
\endgroup

\medskip
\noindent\textbf{The remedy: two-sided comparison with random sprinkling.}
A natural way to address the killing sensitivity problem is to attempt a pathwise comparison between  the full field $\varphi$ and $\psi^r$, the field with Dirichlet boundary condition outside $\tB_r$. 

We provide such a comparison, at the cost of a random shift (sprinkling) of the level.
Let us recall the domain Markov property of $\varphi$ (cf.~\eqref{eq:field_decomp} for a precise description). For $\lambda\in(1,\infty]$,
\begin{equation}\label{eq:introDMP}
      \varphi=\psi^{\lambda r}+\eta^{\lambda r},
\end{equation}
where $\lambda\geq1$, $\psi^{\lambda r}$ is independent of the harmonic contribution $\eta^{\lambda r}$ generated by the field outside $\tB_{\lambda r}$.
Roughly speaking, for some $c\in(0,1)$ and all $1\leq r\leq N$,
%, with $b\asymp(r\wedge \frac{h_N}{\log^\star(N/r)})^{-1/2}$, 
one has,
\begin{equation}\label{eq:introshift}
    \P_N\!\big(
        M^{-\frac12 \eta_0^{\lambda r}}_r\!\!(\psi^{\lambda r})
        \leq
        M_r(\varphi)\,
    |\, \eta_0^{\lambda r}\!\geq\! 0\big)\geq c
    \quad\text{and}\quad
    \P_N\!\big(
        M^{-\frac12\eta_0^{\lambda r}}_r\!\!(\psi^{\lambda r})
        \geq
        M_r(\varphi)\,
    |\,\eta_0^{\lambda r}\!<\!0\big)\geq c.
\end{equation}
Thus, up to a random level shift of order $\eta_0^{\lambda r}$, local clusters of the unrestricted field can be compared in both directions with clusters of a field whose ambient scale is itself of order $r$. This comparison extends beyond cluster volume to arbitrary increasing observables of the local level sets. We refer to Proposition~\ref{prop:shift_h0} for the most precise and general statement.

Three features of this comparison are worth stressing. First, it is two-sided and holds for any family of increasing observables of the local level sets, such as cluster volumes, crossing events or capacities; it therefore reduces questions about local level sets of $\varphi$ to the same questions for a field whose ambient scale is $r$, at which point available techniques from dimensions $d\geq3$ apply. Second, it is uniform in $h_N$ and in $1\leq r\leq N$. 

The third feature is the most important, as it reveals the mechanism behind killing sensitivity. We refer to Figure~\ref{fig:tools} for a visual aid.
Two scales are in fact at play. The critical window of the clusters of $\psi^r$ is of order $a_{\mathrm{crit}}(r)\defeq(r\wedge h_N)^{-1/2}$, see Corollary~\ref{the:mainMra}. This means that the clusters $\{\psi^r\geq a\}$ are still critical for $|a|\lesssim a_{\mathrm{crit}}(r)$ but off-critical when $|a|\gg a_{\mathrm{crit}}(r)$. On the other hand, the typical size of $\eta^{\lambda r}_0$ is of order $a_{\mathrm{harm}}(r)\defeq(r\wedge \frac{h_N}{\log^\star(N/r)})^{-1/2}$, see~\eqref{eq:varEta0}. Moreover, it turns out that
\begin{equation}\label{eq:shift_gap}
  a_{\mathrm{harm}}(r)\,\,\begin{cases}
    \gg a_{\mathrm{crit}}(r), & \text{under }\eqref{def:AntiCon},\\
    \asymp a_{\mathrm{crit}}(r), & \text{under }\eqref{def:Con}.
  \end{cases}
\end{equation}
Together with~\eqref{eq:introshift}, the gap
$a_{\mathrm{harm}}(r)\gg a_{\mathrm{crit}}(r)$
under~\eqref{def:AntiCon} explains the mechanism behind killing
sensitivity. The ambient critical field, viewed locally,
behaves like the killed field at an off-critical level. On the other hand, the fact that $ M^{-\frac12\eta_0^{\lambda r}}_r\!(\psi^{\lambda r})$ remains critical under~\eqref{def:Con} explains exactly why killing sensitivity is absent in $d\geq 3$.

The starting point for the proof of \eqref{eq:introshift} is~\eqref{eq:introDMP}. The key is to show that, when $\lambda$ is sufficiently large, $\eta^{\lambda r}$ is almost constant throughout $\tB_r$ with high probability. Thus, on the event that $\eta^{\lambda r}_0$ is positive, (resp.~negative), the field $\varphi$ inside $\tB_r$ behaves approximately as $\psi^{\lambda r}$ with a uniformly positive (resp.~negative) shift; see the left picture in Figure~\ref{fig:tools}.
The main input needed to make this precise is Proposition~\ref{prop:gDrad}, which is a gradient estimate for the Green's function on $\Slab$. Its proof uses the decomposition of $g_N$ developed in \cite{rodriguez_2025_on}, which separates the contribution to $g_N$ into a topologically trivial part and a part that captures the periodicity in the third direction.
Proposition~\ref{prop:gDrad} yields a bound on the canonical metric of the Gaussian field
$\eta^{\lambda r}_x-\eta^{\lambda r}_0$, from which a standard entropy estimate (see, e.g. \cite[Section 1.3]{_2007_rando}) shows that its oscillation over $\tB_r$ is of order at most $(r\wedge h_N)^{-1/2}$. The delicacies of the slab geometry alluded to above again kick in when making all of this rigorous.

\medskip
\noindent\textbf{The other edge: access to limiting law.}
The gap $a_{\mathrm{harm}}(r)\gg a_{\mathrm{crit}}(r)$
under~\eqref{def:AntiCon} in~\eqref{eq:shift_gap} now works in our favour.
The off-critical description provided by~\eqref{eq:introshift} is exactly what allows us to identify the precise limiting law of $\overline{M}_r$.
Indeed, under the assumptions of Theorem~\ref{thm:limit}, we show in
Proposition~\ref{pro:convergences} that, for
$a\gtrsim  a_{\mathrm{harm}}(r)$,
\begin{equation}\label{eq:intro_supcrit_lim}
  \frac{M^{-a}_r(\psi^{2r})}
  {|B_r|\big(2\Phi(a/\sqrt{g_N^{2r}(0,0)})-1\big)}
  \stackrel{\P_N}{\longrightarrow}1.
\end{equation}
Thus, the largest supercritical cluster has an asymptotically
deterministic density, explicitly determined by its level.

Combining~\eqref{eq:intro_supcrit_lim} with the comparison \eqref{eq:introshift} 
identifies the asymptotic volume on $\{\eta^{2r}_0>0\}$.
On $\{\eta^{2 r}_0<0\}$,~\eqref{eq:introshift} shows that $M_r$ is dominated by the largest cluster of the killed field, whose expectation is of strictly smaller order than $\E_N[M_r]$ under \eqref{def:AntiCon}.
Since $\eta^{2r}_0$ is a centred Gaussian, each event has
probability $\frac12$. This explains both the value $\frac12$
in Theorem~\ref{thm:antiCon} and the mixture structure of the
limiting laws in Theorem~\ref{thm:limit}.

The proof of~\eqref{eq:intro_supcrit_lim} relies on the
isomorphism between the killed field and random interlacements;
see~\eqref{eq:isomIu}. Roughly speaking, under this coupling,
the level set at level $-a$ contains every vertex whose
critical sign cluster is hit by an independent interlacement
process at level $u=a^2/2$.
Let $N^u$ denote the number of such vertices in $B_r$.
The count $N^u$ is more tractable and, through the isomorphism,
can be shown to have an asymptotically deterministic density, see Lemma~\ref{lem:Ncvg}.
The remaining task is to compare $M^{-a}_r(\psi^{2r})$ with
$N^u$, using a local uniqueness estimate for random interlacement to show that almost
all these vertices belong to a single cluster inside $\tB_r$.
The massive nature of our model introduces a further difficulty. Interlacement
trajectories may be killed before accumulating enough capacity
for the local uniqueness argument to apply.
We therefore establish in Proposition~\ref{lem:LU} a local uniqueness estimate for forward-surviving
trajectories.
Restriction to forward-surviving
trajectories, however, causes a loss that needs to be controlled.
Fortunately, the loss turns out to be of lower order. See Lemma~\ref{lem:h-isom} and~\eqref{eq:count_less_exp}.

\medskip
\noindent\textbf{Lower tightness.}
We finally discuss the lower tightness estimate in~\eqref{eq:tailEstlargestVolgena=0}, for which we need to strengthen~\eqref{eq:proofstep1} to a lower bound
close to one. It turns out that only under \eqref{def:Con} does the model have a sufficiently good decoupling property for such a bound to be true.
We use \eqref{eq:introDMP} to obtain
independent local fields in disjoint balls of radius
$\zeta\rho$ inside $B_r$. Each local field offers positive probability of producing a large cluster in the smaller ball of radius $\rho$.
The argument therefore requires a delicate balance between
the cost of decoupling and the independent contribution from local
fields. Condition~\eqref{def:Con} allows us to pick the scale $\zeta$ sufficiently large to decouple while
retaining enough disjoint balls inside $\tB_r$. After balancing these terms, the
contribution from having bad harmonic averages decays only
polynomially, in contrast to the exponential decay available
when $d\geq 3$. We refer to
Proposition~\ref{lem:ConLb} and the surrounding discussion for details.

\subsection{Applications: scaling relations at criticality and the volume tail}

{ We now discuss applications of Theorem~\ref{thm:antiCon} and the tools described above. Scaling theory predicts that critical observables of percolation models follow power laws, whose exponents satisfy a universal set of {scaling relations}. For the cable-graph Gaussian free field in dimensions $d\geq3$, these questions have been extensively studied; see \cite{lupu_2016_from,drewitz_2023_criti,drewitz_2022_clust,drewitz_2024_clust,drewitz_2024_criti,drewitz_2026_arm,drewitz_2026_fract,cai_2024_onear,cai_2024_oneara,cai_2025_quasi}. We show that, formulated directly in terms of the observables as in Kesten's seminal paper \cite{kesten_1987_scali}, the scaling relations persist on $\Slab$, even though the observables themselves are not power laws; see Table~\ref{t:critExp} and below.

The observables we compare are the one-arm probability, the two-point function, the volume tail, and the characteristic volume at scale $r$. We take the latter to be $\E_N[M_r]$ rather than a typical value of $M_r$, since by Theorem~\ref{thm:antiCon} the latter need not lie on a single scale. 
We denote by $\K$ the cluster of the origin in $\{\varphi\geq0\}$.
By the one-arm estimate of \cite[Theorem~1.1]{rodriguez_2025_on} (see \eqref{eq:one-arm}) and since $|B_r|\asymp r^2(r\wedge h_N)$, one has $\bv(r)\asymp |B_r|\,\P_N(\operatorname{rad}(\K)\geq r)$ in both regimes of \eqref{eq:thmM0rExp}, in line with the heuristic that the largest cluster in $B_r$ consists of the points connected to distance $r$. We start with the volume tail.}

\smallskip

\begin{theorem}\label{thm:volTailMain}
For all $N\geq 1$ and $1\leq n\leq \Cr{range}\bv(N)$,
\begin{equation}
\label{eq:volume_tail_main}
    \P_N\big(|\K|\ge n\big)
    \asymp
    \begin{cases}
        n^{-1/5}, 
        & \displaystyle n^{2/5}\leq \Dcross,
        \\[1.5ex]
        \displaystyle
        \Bigg(\frac{g_N\,h_N }{ \log^{\star}\!\Big( \frac{N h_N^{1/4}}  {n^{1/2}(g_N)^{1/4}} \Big) }\Bigg)^{-1/2},
        & \displaystyle n^{2/5}> \Dcross .
    \end{cases}
\end{equation}
In particular,
\begin{equation}
  \P_N\big(|\K|\ge \bv(N)\big)\asymp (\log(N)\vee h_N)^{-\frac12}.
\end{equation}
\end{theorem}
 In the case $h_N=N$, corresponding to three spatial dimensions, one has $\Dcross\asymp N$, so that throughout the relevant range we have $n^{2/5}\lesssim \Dcross$. Hence \eqref{eq:volume_tail_main} matches the result on $\Z^3$ obtained in \cite{drewitz_2024_clust}, and separately by \cite{cai_2025_quasi}, by which the decay is $n^{-(d-2)/(d+2)}$ for $3\leq d\leq5$. For $d>6$, the mean-field volume tail $\P(|\K|\geq n)\asymp n^{-1/2}$ was established in \cite{cai_2024_onear}.

 A somewhat surprising feature of { the one-arm probability $\P_N(\operatorname{rad}(\K)\geq r)\asymp \bv(r)/|B_r|$} is the
\emph{plateau} phenomenon observed in \cite[Remark~1.2,2)]{rodriguez_2025_on}.
In the delocalized regime, that is, when $g_N\to\infty$ or equivalently
$h_N=o(\log N)$, the one-arm probability no longer depends on the strength of
the delocalization. Indeed, $ g_N h_N\asymp \log N$
throughout this regime, so the dependence on $h_N$ cancels from { the one-arm probability}.
No analogous plateau occurs for the volume tail. For instance, taking $n=\frac{N^2}{\sqrt{\log N}}$,
\eqref{eq:volume_tail_main} and $g_N\asymp \log(N)/h_N$ yield, throughout
the delocalized regime,
\[
    \P_N\big(|\K|\geq n\big)
    \asymp
    \sqrt{{\log^\star(h_N)}/{\log N}}.
\]
Thus the volume tail continues to distinguish different rates of
delocalization. This is natural
in view of the volume growth
$|B_r|\asymp r^2(r\wedge h_N)$, through which the slab height enters
explicitly into the relation between cluster radius and volume.
%\item when do things enter the flat regime and 3d regime?

\begin{figure}[t]
\centering
\resizebox{1\textwidth}{!}{%
\begin{tikzpicture}[>=stealth, line cap=round, line join=round]

% ============================================================
% Left panel: one-arm plateau
% ============================================================
\begin{scope}

\def\xR{4.7}
\def\yP{2.15}

% axes
\draw[->,thick] (0,0) -- (5.25,0)
    node[right,font=\small] {$h_N$};
\draw[->,thick] (0,0) -- (0,3.45);

\node[below=4pt,font=\small] at (0,0) {$1$};
\draw (\xR,0.08) -- (\xR,-0.08)
    node[below=4pt,font=\small] {$\log N$};

% observable
\node[font=\small] at (1.55,3.25)
{$
\P_N\!\big(\rad(\K)\geq r\big)
$};

% plateau value
\node[left=5pt,font=\scriptsize] at (0,\yP)
{$\displaystyle
\sqrt{\tfrac{\log(N/r)}{\log N}}
$};

% vertical guide at h_N = log N
\draw[densely dotted,thick]
    (\xR,0) -- (\xR,\yP);

% plateau
\draw[blue,very thick]
    (0,\yP) -- (\xR,\yP);

\fill (0,\yP) circle (0.05);
\fill (\xR,\yP) circle (0.05);

\end{scope}

% ============================================================
% Right panel: volume tail
% ============================================================
\begin{scope}[shift={(8.3,0)}]

\def\xR{4.7}
\def\yL{1.15}
\def\yR{2.85}

% axes
\draw[->,thick] (0,0) -- (5.25,0)
    node[right,font=\small] {$h_N$};
\draw[->,thick] (0,0) -- (0,3.45);

\node[below=4pt,font=\small] at (0,0) {$1$};
\draw (\xR,0.08) -- (\xR,-0.08)
    node[below=4pt,font=\small] {$\log N$};

% observable
\node[font=\small] at (1.75,3.25)
{$
\P_N(
|\K|\geq \tfrac{N^2}{\sqrt{\log N}}
)
$};

% upper guide
\draw[densely dotted,thick]
    (0,\yR) -- (\xR,\yR);
\draw[densely dotted,thick]
    (\xR,0) -- (\xR,\yR);

% endpoint values
\node[left=5pt,font=\scriptsize] at (0,\yL)
{$\displaystyle
\tfrac{1}{\sqrt{\log N}}
$};

\node[left=5pt,font=\scriptsize] at (0,\yR)
{$\displaystyle
\sqrt{\tfrac{\log\log N}{\log N}}
$};

% increasing concave curve
\draw[blue,very thick]
    (0,\yL)
    .. controls (0.35,1.75) and (1.10,2.20) ..
    (2.25,2.45)
    .. controls (3.20,2.656) and (4.10,2.78) ..
    (\xR,\yR);

\fill (0,\yL) circle (0.05);
\fill (\xR,\yR) circle (0.05);

\end{scope}

\end{tikzpicture}%
}

\caption{Dependence on $h_N$ in the delocalised regime
$1\leq h_N\ll\log N$. In the left figure, the displayed one-arm estimate holds for every $1\leq r\leq N/2$ and is independent of $h_N$. In the right figure, we fix one example of a volume scale $n=N^2/\sqrt{\log N}$, for which the volume tail retains a non-trivial dependence on $h_N$.}
\label{fig:plateau-volume}
\end{figure}
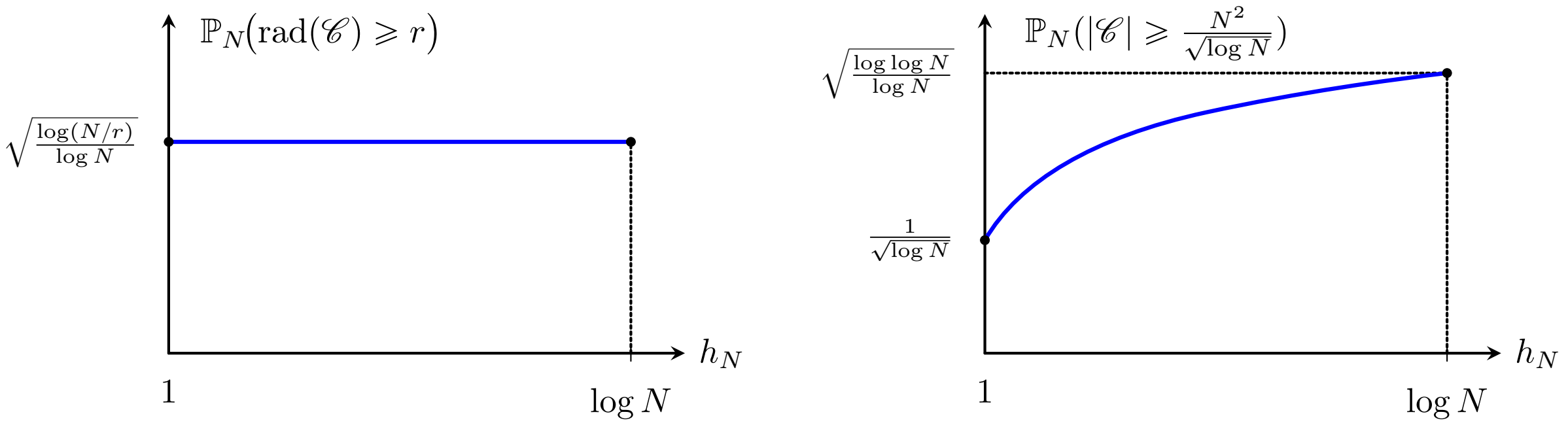

 We now state the aforementioned relations. The estimates below hold uniformly over the stated range of $r$ and $\|x\|$, which may in particular depend on $N$. In physicists' language, they constitute critical hyperscaling relations.

\begin{corollary}\label{thm:critScaling}
For all $N\geq 1$, $1\leq r\leq \Cr{range}N$ and $x\in \Slab$ with $1\leq \|x\|\leq \Cr{range}N$,
\begin{align}
  \label{eq:ptArm}
  \P_{N}\big(0\leftrightarrow x\big) &\asymp \P_{N}\big(\operatorname{rad}(\K)\geq \|x\|\big)^2\\
  \label{eq:vArm}
  \P_N\big(|\K|\geq \E_N[M_r] \big)&\asymp \P_{N}\big(\operatorname{rad}(\K)\geq r \big)\\
  \label{eq:ptvArm}
  \sum_{y\in B_r}\P_{N}\big(0\leftrightarrow y\big)
	 &\asymp 
  \E_N[M_r]\cdot  \P_{N}\big(\operatorname{rad}(\K)\geq r\big).  
\end{align}
\end{corollary}
See Remark~\ref{rem:classicalvgeneralscaling},1) for how Corollary~\ref{thm:critScaling} recovers the classical scaling relations when the relevant observables exhibit power law behaviour.

\begin{table}[t]
\centering
\small
\renewcommand{\arraystretch}{1.9}
\resizebox{\textwidth}{!}{%
\begin{tabular}{lll@{\hspace{0.5cm}}l}
\toprule
\noalign{\vskip -5pt}
Quantity & & Classical scaling & $\Slab$ with $h_N=N/\log N$\\
[-3pt]
\midrule
Volume tail
& $\P_N(|\K|\geq n)$
& $\approx n^{-1/\delta}$
& $\displaystyle
\asymp
\Big(
n^{2/5}\wedge
\frac{h_N}{\log^\star(Nh_N^{1/4}/n^{1/2})}
\Big)^{-1/2}$\\
Two-point function
& $\P_N(0\leftrightarrow x)$
& $\approx \|x\|^{-d+2-\eta}$
& $\displaystyle
\asymp
\Big(
\|x\|\wedge\frac{h_N}{\log(N/\|x\|)}
\Big)^{-1}$\\
One-arm probability
& $\P_N(\operatorname{rad}(\K)\geq r)$
& $\approx r^{-1/\rho}$
& $\displaystyle
\asymp
\Big(
r\wedge\frac{h_N}{\log(N/r)}
\Big)^{-1/2}$\\
Characteristic volume
& $\E_N[M_r]$
& $\approx r^{d_f}$
& $\displaystyle
\asymp
r^2(r\wedge h_N)
\Big(
r\wedge\frac{h_N}{\log(N/r)}
\Big)^{-1/2}$\\
\bottomrule
\end{tabular}%
}
\caption{Comparison between the classical power law scaling of the
critical observables and the corresponding behaviour on $\Slab$ with
$h_N=N/\log N$. The minima capture the crossover between the
three-dimensional and flat regimes.
The estimates hold uniformly for all $N\geq 1$,
$1\leq r\leq N/2$, $x\in\Slab$ with $1\leq\|x\|\leq N/2$,
and $1\leq n\leq c N^{5/2}/\sqrt{\log N}$.
We use the notation $a\wedge b=\min\{a,b\}$.}
\label{t:critExp}
\end{table}

To illustrate Corollary~\ref{thm:critScaling} when the underlying observables are not power laws, Table~\ref{t:critExp} displays the example $h_N=N/\log N$. Verifying \eqref{eq:ptArm}--\eqref{eq:ptvArm} by hand then requires identifying the relevant regime and tracking all logarithmic and iterated-logarithmic contributions; see Remark~\ref{rem:classicalvgeneralscaling},2).
This example also illustrates why it is natural here to formulate scaling relations directly in terms of the underlying observables, rather than through critical exponents. Such formulations go back to the seminal work of Kesten~\cite{kesten_1987_scali}; see also \cite{duminil-copin_2022_plana,duminil-copin_2026_near} for recent results in this direction.
In these works, however, the relevant observables are still expected to obey power laws, whereas this need not be the case for our model.

\subsection{\texorpdfstring{{ Organisation}}{Organisation}}
We now briefly describe how this article is organised. In Section~\ref{sec:facts}, we collect various preliminaries.
The main result in Section~\ref{sec:mainprop} is Proposition~\ref{pro:manyMeso}, which guarantees, with positive probability, the existence
of a supercritical cluster with large volume and capacity. As alluded to around~\eqref{eq:proofstep1},~\eqref{eq:thmM0rExp}, Theorem~\ref{thm:volTailMain} and Corollary~\ref{thm:critScaling} all follow from Proposition~\ref{pro:manyMeso}; we present their proofs in Section~\ref{sec:nonlim_thm}.
In Section~\ref{sec:matchedLemma}, we prove the simpler
matched-scale counterpart to Proposition~\ref{pro:manyMeso}, for which techniques available from higher dimensions readily apply. We defer the proof of Proposition~\ref{pro:manyMeso} to Section~\ref{sec:anticon}.

Section~\ref{sec:anticon} develops the comparison tool outlined in Section~\ref{sec:proofIdea}, which handles killing sensitivity.
Using the Gaussian structure, Proposition~\ref{prop:shift_h0}
compares local clusters of the full and killed fields through
a random sprinkling of the level. Its proof relies on a gradient
estimate for the Green's function, established at the end of the
section. In Section~\ref{sec:proof2}, we apply this comparison
to prove Proposition~\ref{pro:manyMeso}.
In a later section, Proposition~\ref{prop:shift_h0} will also be used crucially in the proof of Theorem~\ref{thm:limit} and to yield the probability
$\frac12$ in Theorem~\ref{thm:antiCon}.

In Section~\ref{sec:amplify}, we use the spatial Markov property
to strengthen the positive-probability estimates of
Section~\ref{sec:mainprop} to the lower tightness result in~\eqref{eq:tailEstlargestVolgena=0}.
The argument requires a delicate balance between
the cost of decoupling and the independent contribution from local
fields, which is only feasible under~\eqref{def:Con}; see Lemma~\ref{lem:loop_crossing_asymp}.

The remaining work concerns the~\eqref{def:AntiCon} regime.
In Section~\ref{sec:Iu}, we take the first step towards identifying
the limiting law by establishing precise volume asymptotics for
the largest killed cluster at a supercritical level, as stated in
Proposition~\ref{pro:convergences}.
The proof combines consequences of the isomorphism between the
field and random interlacements, collected in
Lemma~\ref{lem:isom_consequence}, with a local uniqueness estimate
for forward-surviving trajectories, established in
Proposition~\ref{lem:LU}.

The estimates in Sections~\ref{sec:anticon},~\ref{sec:amplify} and~\ref{sec:Iu} are brought together in
Section~\ref{sec:final}, where we conclude the proofs of
Theorems~\ref{thm:antiCon} and~\ref{thm:limit}.
Finally,
Appendices~\ref{sec:AFbox} and~\ref{A:grad} contain auxiliary estimates
for the capacity scales and the heat kernel respectively.

\medskip
In the sequel, we write $c,c^\prime,\dots$ for strictly positive constants which may change from line to line and write $c_i$ for strictly positive constants which will remain fixed. These constants have no dependence unless we explicitly say so.

\section{Setup and useful facts}\label{sec:facts}

We gather here a few preliminary results concerning cable graphs, random walks, potential theory and cluster capacity estimates, which are used throughout this article.

\subsection{The cable system $\tSlab$ and the diffusion $X$}
Recall the definition of $\Slab$ from \eqref{def:slab}.
We introduce non-negative weights $\lambda_{x,x'}=\lambda_{x',x}$ for $x,x'\in \mathbb{S}_N$ with $\lambda_{x,x'}= \tfrac16$ if and only if $\|x-x'\|=1$, and $\lambda_{x,x'}=0$ otherwise. We further introduce the killing measure $\kappa$ on $\mathbb{S}_N$ with $\kappa_x=  N^{-2}$ for all $x\in \mathbb{S}_N$. These choices imply that
$\lambda_{x}=\sum_{x^\prime} \lambda_{x,x^\prime}+\kappa_{x}=1 +N^{-2}$ for all $x\in \mathbb{S}_N$. We refer to the following weighted graph as a \emph{slab}:
$(\mathbb{S}_N,\lambda,\kappa).$ Its edges are the pairs $\{x,x'\}$ for which $\lambda_{x,x'}\neq 0$; they correspond to the natural product graph structure on $\mathbb{S}_N$.  

To obtain the cable system, $\tSlab^h \equiv \tSlab$, of the slab $\Slab$, we first replace all edges $\{x,x'\}$ by open intervals $I_{x_1,x_2}$ of length $(2\lambda_{x,x'})^{-1}$ and glue them through their endpoints. We then attach a half-open interval $I_x$ of length $(2\kappa_x)^{-1}$ to each $x\in \mathbb{S}_N$ (see \cite{drewitz_2024_criti,lupu_2016_from} for a more thorough background). We endow $\tSlab$ with the natural geodesic distance which assigns length 1 to each interval.
For $r\geq1$, we let
\begin{equation}\label{def:tBr}
  \tB_r\defeq B_r
  \cup \!\!\!\bigcup_{x,x'\in B_r, \lambda_{x,x'}>0}\!\!\! I_{x,x'}
  \cup \bigcup_{x\in B_r} I_x.
\end{equation}
Thus $\tB_r$ contains all cables with both endpoints in $B_r$, together with every killing cable attached to a vertex of $B_r$. We define $\tB(x,r)$ analogously, with $B(x,r)$ in place of $B_r$.
We write $\partial B_r$ for the set of vertices in $\Slab\setminus B_r$ adjacent to $B_r$ and set $\overline B_r=B_r\cup\partial B_r$.

The diffusion $X=(X_t)_{t \geq 0}$ is the continuous-time Markov chain on $\mathbb{S}_N\cup \{\Delta\}$, where $\Delta$ is an absorbing cemetery state, which jumps from $x$ to $y$ at rate $1$, with $\lambda_{x,\Delta}= \kappa_x$. 
We denote by $P_x$ and $E_x$ its law and expectation when started at $x\in \Slab$. $X$ extends naturally to a diffusion on $\tSlab$ that's killed when reaching the open end of a killing cable; see \cite[Section~2]{drewitz_2022_clust,prevost_2023_perco} for its precise definition.
With a slight abuse of notation, we also write $P_x$ and $E_x$ for the law and expectation of diffusion on $\tSlab$. 
For $y\in \tSlab$,
writing $\ell_y(t)$ for the local time of $X$ at $y$ up to time $t$, with respect to Lebesgue measure on the cables and stopped upon getting killed, we define
\begin{equation}\label{def:green}
  g_N^h(x,y)\defeq E_x[\ell_y(\infty)],
  \qquad x,y\in\tSlab.
\end{equation}
On $\Slab\times\Slab$, this agrees with the Green's function of the walk on $\Slab$. We usually suppress the height profile $h$ and write $g_N=g_N^h$. The Green function is symmetric and, on $\Slab$, satisfies $g_N(x,y)=g_N(0,y-x)$ by translation invariance.

For a nonempty compact set $K\subset\tSlab$, let
$H_K=\inf\{t\geq0: X_t\in K\}$. Its equilibrium measure $e_{K,N}$
is characterised by the last-exit formula
\begin{equation}\label{eq:last_exit}
  P_x(H_K<\infty)=\int_K g_N(x,y)\,e_{K,N}(dy),
  \qquad x\in\tSlab.
\end{equation}
The total mass of $e_{K,N}$ is the capacity $\capacity(K)$, which also satisfies
\begin{equation}\label{eq:CapVarForm}
  \capacity(K)=
  \Big(\inf_\mu\iint
    g_N(x,y)\,\mu(dx)\,\mu(dy)\Big)^{-1},
\end{equation}
where the infimum is over probability measures on $K$.
See \cite[Section~2]{drewitz_2022_clust} for more background.

\smallskip
We now recall some useful estimates on $\Slab$ which will be used repeatedly in the sequel. The following is a 
restatement of \cite[Theorem 2.1 and Proposition 3.2]{rodriguez_2025_on}.
For $1\leq r\leq \Cr{range}N$, let
\begin{equation}\label{def:Fbox}
  F(r)\defeq r \wedge  \frac{h_N}{\log^\star(N/r)}.
\end{equation}
By \cite[(3.9)]{rodriguez_2025_on}, one has that for all $0 \leq r \leq \Cr{range}N$ and $x\in \Slab$ with $\|x\|=r$,
\begin{align}
  \label{eq:capB}
  c F(r)\leq &\capacity(B_r) \leq c^\prime F(r),\\
\label{eq:g2-asymp-unif}
cF(r)^{-1} \leq &g_N(0,x)\leq c'F(r)^{-1} .
\end{align}
Note that combining \eqref{eq:g2-asymp-unif} and \cite[Proposition 5.2]{lupu_2016_from}, we immediately obtain the two-point function estimate for all $x,y\in \Slab$ such that $1\leq \|x\|,\|y\|\leq \Cr{range}N$.
\begin{equation}\label{eq:two-ptcrit}
	\P_N(x\leftrightarrow y)\asymp {g_N(x,y)}/{g_N(0,0)}\asymp {F(1)}/{F(\|x-y\|)}.
\end{equation}
Using \eqref{def:Dcross} and \eqref{def:Fbox}, \cite[Theorem~1.1]{rodriguez_2025_on} can equivalently be written, for all $1\leq r\leq \Cr{range}N$, as
\begin{equation}\label{eq:one-arm}
    \P_N\big(\operatorname{rad}(\K)\geq r\big)
    \asymp
    \begin{cases}
        r^{-1/2},
        & r\leq \Dcross,
        \\[1.5ex]
        \displaystyle
        \Big(
        \frac{g_N h_N}{\log^\star(N/r)}
        \Big)^{-1/2},
        & r>\Dcross.
    \end{cases}
\end{equation}

The following estimate on exit time will be useful in the sequel,
\begin{lemma}
For all $1\leq r\leq \Cr{range}N$  
\begin{equation}\label{eq:exitProb}
  P_0(H_{\partial B_r}<\infty)\geq c.
\end{equation}
\end{lemma}
\begin{proof}
   Let $\tau_N$ be an independent exponential random variable with rate $N^{-2}$. First note that using \cite[(2.51)]{lawler_2010_rando} and the tail of an exponential random variable, we have that
\begin{equation}\label{eq:hitting_lb}
	P_0(H_{\partial B_{r}}<\infty)\geq \tilde P_0(H_{\partial B_{r}}<r^2)P(\tau_N\geq r^2)\geq ce^{-\frac{r^2}{N^2}},
\end{equation}
where $\tilde P_0$ is the law of an unkilled simple random walk on $\Slab$ starting at $0$.
\end{proof}

\subsection{Killing set and Markov property}
We now allow additional killing on a set $U\subset\tSlab$.
We write $P_x^U$ and $E_x^U$ for the law and expectation of the diffusion started at $x\in\tSlab$ and
killed at $H_U=\inf\{t\geq0: X_t\in U\}$. The corresponding Green
function is
\begin{equation}\label{eq:gNK}
  g_N^U(x,y)\defeq E_x[\ell_y(H_U)],
  \qquad x,y\in\tSlab,
\end{equation}
so that $g_N^\emptyset=g_N$. We write $g_N^U$ for $g_N^U(0,0)$ whenever
$0\notin U$. For compact $K$ in the domain of the killed diffusion,
we denote by $e_{K,N}^U$ its equilibrium measure and write
$\capacity^U(K)=e_{K,N}^U(K)$ for its capacity.

For $r\geq1$, we abbreviate $g_N^r(x,y)=g_N^{\partial B_r}(x,y)$,
$g_N^r=g_N^r(0,0)$ and $\capacity^r(K)=\capacity^{\partial B_r}(K)$.
By applying the strong Markov property at $H_U$, one has
\begin{equation}\label{eq:gDecomp}
  g_N(x,y)=g_N^U(x,y)
  +E_x[H_U<\infty, g_N(X_{H_U},y)].
\end{equation}

\smallskip
We now restate \cite[(5.2) and (5.4)]{rodriguez_2025_on}.
For $1\leq r'\leq r\leq \Cr{range}N$, let
\begin{equation}\label{eq:defn_fbox}
  F^{r}(r')=r' \wedge \frac{h_N}{\log^{\star}\!\big(\frac{r}{r'\vee h_N}\big)}.
\end{equation}
   For all $N\geq 1$ and $2 \leq r \leq\Cr{range} N$, $x_0 \in \Slab$  the following hold:
   \begin{enumerate}
   \item[(i)] if $U= \Slab \setminus B(x_0,2r)$, one has that for all $x_1,x_2,\in B(x_0,r)$,
   \begin{align}\label{eq:killed_green_smallh} 
	\begin{split}
		  cF^{r}(\|x_1-x_2\|\vee1)^{-1}    
		 &\leq g_N^U(x_1,x_2) 
            \leq   c'F^{r}(\|x_1-x_2\|\vee1)^{-1}   ;
	\end{split}
   \end{align}
 \item[(ii)] if $1\leq r'\leq r\leq N$ and $x\in \Slab$ is such that $B(x,r')\subseteq B(x_0,r)$ and $U= \Slab \setminus B(x_0,2r)$,
      \begin{equation} \label{eq:cap-killed}
        \Cl[c]{killedBall_lb}F^r(r')
        \leq 
        \capacity^U\big(B(x,r')\big) \leq 
        {\Cl[c]{killedBall_ub}}F^r(r').        
      \end{equation}
      \end{enumerate}
 Note that \eqref{eq:cap-killed} is slightly more general than \cite[Corollary]{rodriguez_2025_on}, since the balls don't need to be centered at $x_0$. However one can obtain \eqref{eq:cap-killed} with the same proof.

Using the same proof as for \cite[Proposition 5.4]{rodriguez_2025_on}, one obtains the following control on the killed capacity of the range of random walk before exiting a ball.
For $T \geq 0$ (possibly random), we write $X_{[0,T]}= \{ x \in \mathbb{S}_N: x = X_s \text{ for some } 0 \leq s \leq T\}$ for the range of $X$ until time $T$.
\begin{proposition} \label{lem:cap_walk_space_2}
    There exists $\Cl[c]{rw_cap}>0$ such that for all $x\in \mathbb{S}_N$, $1\leq r'\leq r\leq \Cr{range}N$ such that $B(x,r')\subset B_r$, $2\leq s\leq \frac{r'}{2}$, we have that, with $T=H_{B(x,r')^{\text{c}}}$,
    \begin{equation}\label{eq:rw_cap_with_killing}
        P_x\left( \capacity^{2r}\big(X_{[0,T]}\big)
        \leq  \frac{\Cr{rw_cap}}{s}\capacity^{2r}\big(B(x,r')\big) ,T<\infty \right)
        \leq Ce^{-cs}.
    \end{equation}
\end{proposition}

For $U$ as before, we write $\P_N^U$ for the law of the centered Gaussian free field $\varphi$ on $\tSlab$ with covariance $g_N^U(\cdot,\cdot)$, so $\varphi_x=0$ for all $x \in U$.
We write $\P_N=\P_N^\emptyset$ and $\P_N^r=\P_N^{\partial B_r}$. For $K\subset\tSlab$, we introduce the decomposition
\begin{equation}\label{eq:field_decomp}
\eta^K_x=E_x^U\big[\varphi_{X_{H_K}}\1\{H_K<\infty\}\big] \quad\text{and}\quad \psi_x^K=\varphi_x-\eta^K_x, \quad x \in \tSlab.
\end{equation}
The strong Markov property asserts that if $\varphi$ is a sample of $\P_N^U$, then the fields $\eta^K_{\cdot}$ and $\psi_{\cdot}^K$ are independent centered Gaussian fields under $\P_N^U$, and $\psi_{\cdot}^K$ has law $\P_N^{U\cup K}$.

\subsection{Cluster capacity}\label{subsec:clustCap}
Let us first introduce, for any continuous function $\phi:\tSlab\to\R$, $a\in\R$, $x\in \tSlab$ and closed set $\tB\subset\tSlab$, the local cluster of $x$ in $\{\phi \geq a \}\cap \tB$, defined as
\begin{equation}
\label{eq:Kr}
{\K}^{a}_{\tB}(x,\phi) =\big\{ y \in  \tB: y \leftrightarrow x \text{ in } \{\phi\geq a \}\cap \tB \big\}.
%\qquad \nn{{\K}^{a}_{z,r}(x,\phi)\defeq{\K}^{a}_{\tB(z,r)}(x,\phi),\quad {\K}^{a}_{r}(x,\phi)\defeq{\K}^{a}_{0,r}(x,\phi).}
\end{equation}
If $\tB=\tB(x,r)$ is a { ball}, we will just use $x,r$ as subscripts, and omit $x$ altogether if $x=0$.
Note that ${\K}^{a}_{\tB}(x,\phi)$ is empty unless $x\in \tB$. When $\phi=\varphi$, we abbreviate
${\K}^{a}_{\tB}(x)={\K}^{a}_{\tB}(x,\varphi)$, and similarly ${\K}^{a}_{x,\nn{r}}(x)$, ${\K}^{a}_{\nn{r}}(x)$; we also drop all superscripts when $a=0$.
We write $\K^a(x)=\K^a_{\tSlab}(x)$ for the cluster of $x$ in
$\{\varphi\geq a\}$ on the entire cable system.

\smallskip

For any weighted graph $\G$ and its associated metric graph $\tG$ where $(\varphi_{x})_{x\in \tilde\G}$ is well-defined, we write $\P_{\tG}$ for the law of the field on $\tG$ and write $g_{\tG}(\cdot,\cdot)$ and $\capacity_{\tG}(\cdot)$ for the Green's function and capacity with respect to the diffusion on ${\tG}$. 
It  follows from \cite[Corollary~3.3,(1), Lemma~3,4,(2) and Theorem~3.7]{drewitz_2022_clust} that if $\sup_{x\in{G}}g_{{{\tG}}}(x,x)<\infty$, the density of $\capacity_{\tG}( \K^{a}(x))$ under $\P^{\tG}$ is given by \cite[(3.8)]{drewitz_2022_clust} and that for any $s>g_{\tG}(x,x)^{-1}$, one has the tail estimate
\begin{equation}
\label{eq:capTail}
    c \big({g_{\tG}(x,x)}s\big)^{-\frac12}e^{-c'a^2s} \leq \P_{\tG}\!\big(s \leq \capacity_{\tG}( \K^{a}(x)) < \infty 
    \big)
    \leq c' \big({g_{\tG}(x,x)}s\big)^{-\frac12} e^{-c'a^2s}.
\end{equation}
We abbreviate $\P^{r}=\P_{\tilde{\Z}^d\setminus \tB_r^c}$ and with a light abuse of notation, let $\capacity^r=\capacity_{\tilde{\Z}^d\setminus \tB_r^c}$.
The following Lemma is essentially \cite[Lemma 3.1]{drewitz_2024_clust}, which is itself a simple consequence of \eqref{eq:capTail}.
\begin{lemma}
\label{lem:clustCapDPR}  
For all $d\in\{3,4,5\}$, there exists $\Cl[c]{DPR}(d)>0$ such that for all $1\leq r' \leq r/4$, $y\in B_{r/4}$ and $0<s\leq\Cr{DPR}(r')^{d-2}$,
\begin{equation} \label{eq:clustCapDPR}   
  c(d)\,(1\vee s)^{-1/2}\leq \P^{r}\!\Big(\capacity^{r}\!\big(\K^{0}_{\nn{y,r'}}(y)\big)\geq s\Big)\leq c'(d)\,(1\vee s)^{-1/2}.
\end{equation}
\end{lemma}
Indeed, the upper bound simply follows from \eqref{eq:capTail} upon replacing $\K^0_{\nn{y,r'}}(y)$ by $\K^0(y)$, and the lower bound follows by the argument of \cite[Lemma 3.1]{drewitz_2024_clust} with $a=0$.
The same argument, using \eqref{eq:capTail} for the field with additional killing on $\partial B_r$, yields the following estimate on $\tSlab$.
\begin{lemma}
\label{lem:clustCapNoKill}  
There exists $\Cl[c]{C:lawcaplocNK}>0$ such that for all $1\leq r'\leq r/4\leq \Cr{range}N$, $y\in B_{r/4}$ and $0<s\leq \Cr{C:lawcaplocNK}\capacity^r(B_{r'})$,
\begin{equation} \label{eq:clustCapNoKill}   
  c\big( (g_N^r\, s)\vee1 \big)^{- 1/2}\leq \P_N^r\Big(\capacity^{r}\!\big(\K^{0}_{\nn{y,r'}}(y)\big)\geq s\Big)\leq c'\big( (g_N^r\, s)\vee1\big)^{-1/2}.
\end{equation}
\end{lemma}
\begin{proof}
  Note that the upper bound follows from \eqref{eq:capTail} upon replacing $\K^0_{\nn{y,r'}}(y)$ by $\K^0(y)$. For the lower bound, we first claim that, for all $x\in B_r$ such that $B(x,r')\subseteq B_{r/2}$,
  \begin{equation}\label{eq:onearm_killed}
    \P^{r}_N(x\stackrel{\geq 0}{\leftrightarrow}\partial B(x,r')) \leq c\big(g^r_N F^r(r')\big)^{-1/2}.
  \end{equation}
The argument proving the upper bound in \cite[(1.9)]{rodriguez_2025_on} readily adapts to yield \eqref{eq:onearm_killed}, with $N=r$ and $R=r'$ therein. In particular, the Harnack inequality required in the proof of \cite[Lemma 6.4]{rodriguez_2025_on} remains available by \cite[Corollary 5.2]{rodriguez_2025_on}; the generalized Lupu formula used in the proof of \cite[Lemma 6.6]{rodriguez_2025_on} is provided by \cite[Proposition 2.1]{drewitz_2024_clust}; and the necessary Green's function and capacity estimates follow from \eqref{eq:cap-killed} and \eqref{eq:killed_green_smallh}. 
Hence, combining \eqref{eq:capTail}, \eqref{eq:cap-killed} and \eqref{eq:onearm_killed}, and taking $s\leq \Cr{C:lawcaplocNK}\capacity^r(B_{r'})$, where $\Cr{C:lawcaplocNK}$ is a small enough constant (cf.~\cite[(3.11)]{drewitz_2024_clust} for a similar argument), we have that,
\begin{align*}
\P_N^r(\capacity^r\!(\K^{0}_{\nn{y,r'}}(y))\geq s)&\geq 
\P_N^{r}\big(\capacity^r\!(\K^0(y))\geq s\big)-\P^{r}_N(y\stackrel{\geq 0}{\leftrightarrow}\partial B(y,r'))\\
&\geq c\big((g^r_N s)\vee 1\big)^{-\frac12}-c\big(g^r_N F^r(r')\big)^{-\frac12}\geq c\big((g^r_N s)\vee 1\big)^{-\frac12} .
\end{align*}
\end{proof}

\section{Volume of supercritical cluster}\label{sec:mainprop}
In this section, we state Proposition~\ref{pro:manyMeso}, which
guarantees, with positive probability, the existence
of a supercritical cluster with large volume and capacity.
As in~\cite{drewitz_2024_clust}, this estimate provides the key
input for proving~\eqref{eq:thmM0rExp},
Theorem~\ref{thm:volTailMain} and Corollary~\ref{thm:critScaling}.
We give the proofs of these results using Proposition~\ref{pro:manyMeso} in Section~\ref{sec:nonlim_thm}.

We will not prove Proposition~\ref{pro:manyMeso} in this section. Instead, we focus on a simpler statement, Lemma~\ref{pro:manyMesoKilled}, which is identical to Proposition~\ref{pro:manyMeso}, except that the ambient field is sampled under a killed measure. The proof of Lemma~\ref{pro:manyMesoKilled} is presented in Section~\ref{sec:matchedLemma}.
We will pass from 
Lemma \ref{pro:manyMesoKilled} to Proposition~\ref{pro:manyMeso} in Section~\ref{sec:anticon} via the comparison tool we introduced in Section~\ref{sec:proofIdea}.

Recall the notation from~\eqref{eq:Kr}.
\begin{proposition}\label{pro:manyMeso} 
There exist constants  $\Cl[c]{c:LBhatK},\Cl[c]{MarCap}\in(0,1)$ and $\Cl[c]{proprLB},\Cl[c]{C:LBLneg}<\infty$ such that for all $a\in{(0,\sqrt{g_N}]}$, $\Cr{proprLB}\leq r\leq \Cr{range}N$ with $a^2 \Fbox(r)\geq \Cr{C:LBLneg}$, the following hold:
\begin{equation} \label{eq:LBL2} 
        \P_N \bigg(\exists\, y\in B_r:\, \big| {\K}^{-a}_{\nn{r}}(y)\big| \ge \Cr{c:LBhatK}\frac{ar^2(r\wedge h_N)}{\sqrt{g_N}},
         \, \capacity^{r}\big({\K}^{-a}_{\nn{r}}(y)\big)\geq \Cr{MarCap} (r\wedge h_N) \bigg) \ge 
    c.
\end{equation}
\end{proposition}
\begin{comment}
\end{comment}

We now turn to the killed counterpart of Proposition~\ref{pro:manyMeso}, which we will prove in this section.
\begin{lemma}\label{pro:manyMesoKilled} 
For all $2\leq\lambda \leq \Cr{range}$, there exist constants  $\Cl[c]{c:LBhatK2},\Cl[c]{MarCap2}\in(0,1)$ and $\Cl[c]{proprLB2},\Cl[c]{C:LBLneg2}<\infty$ such that for all $a\in{(0,\sqrt{g_N^{\lambda r}}]}$, $\Cr{proprLB2}\leq r\leq \Cr{range}N$ with $a^2 F^r(r)\geq \Cr{C:LBLneg2}$, the following hold:
\begin{equation} \label{eq:LBL2Killed} 
        \P^{\lambda r}_N \bigg(\exists\, y\in B_r:\, \big| {\K}^{-a}_{\nn{r}}(y)\big| \ge \Cr{c:LBhatK2}\frac{ar^2(r\wedge h_N)}{\sqrt{g^{\lambda r}_N}},
         \, \capacity^{r}\big({\K}^{-a}_{\nn{r}}(y)\big)\geq \Cr{MarCap} (r\wedge h_N) \bigg) \ge 
    c.
\end{equation}
\end{lemma}
Crucially, for Lemma~\ref{pro:manyMesoKilled}, the scale at which the field is killed matches the scale of the cluster under consideration: both are of order $r$. In this matched-scale setting, the proof follows the argument of \cite[Section~4]{drewitz_2024_clust}, up to some minor modifications that we explain in Section~\ref{sec:matchedLemma}. We refer to Remark~\ref{rem:why_difficult} for more on this.

The following result is a simple consequence of Proposition~\ref{pro:manyMeso}, Lemma~\ref{pro:manyMesoKilled}
and the change of measure formula~\cite[(3.7)]{drewitz_2024_clust}. 
Recall the notation from~\eqref{def:Mr_general}.
\begin{corollary}\label{the:mainMra}
	For all $1\leq r\leq \Cr{range}N$, $2\leq\lambda\leq \Cr{range}$, let $K\in \{\partial B_{\lambda r},\emptyset\}$ and $|a|\leq \sqrt{g_N^K}$. Then
	\begin{equation}\label{eq:tailEstlargestVolgen}
		\P_N^{K}\bigg(M^{a}_r\geq c\frac{r^2(r\wedge h_N)}{\sqrt{g_N^K\,(a^{-2} \wedge \capacity^K(r)) }}\bigg)\geq c\exp\big\{\!-\!c'(a\vee 0)^2 \capacity^K(r)\big\}.
	\end{equation}
\end{corollary}
We omit the proof of Corollary~\ref{the:mainMra} and refer to~\cite[(2.19)]{drewitz_2024_clust} for details. We now explain how to conclude the proofs of \eqref{eq:thmM0rExp}, Theorem~\ref{thm:volTailMain} and Corollary~\ref{thm:critScaling} using the above estimates.

\subsection{Proof of \eqref{eq:thmM0rExp}, Theorem~\ref{thm:volTailMain} and Corollary~\ref{thm:critScaling}}\label{sec:nonlim_thm}
We start with Lemma~\ref{pro:upperMRa}, a general bound on the mean of $M^0_r$ where the law may involve a killing set. 
\begin{lemma}
\label{pro:upperMRa}
For all $1\leq r\leq \Cr{range}N$ and $2\leq \lambda \leq \Cr{range}$, let $K\in \{\partial B_{\lambda r}, \emptyset\}$. Then $\E_N^{K}[\LV_r^2]/\E_N^{K}[\LV_r]^2\leq c$ and
\begin{equation}
\label{eq:upperMR0}
\begin{split} 
\E_N^{K}\big[\LV_r\big]  \asymp \frac{|B_r|}{\sqrt{g_N^{K}\capacity^K\!(B_r)}}.
\end{split}
\end{equation}
\end{lemma}
Since $|B_r|\asymp r^2(r\wedge h_N)$,
by~\eqref{def:Fbox} and the definition of $\bv$ from~\eqref{eq:thmM0rExp}, it's easy to see that $\bv(r)\asymp |B_r|(g_N\,F(r))^{-\frac12}$, hence~\eqref{eq:thmM0rExp} is an immediate consequence of~\eqref{eq:capB} and Lemma~\ref{pro:upperMRa} in the case $K=\emptyset$.
\begin{proof}
The lower bound in~\eqref{eq:upperMR0} is a consequence of \eqref{eq:tailEstlargestVolgen} with $a=0$ therein.
For the upper bound, proceeding as in the proof of \cite[Proposition 2.4]{drewitz_2024_clust}, we have that
\begin{equation}\label{eq:MrTrick}
  M_r^2\leq \sum_{x\in B_r}|\K^0_r(x)|
  \quad\text{and}\quad
  \E_N^{K}\big[\LV_r\big]\leq \E_N^{K}\Big[\sum_{x\in{B_r}}|\K_{\nn{r}}^0(x)|\Big]^{\frac12}\leq c|B_r|^{\frac12}\sup_{x\in B_r}\E^{K}_N\big[|\K^0_{\nn{x,2r}}(x)|\big]^{\frac12}.
\end{equation}
By \eqref{eq:two-ptcrit}, \eqref{eq:killed_green_smallh} and a calculation similar to \cite[(3.26) and (3.27)]{rodriguez_2025_on}, we have for all $x\in B_r$ that,
\begin{equation}\label{eq:volExp}
\E_N^{K}\big[|\K^0_{\nn{x,2r}}(x)|\big]\leq c\sum_{y\in{B(x,2r)}}\frac{g_N^{K}(x,y)}{g_N^{K}}\leq c\frac{|B_r|}{g_N^{K}\capacity^K\!(B_r)}.
\end{equation}
Combining the above two displays yields $\E_N^{K}[\LV_r^2]/\E_N^{K}[\LV_r]^2\leq c$ and the upper bound in \eqref{eq:upperMR0}.
\end{proof}

\begin{proof}[Proof of { Corollary} \ref{thm:critScaling}]
The relation \eqref{eq:ptArm} follows directly from \eqref{eq:two-ptcrit} and \eqref{eq:one-arm}. The upper bound in \eqref{eq:ptvArm} follows from a calculation similar to \eqref{eq:volExp}, combined with \eqref{eq:two-ptcrit}, \eqref{eq:one-arm} and \eqref{eq:thmM0rExp}. Now, note that by \eqref{eq:two-ptcrit} we have $\P_N(0\leftrightarrow y)\geq (g_N \Fbox(r))^{-1}$ for all $y\in B_r$. Combining this with \eqref{eq:thmM0rExp}, \eqref{eq:one-arm} and the fact that $|B_r|\geq cr^2(r\wedge h_N)$ yields the lower bound in \eqref{eq:ptvArm}.

It remains to show \eqref{eq:vArm}. We follow the proof of \cite[Corollary 1.6]{drewitz_2024_criti}. By \eqref{eq:one-arm}, Markov's inequality and \eqref{eq:volExp} we have that for $r\geq 1$,
\begin{align*}
	\P_N(|\K^0|\geq \bv(r))
	&\leq \P_N(\rad(\K^0)\geq r)+\frac{1}{\bv(r)}\E_N[|\K^0\cap B_{r}|]\\
	&\leq c\big(g_N\,F(r)\big)^{-\frac12}\!+\!\frac{cr^2(r\wedge h_N)}{\bv(r)g_N\Fbox(r)}\!\leq\! c'\big(g_N\,F(r)\big)^{-\frac12}\!\leq\! c\,\P_N(\rad(\K^0)\geq r).
\end{align*}
For the lower bound, note that by \eqref{def:df}, we have that $\bv(r)\1\{M_r\geq \bv(r)\}\leq \sum_{x\in B_r}\1\{|\K(x)|\geq \bv(r)\}$.
It then follows from \eqref{eq:tailEstlargestVolgen}, \eqref{eq:one-arm} and translation invariance that
\[
\P_N(|\K|\geq \bv(r))\geq \frac{c\,\bv(r)}{|B_r|}\geq c\big(g_N\,F(r)\big)^{-\frac12}\geq c\,\P_N(\rad(\K^0)\geq r).
\]

\end{proof}

It remains to show Theorem \ref{thm:volTailMain}, the proof of which is purely computational.

\begin{proof}[Proof of Theorem \ref{thm:volTailMain}]
It suffices to combine \eqref{eq:vArm} with \eqref{eq:one-arm} and identify the scale $r$ such that $\bv(r)\asymp n$.
Assume first that $n^{2/5}\leq\Dcross$. In this case, $g_N\asymp1$ and $F(n^{2/5})=n^{2/5}$, hence $\bv(n^{2/5})\asymp n$. It then follows from \eqref{eq:one-arm} that,
$
\P_N(|\K^0|\geq n)\asymp n^{-1/5},
$
which concludes the proof of \eqref{eq:volume_tail_main} in this case.

Now assume $n$ is such that $n^{2/5}>\Dcross$ and $n^{1/3}
(\frac{g_Nh_N}{\log^\star(N/h_N)})^{1/6}
\leq h_N$. Let us pick
\[
r= n^{1/3}
\Big(\frac{g_Nh_N}{\log^\star(N/h_N)}\Big)^{1/6}
(\leq h_N).
\]
Then $r\gtrsim\Dcross$ and $F(r)\asymp h_N/\log^\star(N/h_N)$, whence $\bv(r)\asymp n$ and,
\[
\frac{N h_N^{1/4}}{n^{1/2}g_N^{1/4}}
\stackrel{n\asymp \bv(r)}{\asymp} \frac{N h_N^{1/4}g_N^{1/4}\Fbox(r)^{1/4}}{r^{3/2}g_N^{1/4}}  \asymp
\frac Nr\Big(\frac{h_N}{r}\Big)^{1/2}
\log^\star(N/h_N)^{-1/4}.
\]
Now using $h_N^\star\!\lesssim \!r\!\leq \!h_N$, which implies $\log^\star(N/r)\!\asymp\! \log^\star(N/h_N)$ and $c\leq (h_N/r)^{1/2}\!\leq \!\log^\star(\frac N{h_N})^{1/2}$, we have that $\log^\star\!(
\frac{N h_N^{1/4}}{n^{1/2}g_N^{1/4}}
)
\!\asymp\!
\log^\star(\frac Nr).$
Thus \eqref{eq:one-arm} gives the second line of
\eqref{eq:volume_tail_main}.
Finally, $n^{2/5}>\Dcross$ and $n^{1/3}
(\frac{g_Nh_N}{\log^\star(N/h_N)})^{1/6}
> h_N$.
Set $x=n^{1/2}g_N^{1/4}h_N^{-1/4}$ and $r=\frac{x}{\log^\star(N/x)^{1/4}}.$
The fact that $n^{1/3}
(\frac{g_Nh_N}{\log^\star(N/h_N)})^{1/6}
> h_N$ gives $x\gtrsim h_N\log^\star(N/h_N)^{1/4}$, and hence $r\gtrsim h_N$. It is then easy to see that $F(r)\asymp
\frac{h_N}{\log^\star(N/x)}$ and $
\bv(r)\asymp n.$
Applying \eqref{eq:one-arm} once more yields the second line of
\eqref{eq:volume_tail_main}.
\end{proof}

\Needspace{4\baselineskip}
\begin{remark}\label{rem:classicalvgeneralscaling}\leavevmode
\begin{enumerate}[label={\arabic*)}]
 \item  To see the relationship between Corollary~\ref{thm:critScaling} and the conventional language of critical exponents, recall that classical scaling theory predicts that the near-critical behaviour of the model is encoded by a collection of \emph{critical exponents} $\delta,\eta,\rho,d_f,\beta,\gamma,\nu$, introduced in Table~\ref{t:critExp}. We refer the reader to \cite[Chapter 9]{grimmett_1999_perco} for a more detailed discussion on the exponents $\delta,\eta,\rho,\beta,\gamma,\nu$.
Classical scaling theory predicts the relations
\begin{equation}\label{eq:classicalScaling}
  \Delta=\delta\beta ,
  \qquad
  2-\alpha=\gamma+2\beta=\beta(\delta+1),
  \qquad
  \gamma=\nu(2-\eta),
\end{equation}
and, at and below the upper critical dimension, the hyperscaling relations
\begin{equation}\label{eq:classicalHyperscaling}
    d\nu=2-\alpha,
    \qquad
    d\rho=\delta+1, 
    \qquad
    d_f=d/(1+1/\delta),
  \end{equation}
see \cite[Chapter~9]{grimmett_1999_perco} and \cite[(42a)]{stauffer_1979_scali}. 
Of particular relevance to us are the consequences of
\eqref{eq:classicalScaling} and \eqref{eq:classicalHyperscaling} involving
only observables at criticality. The hyperscaling relations give $\rho=(\delta+1)/d$ and
$d_f=d\delta/(\delta+1)$, while, using
$\gamma=\beta(\delta-1)$, $\gamma=\nu(2-\eta)$ and
$d\nu=2-\alpha=\beta(\delta+1)$, one further obtains
$2-\eta=d(\delta-1)/(\delta+1)$.
These identities in turn yield
\begin{equation}\label{eq:criticalHyperscaling}
    d-2+\eta=2/\rho,
    \qquad
    \delta=\rho\times d_f,
    \qquad
     2-\eta=d_f-1/\rho.
\end{equation}
Thus, according to the hyperscaling picture,
fixing any one of $\delta,\eta,\rho$ and $d_f$ determines the other three exponents associated with critical observables. Corollary~\ref{thm:critScaling} gives a generalized form of \eqref{eq:criticalHyperscaling}; under the additional assumption that the relevant observables scale as power laws with the corresponding critical exponents, it reduces to \eqref{eq:criticalHyperscaling}.
\item 
We verify Corollary~\ref{thm:critScaling} by hand in the example $h_N=N/\log N$ of Table~\ref{t:critExp}.
In this case, $g_N\asymp1$ and $\Dcross\asymp \frac{N}{\log N\log\log N}$. As displayed in Table~\ref{t:critExp},
the estimates can be written as the minimum between two distinct scales. For each observable, this minimum is attained by the term involving $h_N$ in the flat regime and by the power law term in the three-dimensional regime.

The relation \eqref{eq:ptArm} is straightforward.
For \eqref{eq:vArm}, suppose first that $r\leq\Dcross$.
Then $\E_N[M_r]\asymp r^{5/2}$, and the first term in the minimum
in the volume-tail estimate gives
\[
\P_N\big(|\K|\geq\E_N[M_r]\big)
\asymp r^{-1/2}\asymp \P_N\big(\operatorname{rad}(\K)\geq r\big).
\]
Next, suppose that $\Dcross\leq r\leq h_N$. In this regime,
$\log(N/r)\asymp\log\log N$, so $\E_N[M_r]\asymp r^3\sqrt{{\log\log N}/{h_N}}.$
Substituting this scale into the logarithm in the volume tail estimate
yields
\[
\log^\star\!\Big(
\frac{N h_N^{1/4}}
{r^{3/2}(\log\log N/h_N)^{1/4}}
\Big)
=
\Big(
\log\frac Nr+\frac12\log\frac{h_N}{r}
-\frac14\log\log\log N
\Big)\vee1
\asymp\log\log N,
\]
Consequently, since $r\geq \Dcross \asymp h_N/\log\log(N)$,
\[
\P_N\big(|\K|\geq\E_N[M_r]\big)
\asymp \Big(
r^{\frac65}\Big(\frac{\log\log N}{h_N}\Big)^{\frac15}
\!\wedge\!\frac{h_N}{\log\log N}\Big)^{-\frac12}
\!\asymp\! \Big(\frac{h_N}{\log\log N}\Big)^{-\frac12}
\!\asymp  \P_N\big(\operatorname{rad}(\K)\geq r\big).
\]
Finally, suppose that $h_N\leq r\leq N/2$. We now have
$\E_N[M_r]\asymp r^2\sqrt{h_N\log(N/r)}$, and hence
\begin{align*}
\begin{split}
\P_N\big(|\K|\geq\E_N[M_r]\big)
&\!\asymp\!  \bigg(
r^{\frac45}\big(h_N\log(N/r)\big)^{\frac15}
\wedge
\frac{h_N}{
\log^\star\!\Big(
\frac{N h_N^{1/4}}
{r(h_N\log(N/r))^{1/4}}
\Big)}\bigg)^{-\frac12}\!\asymp\!  \Big(\frac{h_N}{\log(N/r)}\Big)^{-\frac12}.
%\!\asymp\! \P_N\big(\operatorname{rad}(\K)\geq r\big).
\end{split}
\end{align*}
This completes the verification of \eqref{eq:vArm} in all three regimes.
Notice that keeping track of the precise form of $h_N$ and not throwing away its logarithmic correction is essential for this relation. Indeed, we need the factor $h_N^{1/4}$ to cancel exactly with the corresponding contribution from the volume tail.
Finally, for \eqref{eq:ptvArm}, using $|\partial B_k|\asymp k(k\wedge h_N)$,
\begin{align*}
  \text{for $r\leq\Dcross$,}\quad\sum_{y\in B_r}\P_N(0\leftrightarrow y)
&\asymp \sum_{1\leq k\leq r}k^{2-1}
\asymp r^2\asymp \E_N[M_r]\cdot
\P_N\big(\operatorname{rad}(\K)\geq r\big),\text{ and} \\
\text{for $r>\Dcross$,}\quad\sum_{y\in B_r}\P_N(0\leftrightarrow y)
&\!\asymp \!\!\sum_{1\leq k\leq\Dcross}\!\!k\,
+\frac1{h_N}\sum_{\Dcross<k\leq r}k(k\wedge h_N)\log(N/k)\\
&\!\asymp \frac{r^2(r\wedge h_N)}{h_N}\log(N/r)
\asymp \E_N[M_r]\cdot
\P_N\big(\operatorname{rad}(\K)\geq r\big).
\end{align*}
\end{enumerate}
\end{remark}

\black

\subsection{Proof of Lemma~\ref{pro:manyMesoKilled}}\label{sec:matchedLemma}
We will now show Lemma~\ref{pro:manyMesoKilled} for $\lambda=2$. The argument only uses that the killing scale is comparable to $r$, and therefore applies without change, up to adjusting constants, for $2\leq \lambda\leq \Cr{range}$.
We first estimate the number of points $y\in B_{r/4}$ whose critical cluster in $\tB_{r/2}$ has large $\capacity^r$-capacity. Via \cite[(Isom)]{drewitz_2022_clust}, we will then use interlacements to connect many of these critical clusters into the cluster ${\K}^{-a}_{\nn{r}}(y)$ appearing in \eqref{eq:LBL2Killed}. The use of $\capacity^r$ below is crucial here: since we are considering local clusters at scale $r$, the subsequent connection argument must likewise take place before the relevant interlacement trajectories exit $\tB_r$. 
For $r,s>0$, define
\begin{equation}
\label{eq:defD1}
D(y;r,s)
\defeq
\Big\{
\capacity^r\!\big(\K^0_{\nn{r/2}}(y)\big) \ge s
\Big\} \quad \text{and}\quad
\mathcal D(r,s)
\defeq
\big\{y\in B_{r/4}: D(y;r,s)\text{ occurs}\big\}.
\end{equation}
%Here for $s>0$, the event $D(y;r,s)$ is empty unless $y\in \tB_{r/2}$.
We then set
\begin{equation}
\label{eq:Edef}
\mathcal E(r,s)
\defeq
\Big\{
|\mathcal D(r,s)|
\ge
\tfrac12\,\E_N^{2r}\big[|\mathcal D(r,s)|\big]
\Big\}.
\end{equation}
\begin{lemma}\label{lem:firstmommatched}
  Let $2\leq r\leq \Cr{range}N$. Then for
  all $0<s\leq \Cr{C:lawcaplocNK}\capacity^r\!(B_{r/2})$, one has that,
  \begin{align}
	\label{eq:secondmomentK}
	\E_N^{2r} \big[| \mathcal{D}(r,s)|\big]\geq
  c\frac{r^{2}(r\wedge h_N)}{\sqrt{(g_N^{2r}\,s)\vee 1}}
	\quad & \text{and}\quad
	\E_N^{2r} \big[| \mathcal{D}(r,s)|^2\big]  \le   c\frac{r^{4}(r\wedge h_N)^2}{{(g_N^{2r}\,s)\vee 1}}.
  \end{align}
  Moreover,
\begin{equation} \label{eq:ELB}
   \P_N^{2r} \big(\mathcal E(r,s )\big) \ge c .
\end{equation}
\end{lemma}
The proof of Lemma~\ref{lem:firstmommatched} is essentially identical to that of \cite[Lemma 4.2]{drewitz_2024_clust}, except for the treatment of the term
$\P^{\partial B_{2r}\cup\K^0(y)}_N\big(\capacity^r(\K^0_{\nn{r/2}}(z))\geq s\big)$,
for $y,z\in B_{r/2}$ with $z\notin\K^0(y)$, arising in the second-moment estimate. In \cite[Lemma 4.2]{drewitz_2024_clust}, the analogous term on $\Z^d$ is handled using the comparison $\capacity^r(\cdot)\leq c\capacity^{\tB_{2r}^c\cup\K^0}(\cdot)$ together with a version of \eqref{eq:clustCapDPR}, which induces the factor $g^{\tB_{2r}^c\cup\K^0}\!(0,0)^{-\frac12}$. This is harmless in $d\geq 3$, where the field is localized, but is undesirable in our setting.
\begin{proof}[Proof of Lemma~\ref{lem:firstmommatched}]
The lower bound on the first moment is a direct consequence of \eqref{eq:clustCapNoKill} and the fact that $|B_{r/4}|\geq cr^2(r\wedge h_N)$.
For the second moment, applying the strong Markov property of the field (see \eqref{eq:field_decomp}) yields that for all $y,z\in B_{r/4}$ with $y\neq z$,
\begin{align}\label{eq:secondMomStart}
\begin{split}
	&\P_N^{2r}\Big(\capacity^r\!\big({\K}^0_{\nn{r/2}}(y)\big) \geq s, \capacity^r\!\big({\K}^0_{\nn{r/2}}(z)\big) \geq s\Big)\\
	\leq &\P_N^{2r}\big(z\!\in\! \K^0(y)\big)
	+\E_N^{2r}\Big[\1\big\{z\!\notin\! \K^0(y), \capacity^r\!\big({\K}^0_{\nn{r/2}}(y)\big) \!\geq\! s\big\}\P^{\partial B_{2r}\cup\K^0(y)}_N\Big( \capacity^r\!\big({\K}^0_{\nn{r/2}}(z)\big) \!\geq\! s\Big)\Big].
\end{split}
\end{align}
Using the coupling~\cite[Proposition 2.1]{lupu_2016_from} combined with the restriction property of loop soups \cite[Section 6]{fitzsimmons_2014_marko}, one can infer that the law of $\K^{0}(z)$ under $\P_N^{2r}$ stochastically dominates the law of $\K^{0}(z)$ under $\P_N^{\partial B_{2r}\cup\K^0(y)}$. It then follows from this observation that the expectation in \eqref{eq:secondMomStart} is at most
\begin{align}\label{eq:secondMomDisc}
\begin{split}
	&\E_N^{2r}\Big[\1\{z\!\notin\! \K^0(y), \capacity^r\!\big({\K}^0_{\nn{r/2}}(y)\big) \geq s\}\P_N^{2r}\big( \capacity^r\big({\K}^0_{\nn{r/2}}(z)\big) \geq s\big)\Big]\\
	\stackrel{\eqref{eq:clustCapNoKill}}{\leq}
	& (g_N^{2r}\,s)^{-\frac12}\times \P_N^{2r}\Big(z\notin \K^0(y), \capacity^r\big({\K}^0_{\nn{r/2}}(y)\big) \geq s\Big)\stackrel{\eqref{eq:clustCapNoKill}}
  {\leq}(g_N^{2r}\,s)^{-1}.
\end{split}
\end{align}
Moreover, it follows from \cite[Proposition 5.2]{lupu_2016_from} that $\P_N^{2r}(z\in \K^0(y))\leq cg_N^{2r}(y,z)/g_N^{2r}$.
Combining this with the above two displays and a calculation similar to \cite[(3.26) and (3.27)]{rodriguez_2025_on} gives
\begin{align}\label{eq:secondmomD}
\begin{split}
\E_N^{2r} \big[| \mathcal{D}(r,s)|^2\big]
&\!=\!|B_r|\!+\!\!\!\sum_{y,z\in B_{r/4}, y\neq z}\!\P_N^{2r}\Big(\capacity^r\!\big({\K}^0_{\nn{r/2}}(y)\big) \geq s, \capacity^r\!\big({\K}^0_{\nn{r/2}}(z)\big) \geq s\Big)\\
&\!\leq\! |B_r| \Big(1\!+\!\frac{|B_r|}{(g_N^{2r} \,s)\!\vee\! 1}\Big) +\!\!
\sum_{y,z\in B_{r/4}}\!\frac{c'g_N^{2r}(y,z)}{g_N^{2r}}
\leq  |B_r|\Big(1\!+\!
\frac{|B_r|}{(g_N^{2r} \,s)\!\vee\! 1}\Big) .
\end{split}
\end{align}
To conclude, it suffices to show that $g_N^{2r} s\leq c r^2(r\wedge h_N)$.
By~\eqref{eq:killed_green_smallh} and~\eqref{eq:CapVarForm},
\begin{equation}\label{eq:gfBound}
  g_N^{2r} s
  \leq c\Big(1+\frac{\log r}{h_N}\Big)(r\wedge h_N)
  \leq c(r+\log r)
  \leq c r^2(r\wedge h_N).
\end{equation}
This concludes the proof of \eqref{eq:secondmomentK}, which also yields the lower bound on $\P_N^{2r}\!\big(\mathcal E(r,s )\big)$ by the Paley-Zygmund inequality.
\end{proof}

The remainder of the argument to deduce Lemma \ref{pro:manyMesoKilled} from Lemma \ref{lem:firstmommatched} is almost identical to that of \cite[Lemma 4.3]{drewitz_2024_clust}, so we only sketch the main ideas.

We will use the isomorphism with random interlacement (cf.~\cite[(Isom)]{drewitz_2022_clust}).
Note that by \eqref{eq:ELB}, \eqref{eq:secondmomentK} and the argument below \eqref{eq:secondmomD}, if $r\geq c$ for some large enough constant $c$, then on the event $\mathcal{E}(r,\Cr{C:lawcaplocNK}\capacity^r(B_{r/2})) $, there exists at least one $y\in B_{r/4}$ such that $\capacity^r(\K^0_{\nn{r/2}}(y))\geq \Cr{C:lawcaplocNK}\capacity^r(B_{r/2})$. Fix one such $y$ and let $\hat{\K}=\K^0_{\nn{r/2}}(y)$.
  We further write $\hat{\I}^u_r(\hat{\K})$ for the union of the forward parts of the trajectories in $\I^u$ which intersect $\hat{\K}$, started upon first entering $\hat{\K}$ and stopped upon first exiting $\tilde B_r$.
On the event $\mathcal{E}(r,\Cr{C:lawcaplocNK}\capacity^r(B_{r/2}))\cap \mathcal{E}(r,u^{-1})$, by the thinning property of Poisson variables, for any $x\in \mathcal{D}(r,u^{-1})$, the number of trajectories contributing to $\hat{\I}^u_r(\hat{\K})$ under $\bar{\P}^{2r}_N$ which intersect $\K^0_{\nn{r/2}}(x)$ is a Poisson variable with parameter at least
	\begin{align}\label{eq:matchedPoisThin}
	u\, \Cr{C:lawcaplocNK}F^r\!(r)\times \inf_{z\in{\tB_{r/2}}}P_z^{2r}\big(H_{\K^0_{\nn{r/2}}\!\!(x)}\!<H_{\tB_r^{\mathsf c}} \big)
	&\geq c F^r\!(r)\times \inf_{x,z\in B_{r/2}}g^{r}_N(x,z)
	\geq c,
	\end{align}
	where we used \eqref{eq:last_exit} and \eqref{eq:cap-killed} in the first inequality, and \eqref{eq:killed_green_smallh} in the last inequality. The Green's function $g_N^r$ appears here because the trajectories are stopped upon first exiting $\tB_r$.

	We deduce, using the previous observations along with \eqref{eq:Edef} and \eqref{eq:secondmomentK}, that for $r \ge c$,
\begin{align}
\label{eq:firstmoment2}
\begin{split}
    \E_N^{2r} \otimes \bar{\E}_N^{2r}&\Big[ 1\big\{\capacity^r\!( \hat{\K})
    \geq \Cr{C:lawcaplocNK}\capacity^r(B_{r/2})\big\} \cdot
	\sum_{x \in \mathcal{D}(r,u^{-1})}
  \1\{\K^0_{\nn{r/2}}(x)\cap  \hat{\I}^u_{r}(\hat{\K}) \ne \emptyset \}
	\Big] \\
    &\ge \E_N^{2r} \Big[ \1\{\mathcal{E}(r,\Cr{C:lawcaplocNK}\capacity^r(B_{r/2}))\} \cdot
	\sum_{x \in \mathcal{D}(r,u^{-1})} \bar{\P}_N^{2r}\big(
	\K^0_{\nn{r/2}}(x)\cap  \hat{\I}^u_{r}(\hat{\K}) \ne \emptyset
	\big)
	\Big] \\
    &\ge \P_N^{2r} \big( \mathcal{E}(r,\Cr{C:lawcaplocNK}\capacity^r(B_{r/2}))
	\,\cap\, \mathcal{E}(r,u^{-1})\big)
	\times \frac{cr^{2}(r\wedge h_N)}{\sqrt{g_N^{2r}\, u^{-1}}}.
\end{split}
\end{align}
In view of \eqref{eq:defD1}, \eqref{eq:defD1} and \eqref{eq:Edef}, the events $\mathcal{E}(r,\Cr{C:lawcaplocNK}(r\wedge h_N))$ and $\mathcal{E}(r,u^{-1})$ are increasing  in $\varphi$. Therefore, by the FKG inequality and \eqref{eq:ELB}, one can bound the probability in the last line of \eqref{eq:firstmoment2} from below by a constant.
The isomorphism \cite[(Isom)]{drewitz_2022_clust} then asserts that all the $\K^0_{\nn{r/2}}(x)$ clusters visited by trajectories contributing to $\hat{\I}^u_r(\hat{\K})$
can be viewed as part of the large cluster of interest, which is indeed connected at scale $r$ by definition.
One can then use the above display to obtain a first moment lower bound on the number of sites in $B_r$ which satisfy the volume and capacity condition in \eqref{eq:LBL2Killed}. The corresponding second moment is bounded using \eqref{eq:secondmomentK}. An application of the Paley--Zygmund inequality then concludes the proof. We refer to the proof of \cite[Lemma 4.3]{drewitz_2024_clust} for details.

\begin{remark}\label{rem:why_difficult}
We now comment on why Proposition~\ref{pro:manyMeso} is much more difficult than Lemma~\ref{pro:manyMesoKilled} in our setting.
Extending the proof of Lemma~\ref{pro:manyMesoKilled} in
Section~\ref{sec:matchedLemma} to Proposition~\ref{pro:manyMeso} would
require an analogue of Lemma~\ref{lem:firstmommatched} under $\P_N$.
In particular, one would need to estimate
\begin{equation}\label{eq:clustCapLocal}
  \P_N\Big(\capacity^r\!\big(\K^0_{\nn{y,r/4}}(y)\big)\geq s\Big)
\end{equation}
for any $y\in B_{r/4}$.
The main difficulty of proving \eqref{eq:clustCapLocal} stems from the fact that the cluster is sampled from $\P_N$ whereas $\capacity^r$ is measured with respect to the walk with additional killing upon hitting $\partial B_r$.
The mismatch is not an issue in dimensions $d\geq 3$ and \eqref{eq:clustCapLocal} follows readily from standard cluster capacity estimates. 
Recall from \eqref{eq:clustCapDPR} and \eqref{eq:clustCapNoKill} that a result akin to \eqref{eq:clustCapLocal} where $\P_N$ is replaced by $\P_N^r$ can be obtained directly from the standard capacity estimate \cite[(3.8)]{drewitz_2022_clust} and the one arm estimate.
Using the fact that $g_{\tilde{\Z}^d}(x,y)\asymp g_{\tilde{\Z}^d\setminus \tB_{r}^c}(x,y)$ for all $x,y\in \tB_{r/2}$ (see e.g.~\cite[(5.6)]{drewitz_2023_criti}), one has that by \eqref{eq:CapVarForm}, $\capacity_{\tilde{\Z}^d}\big(\K^{0}_{\nn{y,r/4}}(y)\big)\asymp \capacity_{\tilde{\Z}^d\setminus \tB_r^c}\big(\K^{0}_{\nn{y,r/4}}(y)\big)$, so the unrestricted estimate from \cite[Lemma~3.1]{drewitz_2024_clust} gives, for all $1\leq s\leq c'r^{d-2}$,
\begin{equation}\label{eq:localkilledHighD}
    cs^{-1/2}\leq \P^{\tilde{\Z}^d}\!\Big(\capacity^r\big(\K^{0}_{\nn{y,r/4}}(y)\big)\geq s\Big)\leq c's^{-1/2}.
\end{equation}
In view of~\eqref{eq:g2-asymp-unif} and~\eqref{eq:killed_green_smallh}, the above strategy fails whenever $\Fbox(r)\ll F^r(r)$. This corresponds in particular to \eqref{def:AntiCon} or the entire flat regime where $\Dcross\ll r\ll N$ (cf.~\eqref{eq:introAntiCon2}).
\end{remark}

\section{Comparison with killed clusters and gradient estimates}\label{sec:anticon}
The main result of this section is Proposition~\ref{prop:shift_h0},
which provides a way to handle the killing sensitivity phenomenon discussed
in Section~\ref{sec:proofIdea}. 
It gives a two-sided comparison between
local clusters of the full field and their killed counterparts, at the
cost of a random sprinkling governed by the harmonic contribution from
outside. 
As we will see in Section~\ref{sec:final}, the probability $\frac12$ appearing in Theorem~\ref{thm:antiCon} directly stems from the fact that the harmonic contribution is constructive or destructive with asymptotic probability $\frac12$ each in the \eqref{def:AntiCon} regime.

We will first show Proposition~\ref{prop:shift_h0} and, as an application, supply the proof of
Proposition~\ref{pro:manyMeso} in Section~\ref{sec:proof2}.
A key ingredient to show Proposition~\ref{prop:shift_h0} is a
gradient estimate for the Green's function on $\Slab$, which controls the oscillation of harmonic averages inside the killed region. We present its proof
in Section~\ref{sec:gradEst}.

\smallskip
Let $\mathcal C(\tSlab)$ denote the space of continuous functions on
$\tSlab$. For any $r\geq 1$, we say that a family of bounded functionals
$(\Fi_r^a)_{a\in \R}$, with
$\Fi_r^a:\mathcal C(\tSlab)\to\mathbb R$, is increasing at scale $r$, if for all $\phi,\phi'\in \mathcal{C}(\tSlab)$ and $a,a'\in \R$,
\begin{equation}\label{def:increasFamily}
	\{x\!\in\! \tB_r:\phi_x\!\geq\! a\}\subseteq \{x\!\in\! \tB_r:\phi'_x\!\geq\! a'\}
	\text{ implies } \Fi^a_r(\phi)\leq \Fi^{a'}_r(\phi').
\end{equation}
\begin{proposition}\label{prop:shift_h0}
There exist $\Cl[c]{rNshift}\in(0,1),\Cl[c]{antiConlb}\in(0,\infty)$ and $\Cl[c]{choiceLam}\in[2,\infty)$ such that for any $1\leq r\leq \Cr{rNshift}N$, $a\in \R$, $(\Fi_r^a)_{a\in \R}$ as in \eqref{def:increasFamily} and $b\in\R$ such that $|b|>\frac{\Cr{antiConlb}}{\sqrt{r\wedge h_N}}$,
\begin{equation}\label{eq:shift_general}
  \P_N\Big(
  \Fi^{a-b+\frac12 |b|}_r (\psi^{\Cr{choiceLam} r}) \leq 
	\Fi^{a}_r(\varphi)\leq 
	\Fi^{a-b-\frac12|b|}_r (\psi^{\Cr{choiceLam} r}) 
	\,\big|\, \eta_0^{\Cr{choiceLam} r}=b  \Big)\geq 1-\frac{\Cr{antiConlb}}{|b|\sqrt{r\wedge h_N}}.
\end{equation}
Moreover, under~\eqref{def:AntiCon}, for all $\epsilon\in(0,1)$, there exists $\Cl[c]{largeNshift}(\epsilon)<\infty$ such that for all $N\geq \Cr{largeNshift}$,
\begin{equation}\label{eq:shift_anticon}
  \P_N\Big(
  \Fi^{a-b+\epsilon |b|}_r (\psi^{2 r}) \leq 
	\Fi^{a}_r(\varphi)\leq 
	\Fi^{a-b-\epsilon|b|}_r (\psi^{2 r}) 
	\,\big|\, \eta_0^{2 r}=b  \Big)\geq 1-\frac{\Cr{antiConlb}}{\epsilon|b|\sqrt{r\wedge h_N}}.
\end{equation}
\end{proposition}
%Under \eqref{def:AntiCon}, the above comparison already holds with $\Cr{choiceLam}=2$ for all $N$ sufficiently large. 
Recall the notation from~\eqref{eq:Kr}.
A useful heuristic to keep in mind is that $\eta_0^{\Cr{choiceLam} r}\!\approx \!(r\!\wedge \!\frac{h_N}{\log^\star\!(N/r)})^{-\frac12}$
at its typical scale. Thus, upon taking
{\setlength{\abovedisplayskip}{5pt}
 \setlength{\belowdisplayskip}{5pt}
 \setlength{\abovedisplayshortskip}{5pt}
 \setlength{\belowdisplayshortskip}{5pt}
\[
\Fi_r^a(\varphi)\defeq \big|\K^a_{\nn{r}}(0,\varphi)\big|,
\]}
\!\!\!Proposition~\ref{prop:shift_h0}, followed by integration over the Gaussian density of
$\eta_0^{\Cr{choiceLam} r}$, shows that the local
cluster under $\P_N$ can be larger or smaller than its counterpart
for the field with additional killing on $\partial B_{\Cr{choiceLam}r}$, with positive probability in each case.

To show Proposition \ref{prop:shift_h0}, we need the following gradient estimate for the Green's function on $\tSlab$, which is of independent interest. In the sequel, we write $g_N(x)=g_N(0,x)$ for $x\in \tSlab$.
\begin{proposition}[Gradient estimate] \label{prop:gDrad}
For all $N\geq 1$, every coordinate vector $e$ and $x=(y,z)\in \Slab$ such that $\|x\|\leq\Cr{range}N$, one has
\begin{equation}
\label{eq:ggradDisc}
  |g_N(x+e)-g_N(x)|
  \leq
  c\left(
    \frac{1}{(\|x\|\vee 1)^2}
    +
    \frac{1}{h_N(|y|\vee h_N)}
  \right).
\end{equation}
Moreover, for all $\lambda\geq 2$, $1\leq r\leq\Cr{range}N$ and all $x=(y,z),x'=(y',z')\in \tSlab$ such that $\|x\|\leq \Cr{range}N$, $|x|\wedge |x'|\geq \lambda r$ and  $\|x-x'\|\leq r$, we have that
\begin{equation}\label{eq:ggradCable}
	|g_N(x)-g_N(x')|\leq c(\|x-x'\|\vee 1)\left(
    \frac{1}{(\lambda r)^2}
    +
    \frac{1}{h_N((\lambda r)\vee h_N)}
  \right).
\end{equation}
\end{proposition}

We postpone the proof of Proposition \ref{prop:gDrad} to the end of this section and first explain how to show Proposition \ref{prop:gDrad}. 
In view of \eqref{eq:field_decomp}, we abbreviate $\eta^r_x=\eta^{\partial B_{r}}_x$ and $\psi^r_x=\psi^{\partial B_{r}}_x$.
We first derive the following useful consequence of Proposition \ref{prop:gDrad}, which controls the oscillation of the harmonic average inside $\tB_r$.

\begin{lemma}\label{lem:oscillation}
For all $\lambda\geq 2$,
\begin{equation}\label{eq:expOscillation}
	\E_N\Big[\sup_{x\in \tB_r}\big|\eta^{\lambda r}_x-\eta^{\lambda r}_0\big|\Big]\leq \frac{c}{\sqrt{\lambda}\sqrt{ (\lambda r)\wedge h_N}}.
\end{equation}
\end{lemma}
\begin{proof}
  We abbreviate $Z_x:=\eta^{\lambda r}_x-\eta^{\lambda r}_0$.  Since $\eta^{\lambda r}$ is harmonic in $\tB_r$ in the sense of \cite[Definition 5.1]{prevost_2023_perco}, it is affine on every cable in $\tB_r$. This implies $\sup_{x\in\tB_r}|Z_x|=\sup_{x\in B_r}|Z_x|$.  Hence it suffices to show \eqref{eq:expOscillation} with $\tB_r$ replaced by $B_r$.
  We start by controlling 
    $d_c(x,y):=\E_N[(Z_x-Z_y)^2]^{1/2}$,
    the so called ``canonical metric" on $B_r$ (cf.~\cite[1.3.1]{_2007_rando}). It follows from \eqref{eq:field_decomp} and \eqref{eq:gDecomp} that for any $x,y\in B_r$,
	\begin{align*}
		\E_N[(Z_x-Z_y)^2]&=\E_N[(\varphi_x-\psi^{\lambda r}_x-\varphi_y+\psi^{\lambda r}_y)^2]\\
		&=g_N(x,x)-g_N^{\lambda r}(x,x)+g_N(y,y)-g_N^{\lambda r}(y,y)-2\big(g_N(x,y)-g_N^{\lambda r}(x,y)\big)\\
		&=E_x\big[g_N(x,X_{H_{\partial B_{\lambda r}}})-g_N(y,X_{H_{\partial B_{\lambda r}}}),H_{\partial B_{\lambda r}}<\infty\big]\\
		&\qquad\qquad\qquad+E_y\big[g_N(y,X_{H_{\partial B_{\lambda r}}})-g_N(x,X_{H_{\partial B_{\lambda r}}}),H_{\partial B_{\lambda r}}<\infty\big]\\
		&\leq 2\times \sup_{z\in \partial B_{\lambda r}}\big|g_N(x,z)-g_N(y,z)\big|.
	\end{align*}
	By the above display and applying \eqref{eq:ggradCable} to $|g_N(x-z)-g_N(y-z)|$ for any $x,y\in B_r$ and $z\in \partial B_{\lambda r}$, we have that 
	\begin{equation}\label{eq:canonMetricBound}
		d_c(x,y)\leq \E_N[(Z_x-Z_y)^2]^{\frac12}\leq c(\sqrt{\|x-y\|}\vee 1)\left(
    	\frac{1}{\lambda r}+\frac{1}{\sqrt{h_N((\lambda r)\vee h_N)}}\right).
	\end{equation} 
    Let $\diam_{d_c}(B_r)$ denote the diameter of $B_r$ under $d_c$ and for $\epsilon>0$, let $N(\epsilon)$ denote the number of $d_c$ balls of radius $\epsilon$ needed to cover $B_r$.  It follows from \eqref{eq:canonMetricBound}, that
	\begin{equation}
		\diam_{d_c}(B_r)\leq c\left(\frac{1}{\lambda \sqrt{r}}+\frac{1}{\sqrt{\lambda h_N}}\right)\leq \frac{c}{\sqrt{\lambda}\sqrt{ (\lambda r)\wedge h_N}}.
	\end{equation}
	Furthermore, \eqref{eq:canonMetricBound} implies that there exists $\Cl[c]{dcEntropy}>0$ such that for any $x,y\in B_r$ with $\|x-y\|\leq \Cr{dcEntropy}\epsilon^2\lambda r((\lambda r)\wedge h_N)$, we have $d_c(x,y)\leq \epsilon$. Hence $N(\epsilon)\leq c(r/(\epsilon^2\lambda r((\lambda r)\wedge h_N)))^3$. Thus by \cite[Theorem 1.3.3]{_2007_rando} and writing $\beta=\lambda ((\lambda r)\wedge h_N)$, we have
	\[
	\E_N\big[\sup_{x\in B_r}Z_x\big]
    \leq c \int_0^{c'\beta^{-\frac12}} \sqrt{\log\left(\frac{c''}{\epsilon^2\beta}\right)} \,d\epsilon 
	\leq \frac{c}{\sqrt{\beta}}\int_{c'}^{\infty}\sqrt{t}e^{-t/2}\,dt\leq \frac{c}{\sqrt{\beta}}.
	\]
	 where we did the change of variable $t=\log(\tfrac{c''}{\epsilon^2\beta})$ for the second inequality.
\end{proof}

\begin{proof}[Proof of Proposition \ref{prop:shift_h0}]
By \eqref{def:increasFamily} and \eqref{eq:field_decomp}, it's easy to see that for any $\lambda \geq 2$  soon to be fixed and any $\epsilon\in (0,1)$,
\begin{align}\label{eq:field_comp}
\begin{split}
%& 
\Big\{\sup_{x\in \tB_r}|\eta^{\lambda r}_x-\eta^{\lambda r}_0|\leq \epsilon|b|, \eta^{\lambda r}_0=b\Big\}\subseteq 
\Big\{\Fi^{a-b+\epsilon|b|}\big(\psi^{\lambda r}\big)\leq \Fi^{a}_r(\varphi)\leq \Fi^{a-b-\epsilon|b|}\big(\psi^{\lambda r}\big)\Big\} .
%\quad\text{and},\\	& \Big\{\sup_{x\in \tB_r}|\eta^{\lambda r}_x-\eta^{\lambda r}_0|\leq \frac12b, \eta^{\lambda r}_0=-b\Big\}\subseteq \{ \Fi^{a}_r(\varphi)\leq \Fi^{a+\frac12 b}(\psi^{\lambda r})\}.	
\end{split}
\end{align}
We now study the conditional distribution of $\sup_{x\in \tB_r}|\eta^{\lambda r}_x-\eta^{\lambda r}_0|$. Let $\alpha_x\!=\!\alpha_x(\lambda)\!\defeq\!\frac{\cov_N(\eta_x^{\lambda r},\eta_0^{\lambda r})}{\variance_N(\eta_0^{\lambda r})}$. Using the standard formula for the conditional mean of a multivariate Gaussian (see, e.g. \cite[1.2.7]{_2007_rando}), we have that $\E_N[\eta_x^{\lambda r}-\eta_0^{\lambda r}\,|\,\eta^{\lambda r}_0]=(\alpha_x-1)\eta^{\lambda r}_0.$
Now for any $x\in \tB_r$, consider
\begin{equation}\label{eq:def-phi}
\phi^{\lambda r}_x\defeq \eta_x^{\lambda r}-\eta_0^{\lambda r}-(\alpha_x-1)\eta^{\lambda r}_0.
\end{equation}
Since $\E_N[\phi^{\lambda r}_x\eta^{\lambda r}_0]=0$, for all $x\in \tB_r$, $(\phi^{\lambda r}_x)_{x\in \tB_r}$ is a centered Gaussian field that's independent of $\eta^{\lambda r}_0$. Further, conditional on the event $\big\{\eta^{\lambda r}_0=b\big\}$, we have that by \eqref{eq:def-phi} and the reverse triangle inequality,
\begin{equation} \label{eq:phi_eta}
	\Big\{\sup_{x\in \tB_r}|\eta^{\lambda r}_x-\eta^{\lambda r}_0|> \epsilon|b|\Big\}
	\subseteq \Big\{\sup_{x\in \tB_r}|\phi^{\lambda r}_x|> \big(\epsilon -\sup_{x\in \tB_r}|\alpha_x-1|\big) |b|\Big\}.
\end{equation}
For any $x,y\in \tB_r$, $\E_N[(\phi^{\lambda r}_x-\phi^{\lambda r}_y)^2]^{\frac12}\leq \E_N[(\eta^{\lambda r}_x-\eta^{\lambda r}_y)^2]^{\frac12}$, i.e. the canonical metric of $(\phi^{\lambda r}_x)_{x\in \tB_r}$ can be bounded by the canonical metric of the field considered in Lemma~\ref{lem:oscillation}. It then follows from the proof of Lemma \ref{lem:oscillation} that,
\begin{equation}\label{eq:phiExpCtrl}
	\E_N\Big[\sup_{x\in \tB_r}|\phi^{\lambda r}_x|\Big]\leq \frac{c}{\sqrt{\lambda}\sqrt{ (\lambda r)\wedge h_N}}.
\end{equation}

In view of \eqref{eq:phi_eta}, it remains to control $\sup_{x\in \tB_r} |\alpha_x-1|$.
By \eqref{eq:exitProb}, \eqref{eq:field_decomp}, \eqref{eq:hitting_lb} and \eqref{eq:gDecomp}, we have that
\begin{align}\label{eq:varEta0}
\begin{split}
\variance_N(\eta^{\lambda r}_0)&=g_N(0,0)-g_N^{\lambda r}(0,0)
=E_0\big[g_N(0,X_{H_{\partial B_{\lambda r}}}),H_{\partial B_{\lambda r}}<\infty \big]	\\
&\geq \!\inf_{x\in \partial B_{\lambda r}}g_N(0,x)P_0(H_{\partial B_{\lambda r}}<\infty)\geq ce^{-\frac{(\lambda r)^2}{N^2}}\left(\frac{1}{\lambda r}\!+\!\frac{1}{h_N}\log\Big(\frac{N}{(\lambda r)\vee h_N}\!\vee\! 2\Big)\right).
\end{split}
\end{align}
Similarly, combining \eqref{eq:field_decomp}, \cite[(1.49)]{sznitman_2012_topic} and \eqref{eq:ggradCable} yields,
\begin{align*}
	\quad|\cov_N(\eta^{\lambda r}_x,\eta^{\lambda r}_0)-\variance_N(\eta^{\lambda r}_0)|
	&=|g_N(x,0)-g_N^{\lambda r}(x,0)-(g_N(0,0)-g_N^{\lambda r}(0,0))|\\
	&\leq E_0\Big[|g_N(x,X_{H_{\partial B_{\lambda r}}})-g_N(0,X_{H_{\partial B_{\lambda r}}})|,H_{\partial B_{\lambda r}}<\infty\Big]\\
	&\leq \sup_{z\in \partial B_{\lambda r}}|g_N(x,z)-g_N(0,z)|
	\leq \frac{1}{\lambda^2r}+\frac{1}{\lambda h_N}.
\end{align*}
Hence, using the last two displays and picking $1\leq r\leq \Cr{rNshift}N$ with $\Cr{rNshift}\in(0,1)$ small enough, we can choose $\lambda=\Cr{choiceLam}$ for a large enough constant $\Cr{choiceLam}$ such that,
\begin{equation}\label{eq:alphaBound}
\sup_{x\in \tB_r} |\alpha_x(\lambda)-1|\leq ce^{\frac{(\lambda r)^2}{N^2}}\frac{(\lambda r)\wedge h_N}{(\lambda^2r)\wedge (\lambda h_N)}\leq \frac{c}{\lambda }e^{\frac{(\lambda r)^2}{N^2}}\leq \frac{1}{4}.
\end{equation}
Moreover, in view of~\eqref{def:Fbox} and~\eqref{eq:defn_fbox},under~\eqref{def:AntiCon}, we can pick $\Cr{largeNshift}(\epsilon)$ large enough such that for all $N\geq \Cr{largeNshift}$,
\begin{equation}\label{eq:alphaBound_small}
\sup_{x\in \tB_r} |\alpha_x(2)-1|\leq ce^{\frac{r^2}{N^2}}\frac{F(r)}{F^{2r}(r)}\leq \frac\epsilon2.
\end{equation}
We may now use \eqref{eq:alphaBound}, Markov's inequality and \eqref{eq:phiExpCtrl} to obtain that,
\begin{equation}\label{eq:markovPhi}
	\P_N\Big(\sup_{x\in \tB_r}|\phi^{\Cr{choiceLam} r}_x|>\big(\frac{1}{2}-\sup_{x\in \tB_r}|\alpha_x(\Cr{choiceLam})-1|\big) |b|\Big)\leq \frac{c\,}{|b|\sqrt{ r\wedge h_N}}\wedge 1.
\end{equation}
Combining this with~\eqref{eq:field_comp} and~\eqref{eq:phi_eta} with $\epsilon=\frac12$ concludes the proof of~\eqref{eq:shift_general}. The proof of~\eqref{eq:shift_anticon} follows similarly using~\eqref{eq:alphaBound_small} instead of~\eqref{eq:alphaBound}.
\end{proof}

\subsection{Proof of Proposition \ref{pro:manyMeso} via Proposition~\ref{prop:shift_h0}}\label{sec:proof2}

\begin{proof}
For $\phi\in \mathcal C(\tSlab)$ consider, for $a,r$ as in the statement of Proposition~\ref{pro:manyMeso} and $q\in \R$,
\[
\Theta_r^q(\phi)
\defeq
\1\bigg\{
\exists\,y\in B_r:\,
\big|\K^q_{\nn{r}}(y,\phi)\big|
\geq
\Cr{c:LBhatK}\frac{|a|r^2(r\wedge h_N)}{\sqrt{g_N}},
\quad
\capacity^{r}\big(\K^q_{\nn{r}}(y,\phi)\big)
\geq
\Cr{MarCap}(r\wedge h_N)
\bigg\}.
\]
and note that the family $(\Theta_r^q)_{q\in\R}$ is increasing at scale $r$ in the sense of \eqref{def:increasFamily}.
For $r\geq 1$ and $s\geq 1$ soon to be fixed, let $f(s)\defeq \frac{s}{\sqrt{F(r)}}.$
We first claim that, for all $s\geq 1$, $\Cr{proprLB}\leq r\leq \Cr{range}N$, $a\in(0,\sqrt{g_N}]$ such that $a^2F(r)\geq\Cr{C:LBLneg}$ as in Proposition \ref{pro:manyMeso},
\begin{equation}\label{eq:thetaKilledShift}
\P_N\big(
\Theta_r^{-a-f(s)}
(\psi^{\Cr{choiceLam}r})=1
\big)\geq c_0>0.
\end{equation}
To see this, set $\tilde a\defeq (a+f(s))\wedge\sqrt{g_N^{\Cr{choiceLam}r}}$. By \eqref{eq:killed_green_smallh} and \eqref{eq:defn_fbox}, $\sqrt{g_N^{\Cr{choiceLam}r}}F^r\!(r)\geq \Cr{C:LBLneg2}$ for $r\geq \Cr{proprLB}$ large enough. Additionally
$a^2F^r(r)\geq a^2F(r) \geq \Cr{C:LBLneg2}$ for all $a^2F(r)\geq\Cr{C:LBLneg}$ with $\Cr{C:LBLneg}$ large enough. It then follows that $\tilde a F^r(r)\geq \Cr{C:LBLneg2}$ and we may apply Lemma~\ref{pro:manyMesoKilled} to get that with probability at least $c_0$ there exists some $y\in B_r$ such that
\begin{equation}\label{eq:locVolSupcrit}
  \big|\K^{-\tilde a}_{\nn{r}}(y,\psi^{\Cr{choiceLam}r})\big|
\geq c\frac{\tilde a r^2(r\wedge h_N)}{\sqrt{g_N^{r}}},
\qquad
\capacity^{r}\!
\big(\K^{-\tilde a}_{\nn{r}}(y,\psi^{\Cr{choiceLam}r})\big)
\geq c(r\wedge h_N).
\end{equation}
Since $a\leq \sqrt{g_N}$ and ${g_N}\geq g_N^r$, we get by \eqref{eq:killed_green_smallh} that $\tilde a/\sqrt{g_N^{r}}\geq c (\tfrac{a}{\sqrt{g_N}} \wedge 1)\geq c\,a/\sqrt{g_N}$. Combining this with the last display, and the fact that 
$\K^{-\tilde a}_{\nn{r}}(y;\psi^{\Cr{choiceLam}r})
\subseteq \K^{-a-f(s)}_{\nn{r}}(y,\psi^{\Cr{choiceLam}r})$
yields \eqref{eq:thetaKilledShift}.

We now transfer \eqref{eq:thetaKilledShift} to $\varphi$ using Proposition~\ref{prop:shift_h0}. Since $\Fbox(r)\leq r\wedge h_N$ by~\eqref{def:Fbox}, there exists $\Cl[c]{choice_s}<\infty$ such that $f(\Cr{choice_s})\geq \Cr{antiConlb}/(c_0\sqrt{r\wedge h_N})$. We henceforth fix $s=\Cr{choice_s}$.
For $b\geq2f(\Cr{choice_s})$, Proposition~\ref{prop:shift_h0} and the fact that
$(\Theta_r^a)_{a\in\R}$ is increasing give
\begin{equation}\label{eq:applyshiftpropforprop}
\P_N\Big(
\Theta_r^{-a}(\varphi)
<
\Theta_r^{-a-f(s)}
(\psi^{\Cr{choiceLam}r})
\,\Big|\,
\eta^{\Cr{choiceLam}r}_0=b
\Big)
\leq
\frac{\Cr{antiConlb}}{b\sqrt{r\wedge h_N}}\leq c_0/2.
\end{equation}
Since $\psi^{\Cr{choiceLam}r}$ and $\eta^{\Cr{choiceLam}r}_0$ are independent by~\eqref{eq:field_decomp}, \eqref{eq:applyshiftpropforprop} and~\eqref{eq:thetaKilledShift} give, for all $b\!\geq\!2f(\Cr{choice_s})$,
\[\P_N\big(
\Theta_r^{-a}(\varphi)=1
\,\big|\,
\eta^{\Cr{choiceLam}r}_0=b
\big)
\geq
\P_N\big(
\Theta_r^{-a-f(s)}
(\psi^{\Cr{choiceLam}r})=1
\big)
-
\frac{\Cr{antiConlb}}{b\sqrt{r\wedge h_N}}
\geq
c_0/2.\]
To conclude, integrating over $\eta^{\Cr{choiceLam}r}_0$ and using \eqref{eq:varEta0} together with the Gaussian density gives
\[
\P_N\big(\Theta_r^{-a}(\varphi)=1\big)
\geq
c\int_{2f(\Cr{choice_s})}^{4f(\Cr{choice_s})}
\sqrt{F(r)}
\exp\{-cF(r)b^2\}\,db \geq c.
\]
\end{proof}

\subsection{Gradient estimate}\label{sec:gradEst}
In this section, we prove Proposition~\ref{prop:gDrad}.
We first record a standard gradient estimate for the continuous-time heat kernel. Its proof is deferred to Section~\ref{A:grad}.
\begin{lemma}
\label{lem:ct-hk-gradient}
For $d\geq1$ and $r\in(0,1]$, let $p_t^{d,r}(\cdot)$ denote the
continuous-time heat kernel of the simple random walk on $\Z^d$ with
jump rate $r$. There exist
$\Cl[c]{hkg}(r,d),\Cl[c]{hkg_exp}(r,d)\in(0,\infty)$ such that, for all
$y\in\Z^d\setminus\{0\}$, $t\geq |y|$ and every coordinate vector $e$,
\begin{equation}
\label{eq:ct-hk-gradient}
  \big|p_t^{d,r}(y+e)-p_t^{d,r}(y)\big|
  \leq
  \Cr{hkg} t^{-(d+1)/2}
  \exp\left\{-\Cr{hkg_exp}\frac{|y|^2}{t}\right\}.
\end{equation}
\end{lemma}

In the sequel, we write $x=(y,z) \in \mathbb{S}_N$ with $y\in \mathbb{Z}^2$ and $z \in (\Z / h_N \Z)$. 
To show \eqref{prop:gDrad}, we use the same decomposition of $g_N$ as in \cite[(2.6)]{rodriguez_2025_on}, which we briefly recall here.
Let $Y$ and $Z$ be independent continuous-time simple random walks on $\Z^2$ and $\Z$, with jump rates $2/3$ and $1/3$, respectively, and write $P^2$ and $P^1$ for their laws (i.e.~$P_0^2(Y_t\!=\!y)\!=\!p_t^{2,2/3}\!(y)$ and $P_0^1(Z_t\!=\!z)\!=\!p_t^{1,1/3}\!(z)$). Let $P\defeq P_0^2\otimes P_0^1$, then, writing $g_N(x)=g_N(0,x)$,
\begin{equation}\label{eq:gN_decomp}
    g_N(x)=\lambda_x^{-1}\big(g_N^2(x)+g_N^3(x)\big),\quad\text{where}
\end{equation}
\begin{align}
g_N^3(x)
&\defeq
\int_0^\infty
P\big((Y_t,Z_t)=(y,\hat z)\big)
e^{-t/N^2}\,dt,
\label{eq:gN3}\\
g_N^2(x)
&\defeq
\int_0^\infty
P_0^2(Y_t=y)e^{-t/N^2}
\sum_{k\in\Z\setminus\{0\}}
P_0^1(Z_t=\hat z+kh_N)\,dt.
\label{eq:gN2}
\end{align}
We will estimate the gradients of these two terms separately in Lemma~\ref{lem:g3-grad} and Lemma~\ref{lem:g2-grad}.

We need some shorthand for the neighbour of $x$ in direction $e$.
Recall from below \eqref{def:slab} that we denote $\hat z$ as the representative of $z\in\mathbb Z/h_N\mathbb Z$ satisfying
$|\hat z|=d_{\mathbb Z/h_N\mathbb Z}(0,z)\leq \frac{h_N}{2}.$ Now
for any coordinate vector $e$ in $\Z^2$ and $\hat{x}=(y,\hat{z})$, let
\begin{equation*}
\hat{x}_e=(y_e,\hat{z}_e)=
	\begin{cases}
		(y+e,\hat{z}),
		& \text{ if $e\in \{\pm e_1,\pm e_2\}$},\\
		(y,\hat{z+e}),
		& \text{ if $e\in \{\pm e_3\}$}.
	\end{cases}
\end{equation*}

\begin{lemma}\label{lem:g3-grad}
For all $N\geq 1$, coordinate vector $e$ and $x\in\Slab$ such that $\|x\|\leq\Cr{range}N$, one has
\begin{equation}
\label{eq:g3-gradient-edge}
  |g_N^3(x+e)-g_N^3(x)|
  \leq \frac{c}{(\|x\|\vee 1)^2}.
\end{equation}	
\end{lemma}
\begin{proof}
Let $M\geq 1$ to be fixed soon.
Note that we may without loss of generality assume $\|x\|\geq M$.
Abbreviating $q(\hat{x},t):=P((Y_t,Z_t)=\hat{x})$,
it follows from \eqref{eq:gN3} that 
\begin{align*}
\begin{split}
|g_N^3(x+e)-g_N^3(x)|
\leq &\int^\infty_0 |q(\hat{x}_e,t)-q(\hat{x},t)|\,dt\\
\leq &2\int^{|\hat x|}_0 \max\big(q(\hat{x}_e,t),q(\hat{x},t)\big)\,dt+\int^\infty_{|\hat x|} |q(\hat{x}_e,t)-q(\hat{x},t)|\,dt.
\end{split}
\end{align*}
Note that $(Y,Z)$ has the law of a unit rate simple random walk on $\Z^3$ starting at the origin. It then follows from \cite[Lemma A.1]{rodriguez_2025_on} and $|\hat x|=\|x\|$ that
\[
\int^{|\hat x|}_0 \max\big(q(\hat{x}_e,t),q(\hat{x},t)\big)\,dt\leq c\|x\|e^{-c'\|x\|}\leq c\|x\|^{-2},
\]
where the last inequality follows by choosing $M$ to be a large enough constant. Further, since $|\hat{x}-\hat{x}_e|=1$, \eqref{eq:ct-hk-gradient} yields
\[
\int^\infty_{|\hat x|} |q(\hat{x}_e,t)-q(\hat{x},t)|\,dt
\leq 
\int^\infty_{|x|} c t^{-2}e^{-c'|\hat x|^2/t}\,dt
\leq c\|x\|^{-2}.
\]
\end{proof}
We will use the standard heat kernel bound frequently in the sequel (see, e.g.~\cite[Theorem 6.28(b)]{barlow_2017_rando}). For all $d\geq1$, 
$r\in(0,1]$ and $t\geq 1\vee |y|$,
\begin{equation}\label{eq:standard_hk}
	p_t^{d,r}(y)
  \leq
  c(r)\, t^{-d/2}
  \exp\Big\{-c'(r)\frac{|y|^2}{t}\Big\}.
\end{equation}
\begin{lemma}\label{lem:g2-grad}
For all $N\geq 1$, coordinate vector $e$ and $x=(y,z)\in\Slab$ such that $\|x\|\leq\Cr{range}N$, one has
\begin{equation}
\label{eq:g2-grad}
  |g_N^2(x+e)-g_N^2(x)|
  \leq \frac{c}{h_N(|y|\vee h_N)}.
\end{equation}	
\end{lemma}
\begin{proof}
Abbreviating $q^2_t(y)=P_0^2(Y_t=y)$ and $q^1_t(z)=P_0^1(Z_t=y)$, it follows from
\eqref{eq:gN2} that
\begin{align}
\label{eq:g2_vert}
  |g_N^2(x+e)\!-\!g_N^2(x)|
  &\leq 
  \int_0^\infty
  |q^2_t(y+e)\!-\!q^2_t(y)|
  \!\sum_{k\in\mathbb Z\setminus \{0\}}\!q^1_t(\hat z+k h_N)
   \,dt , \text{ if $e\in\{\pm e_1,\pm e_2\}$;}\\
\label{eq:g2_hori}
  |g_N^2(x+e)\!-\!g_N^2(x)|
  &\leq 
  \int_0^\infty
  q^2_t(y)
  \!\sum_{k\in\mathbb Z\setminus \{0\}}\! |q^1_t(\hat z_e+k h_N)\!-\!q^1_t(\hat z+k h_N)|\,dt , \text{ if $e\in\{\pm e_3\}$.}
\end{align}
We start by collecting a few useful estimates on the periodic parts of \eqref{eq:g2_vert} and \eqref{eq:g2_hori}. Let $I_t=\{k\in \Z\setminus \{0\}:|\hat z+kh_N|\leq t\}$. Note that for any $|k|\geq 1$ and $z\in \Z/h_N\Z$, we have that $|\hat{z}+kh_N|\geq (k-\frac12)^2h_N^2$.
Combining this with \eqref{eq:standard_hk} and \cite[Lemma A.1]{rodriguez_2025_on} yields that for any $0<t\leq 2h_N^2$,
\begin{align}\label{eq:sum_smallt_z}
\begin{split}
	\sum_{k\in\mathbb Z\setminus \{ 0\}}\!\!q^1_t(\hat z+k h_N)
	&\!=\! \sum_{k\in I_t}\!q^1_t(\hat z+k h_N)\!+\!\sum_{k\notin I_t}\!q^1_t(\hat z+k h_N)\\
	&\!\leq\! c\Big(\sum_{k=1}^\infty t^{-\frac12}e^{-\frac{c'(k-\frac12)^2h_N^2}{t}}\! +\!\sum_{k=1}^\infty e^{-c''(k-\frac12)h_N} \Big)
	\!\leq\! c\Big( t^{-\frac12}e^{-\frac{c'h_N^2}{t}}\! +\! e^{-c''h_N} \Big).
\end{split}
\end{align}
Similarly, using $|\hat z-\hat z_e|=1$, \eqref{eq:ct-hk-gradient} and \cite[Lemma A.1]{rodriguez_2025_on} yields that for any $0<t\leq 2h_N^2$,
\begin{align}\label{eq:sum_smallt_zgrad}
\begin{split}
	\sum_{k\in\mathbb Z\setminus \{0\}}\!|q^1_t(\hat z_e+k h_N)-q^1_t(\hat z+k h_N)|
	&\leq c\left(\sum_{k=1}^\infty t^{-1}e^{-\frac{c'(k-\frac12)^2h_N^2}{t}} +\sum_{k=1}^\infty e^{-c''(k-\frac12)h_N} \right)\\
	&\leq c\left( t^{-1}e^{-\frac{c'h_N^2}{t}} + e^{-c''h_N} \right).
\end{split}
\end{align}
Furthermore, by \eqref{eq:standard_hk}, \eqref{eq:ct-hk-gradient} and \cite[(2.24)]{rodriguez_2025_on}, we have that for all $t\geq 2h_N^2$,
\begin{equation}\label{eq:sum_bigt_z}
\sum_{k\in\mathbb Z\setminus \{ 0\}}\!\!q^1_t(\hat z+k h_N)\leq ch_N^{-1}
\quad\text{and}\quad
\sum_{k\in\mathbb Z\setminus \{0\}}\!|q^1_t(\hat z_e+k h_N)-q^1_t(\hat z+k h_N)|\leq ct^{-\frac12}h_N^{-1}.
\end{equation}
We are now ready to show \eqref{eq:g2-grad}. For $K\geq 1$ to be fixed soon, we may assume $|y|\vee h_N\geq K$ since otherwise \eqref{eq:g2-grad} can be made trivially true by increasing the constant $c$.  We start with the case \eqref{eq:g2_hori}. By \cite[(2.17)]{rodriguez_2025_on} with $M=2$, and \eqref{eq:sum_smallt_z}, we have that
\begin{align*}
&\int_0^{(2h_N^2)\vee|y|}|q^2_t(y+e)\!-\!q^2_t(y)|
  \!\sum_{k\in\mathbb Z\setminus \{0\}}\!q^1_t(\hat z+k h_N)\,dt\\
\leq &e^{-c(h_N\vee|y|)}+\int_{h_N\vee |y|}^{(2h_N^2)\vee|y|}|q^2_t(y+e)\!-\!q^2_t(y)|
  \times\Big( ct^{-\frac12}e^{-\frac{c'h_N^2}{t}}\! +\! e^{-c''h_N} \Big)\,dt.
\end{align*}
It's clear that $e^{-c(h_N\vee|y|)}\leq c'(h_N(|y|\vee h_N))^{-1}$. By \eqref{eq:ct-hk-gradient}, the integral on the last line is at most
\begin{align}\label{eq:g2_smallt_int}
\begin{split}
&c\int_{h_N\vee |y|}^{(2h_N^2)\vee|y|}t^{-2}e^{-c'\frac{|y|^2+h_N^2}{t}}\,dt
+ce^{-ch_N}\int_{h_N\vee |y|}^{(2h_N^2)\vee|y|}t^{-\frac32}e^{-c'\frac{|y|^2}{t}}\,dt\\
\leq &\frac{c}{|y|^2+h_N^2}+\frac{c'e^{-c''h_N}}{|y|}\leq \frac{c}{h_N(|y|\vee h_N)}.
\end{split}
\end{align}
For the remaining part of the integral, the observation $\int_0^su^{-1/2}e^{-u}\,du\leq c(1\wedge \sqrt{s})$, together with \eqref{eq:sum_bigt_z}, \eqref{eq:ct-hk-gradient} and the substitution $u=\frac{|y|^2}{t}$ gives,
\begin{align}\label{eq:g2_bigt_int}
\begin{split}
&\int_{(2h_N^2)\vee|y|}^{\infty}|q^2_t(y+e)\!-\!q^2_t(y)|
  \!\sum_{k\in\mathbb Z\setminus \{0\}}\!q^1_t(\hat z+k h_N)\,dt\\
\leq &\frac{c}{h_N}\int_{(2h_N^2)\vee|y|}^{\infty}t^{-\frac32}e^{-c'\frac{|y|^2}{t}} \,dt\leq \frac{c}{h_N|y|}\int_{0}^{\frac{|y|^2}{(2h_N^2)\vee|y|}} u^{-\frac12}e^{-u}\,du\leq \frac{c}{h_N(h_N\vee|y|)},
\end{split}
\end{align}
this concludes the proof in the case \eqref{eq:g2_hori}. 
For \eqref{eq:g2_vert}, it suffices to carry out similar estimates to those above, using \eqref{eq:sum_smallt_zgrad} and the gradient estimate in \eqref{eq:sum_bigt_z} instead.
Indeed, by \cite[(2.17)]{rodriguez_2025_on} with $M=2$, \eqref{eq:standard_hk}, \eqref{eq:sum_smallt_zgrad} and the same calculation as in \eqref{eq:g2_smallt_int}, one has
\begin{align*}
&\int_0^{(2h_N^2)\vee|y|}
  q^2_t(y)
  \!\sum_{k\in\mathbb Z\setminus \{0\}}\! |q^1_t(\hat z_e+k h_N)\!-\!q^1_t(\hat z+k h_N)|\,dt\\
\leq &e^{-(h_N\vee |y|)}+\int_{h_N\vee|y|}^{(2h_N^2)\vee|y|}
  t^{-1}e^{-c'\frac{|y|^2}{t}} 
  \times\Big( t^{-1}e^{-\frac{c'h_N^2}{t}} + e^{-c''h_N} \Big)\,dt \leq \frac{c}{h_N(|y|\vee h_N)}.
\end{align*}
Similarly, it follows from \eqref{eq:sum_bigt_z} and \eqref{eq:standard_hk} and \eqref{eq:g2_bigt_int} that 
\[\int_{(2h_N^2)\vee|y|}^\infty
  \!q^2_t(y)\!
  \!\sum_{k\in\mathbb Z\setminus \{0\}}\! |q^1_t(\hat z_e+k h_N)-q^1_t(\hat z+k h_N)|\,dt
\!\leq\! \int_{(2h_N^2)\vee|y|}^\infty \!\frac{ct^{-\frac12}}{h_N}
  t^{-1}e^{-c'\frac{|y|^2}{t}} \,dt\leq \frac{c}{h_N(|y|\vee h_N)}.\]
\end{proof}

We are now ready to finish the

\begin{proof}[Proof of Proposition~\ref{prop:gDrad}]
Note that \eqref{eq:ggradDisc} follows directly from Lemma~\ref{lem:g3-grad} and Lemma~\ref{lem:g2-grad}.
For \eqref{eq:ggradCable}, if $x=x'$, the claim is trivial. Otherwise, let
$\bar x,\bar x'\in\Slab$ be endpoints of cables containing $x$ and $x'$,
respectively, chosen so that $\bar x\neq \bar x'$. For some $m\geq1$, let $x_1=\bar{x}$ and $x_m=\bar{x}'$ and let $(x_i)_{0\leq i\leq m}$ be a nearest-neighbour path in $\Slab$ from
$\bar{x}$ to $\bar{x}'$ with $m\leq c(\|x-x'\|\vee1)$ and $\|x-x_i\|\leq \|x-x'\|$ for all $1\leq i\leq m$. Write $x_i=(y_i,z_i)$. It then follows from this setup that $\|x_i\|\geq \lambda r/2$ for all $1\leq i\leq m$.
Suppose first that we in fact have $\min_{1\leq i\leq m}|y_i|\geq \lambda r/4$. Combining \eqref{eq:ggradDisc} and the above considerations gives,
\begin{align}\label{eq:gradTelescope}
\begin{split}
|g_N(\bar{x})-g_N(\bar{x}')|\leq \sum_{i=1}^{m-1}|g_N(x_i)-g_N(x_{i+1})|
&\leq cm\left(
    \frac{1}{(\lambda r)^2}
    +
    \frac{1}{h_N((\lambda r)\vee h_N)}
  \right)\\
&\leq c\|x-x'\|\left(
    \frac{1}{(\lambda r)^2}
    +
    \frac{1}{h_N((\lambda r)\vee h_N)}
  \right).
\end{split}
\end{align}
On the other hand, if $\min_{1\leq i\leq m}|y_i|< \lambda r/4$, this must mean there exists $i\in \{1,\dots,m\}$ such that $|\hat{z}_i-\hat{z}'_i|\geq \lambda r/4$, which implies that $h_N\geq c\lambda r$. Combining this with \eqref{eq:ggradDisc} and the telescoping series computation above again yields,
\[
|g_N(\bar{x})-g_N(\bar{x}')|\leq cm\left(
    \frac{1}{(\lambda r)^2}
    +
    \frac{1}{h_N^2}
  \right)
\leq c\|x-x'\|\left(
    \frac{1}{(\lambda r)^2}
    +
    \frac{1}{h_N((\lambda r)\vee h_N)}
  \right).
\]
It remains to determine how much the value of $g_N(\cdot)$ can vary within a single cable. Since $g_N(0,\cdot)$ is harmonic on $\tSlab\setminus\{0\}$ in the
sense of \cite[Definition 5.1]{prevost_2023_perco}, it is affine
on every cable that does not contain $0$. Now let $\underline{x}$ denote the other endpoint of the cable adjacent to $\bar x$ that contains $x$. By the harmonicity of $g_N(\cdot)$ and \eqref{eq:ggradDisc}, we have that,
\begin{equation}\label{eq:ggradPassCable}
    |g_N(x)-g_N(\bar{x})|\leq|g_N(\underline{x})-g_N(\bar{x})|
  \leq c\left(
    \frac{1}{(\lambda r)^2}
    +
    \frac{1}{h_N(\lambda r\vee h_N)}
  \right).
\end{equation}
The same argument holds for $x'$ and $\bar{x}'$. Combining this with \eqref{eq:gradTelescope} concludes the proof.
\end{proof}

\section{Lower tightness under \eqref{def:Con}}\label{sec:amplify}
The main result of this section is Proposition~\ref{lem:ConLb}, which supplies the required lower bound of $M_r$ in~\eqref{eq:tailEstlargestVolgena=0}.
Proposition~\ref{lem:ConLb} is obtained by strengthening the positive probability lower bounds in Section~\ref{sec:mainprop} to bounds arbitrarily close to $1$ through a decoupling argument.
As we explain below,~\eqref{def:Con} supplies the decoupling properties
needed to obtain such a lower bound.

Recall that under~\eqref{def:Con}, we write $K=\sup_N\capacity^{2r}(B_r)/\capacity(B_r)$.
\begin{proposition}\label{lem:ConLb}
  Under \eqref{def:Con}, for all $N\geq1$, $1\leq r\leq \Cr{range}N$ and $1\leq t\leq c\E_N[M_r]$,
  \begin{equation}\label{eq:ConLb}
    \P_N(M_r\geq t^{-1}\bv(r))\geq 1- ct^{-c'}
  \end{equation}
  where $c,c'$ depend only on $K$ and $\Cr{range}$.
\end{proposition}

The key to passing from the results in Section~\ref{sec:mainprop} to Proposition~\ref{lem:ConLb} is to use the domain Markov property~\eqref{eq:field_decomp} to obtain independent local fields, each of which contains a large cluster with positive probability. Controlling the accompanying harmonic averages via Lemma~\ref{lem:harm_avg} allows us to transfer these local estimates to the original field. As we explain shortly, the conditions needed for this control provide an indication of why \eqref{def:Con} arises as the natural condition for tightness.
In the sequel we write $t_+\defeq \max\{t,0\}$ for $t\in \R$.
\begin{lemma}\label{lem:harm_avg}
For all $\zeta\geq4$, $1\leq r\leq\Cr{range}N/\zeta$,
$a>0$ and $x\in\Slab$,
\begin{equation}\label{eq:harm_avg}
\P_N\Big(
\inf_{y\in\tB(x,r)}
\eta_y^{\partial B(x,\zeta r)}<-a
\Big)
\leq
\exp\!\bigg\{
-cF(\zeta r)
\Big(a-\frac{c'}{\sqrt{F(\zeta r)}}\Big)^2_+
\bigg\}.
\end{equation}
%Moreover, if $\zeta r\leq R\leq\Cr{range}N$ and $\tB(x,\zeta r)\subseteq\tB_R$, the same estimate holds with $\P_N$ and $F(\zeta r)$ replaced by $\P_N^R$ and $F^R(\zeta r)$, respectively.
\end{lemma}
We postpone the proof of Lemma~\ref{lem:harm_avg} to the end of the section and
now explain how~\eqref{def:Con} enters the argument of \eqref{eq:ConLb}.
Roughly speaking, we first choose a smaller scale $\rho$ at which
a local cluster in $\tB_\rho$ has the desired volume with positive
probability. We then apply \eqref{eq:field_decomp} in disjoint
balls of radius $\zeta\rho$, with $\zeta\geq4$, to obtain independent
local fields. To transfer their large clusters to the ambient field,
we require that the corresponding harmonic averages remain above
$-a$ throughout the smaller balls of radius $\rho$.
To retain a positive
probability of finding a large cluster in the local field at level
$a$, we can afford a shift of order $a\asymp F(\rho)^{-1/2}$. On the other
hand, \eqref{eq:harm_avg} shows that the probability of a bad
harmonic fluctuation is small when $a^2F(\zeta\rho)$ is large.
We thus need to choose $\zeta$ large enough that
$F(\zeta\rho)/F(\rho)$ is large, while keeping $\zeta\rho$
small enough that many disjoint balls of this radius fit inside
$\tB_r$.
In dimensions $d\geq3$, the polynomial growth of the corresponding
capacity scale allows the two requirements of $\zeta$ to be met simultaneously.
On $\Slab$, such a choice is only feasible under \eqref{def:Con}. This is evident in the following Lemma, where \eqref{eq:ConCrossing} shows that $\zeta$ is bounded in terms of $K$ whenever the capacity ratio is bounded by $K$, and the control on $\zeta$
deteriorates as $K$ increases.

\begin{lemma}\label{lem:loop_crossing_asymp}
Consider any sequence $r\!=\!r_N$ with $1\leq r_N\leq \Cr{range}N$.
Suppose \eqref{def:Con} holds and denote $K\defeq \sup_N \capacity^{2r}(B_r)/\capacity(B_r)$.
Then for any $\lambda\in(0,1)$, there exist $\Cl[c]{rlbK}(K,\lambda),\Cl[c]{rrplbK}\in[1,\infty)$ such that for all $\Cr{rlbK}\leq r\leq \Cr{range}N$,
\begin{equation}\label{eq:ConCrossing}
  \text{
$\exists 1\leq r'\leq r$ such that $F(r')\leq \lambda F(r)$ and $\frac{r}{r'}\leq \Cr{range}\exp\{\Cr{rrplbK} K/\lambda \}$.}
\end{equation}
Furthermore, there exists $\Cl[c]{lambdaKilledCap}\in(0,1)$ such that for all $1\leq r'\leq r\leq \Cr{range}N$ and $\lambda\in(0,\Cr{lambdaKilledCap})$,
\begin{equation}\label{eq:lambdaKilledCap}
  \text{$F(r')\leq \lambda F(r)$ implies  $\frac{\capacity^r(r')}{ \capacity(r)}\leq c\lambda $}
\end{equation}
\end{lemma}

The proof of Lemma~\ref{lem:loop_crossing_asymp} is presented in Section~\ref{sec:AFbox}.
We now explain how to obtain Proposition~\ref{lem:ConLb} from Lemma~\ref{pro:manyMesoKilled}.
The argument requires a delicate balance between the cost of decoupling and the independent contribution from local fields.
As we will see in~\eqref{eq:decouple_critical} and \eqref{eq:crossingMend}, this balance depends on the relative constants in the exponents, which we control using Lemma~\ref{lem:loop_crossing_asymp}. The resulting contribution from having bad harmonic averages decays only polynomially, whereas in the case $d\geq 3$ it would decay exponentially.

\begin{proof}[Proof of Proposition~\ref{lem:ConLb}]
Consider $m=m(K)\geq 1, \zeta=\zeta(K)\geq 1$ to be fixed soon and set $\rho=\frac{r}{\zeta m}$. Take $r\geq C(K)\geq 1$ so that $\rho\geq 1$.
In the sequel, we want to apply Lemma~\ref{lem:loop_crossing_asymp} when the capacity ratio is 
${\capacity^{2\zeta\rho}(B_{\zeta\rho})}/
{\capacity(B_{\zeta\rho})}$. We first show the following claim.
\begin{equation}\label{eq:capRatioClaim}
 \frac{\capacity^{2\zeta\rho}(B_{\zeta\rho})}
{\capacity(B_{\zeta\rho})}
\leq c\big(K+\log m\big). 
\end{equation}
Indeed, 
by \eqref{eq:capB} and \eqref{eq:cap-killed},
\[
\frac{\capacity^{2\zeta\rho}(B_{\zeta\rho})}
{\capacity(B_{\zeta\rho})}
\leq c\frac{\zeta\rho\wedge h_N}{F(\zeta\rho)}
\leq c\Big(
1+\frac{\zeta\rho\wedge h_N}{h_N}
\log^\star\!\Big(\frac{N}{\zeta\rho}\Big)
\Big)\leq c\Big(
1+\frac{r\wedge h_N}{h_N}\log^\star\!\Big(\frac Nr\Big)
+\log m
\Big),
\]
where the last inequality follows since $\zeta\rho\leq r$ and $\log^\star\!(\frac{N}{\zeta\rho})
=\log^\star\!(\frac{Nm}{r})
\leq \log^\star\!(\frac Nr)+\log m.$
It then follows from \eqref{def:Fbox} and \eqref{def:Con} that, $\frac{r\wedge h_N}{h_N}\log^\star\!(\frac Nr)
\leq \frac{r\wedge h_N}{F(r)}
\leq c\frac{\capacity^{2r}(B_r)}{\capacity(B_r)}
\leq cK.$
Combining this with the last display yields \eqref{eq:capRatioClaim}. It then follows from Lemma~\ref{lem:loop_crossing_asymp} that for any $\lambda\in(0,1)$,
\begin{equation}\label{eq:apply_crossing_lem}
  \text{$F(\rho)\leq\lambda F(\zeta\rho)$ whenever }
 \zeta\geq f(\lambda)\defeq 
\Cr{range}\exp\Big\{
\frac{c}{\lambda}
\big(K+\log(m)\big)\Big\}. 
\end{equation}

Now note that there exists an integer $n$ with $\frac12cm^2\leq n \leq cm^2$ and $(x_1,\dots,x_n)\subseteq B_r$ such that 
$B(x_i,\zeta\rho)\subseteq B_r$ are pairwise disjoint.
In view of \eqref{eq:field_decomp}, write
$\psi^i=\psi^{\partial B(x_i,\zeta\rho)}$ and
$\eta^i=\eta^{\partial B(x_i,\zeta\rho)}$.
Let $M_{x_i,\rho}^a(\psi^i)$ denote the largest cluster volume in
$\{\psi^i\geq a\}\cap\tB(x_i,\rho)$, as in~\eqref{def:Mr_general}
with $\tB(x_i,\rho)$ in place of $\tB_r$.
For $a\geq 0$ and $i\in\{1,\dots,n\}$, consider the events
$E_i^a=\{M_{x_i,\rho}^a(\psi^i)\geq c\bv(\rho)\}$.
We first claim that there exists $\Cl[c]{lambdaZeta}\in(0,1)$ such that for all $\lambda\in(0,\Cr{lambdaZeta})$, $\zeta\geq f(\lambda)$, $i\in\{1,\dots,n\}$ and $a\geq 0$,
\begin{equation}\label{eq:local_claim}
  \P_N(E_i^a) \geq c\exp\{-c' \lambda a^2F(\zeta\rho)\}.
\end{equation}
Indeed, Corollary~\ref{the:mainMra}
gives
$\P_N(M_\rho^{u}\geq c\bv(\rho))\geq c_0$ with $u=\frac{2c'}{\sqrt{F(\rho)}}$.
By \eqref{eq:harm_avg}, we have that
\[
\P_N\Big(\sup_{y\in \tB(x_i,\rho)}\eta^i_y>u\Big)\leq
 \exp\!\Big\{\!-\!cF(\zeta\rho)\Big(u-\frac{c'}{\sqrt{F(\zeta\rho)}}\Big)^2_+\Big\}\leq 
\exp\{-c\,F(\zeta\rho)/F(\rho)\}\leq c_0/2,
\]
where the last inequality holds by \eqref{eq:apply_crossing_lem} and upon choosing $\zeta\geq f(\Cr{lambdaZeta})$ with $\Cr{lambdaZeta}$ sufficiently small. It then follows from \eqref{eq:field_decomp} that $\P_N(E_i^0)\geq c_0/2$, which, combined with the
change of measure formula \cite[(3.7)]{drewitz_2024_clust} yields $\P_N(E_i^a)\geq c\exp\{-c'a^2F^{\zeta\rho}(\rho)\}$. We may now conclude the proof of \eqref{eq:local_claim} with \eqref{eq:lambdaKilledCap} upon decreasing $\Cr{lambdaZeta}$ if necessary to get $F^{\zeta\rho}(\rho)\leq c\lambda F(\zeta\rho)$.

Since $F(\cdot)$ is increasing and $\rho=r/(\zeta m)$,
\[
\frac{\bv(\rho)}{\bv(r)}
\geq c\Big(\frac{\rho}{r}\Big)^2
\frac{\rho\wedge h_N}{r\wedge h_N}
\sqrt{\frac{F(r)}{F(\rho)}}
\geq c(\zeta m)^{-3}.
\]
Thus, whenever $t\geq c(\zeta m)^3$ the cluster in the event
$E_i^a$ has volume greater than $t^{-1}\bv(r)$. It then follows from \eqref{eq:field_decomp}, \eqref{eq:local_claim} and \eqref{eq:harm_avg} that there exists $\Cl[c]{exp1},\Cl[c]{exp2}>0$ and $\Cl[c]{A}\in[1,\infty)$ such that
 for all $\lambda\in(0,\Cr{lambdaZeta})$, $\zeta\geq f(\lambda)$,
  $a\geq \Cr{A}/\sqrt{F(\zeta\rho)}$, $r\geq \zeta m$ and $t\geq c(\zeta m)^3$,
\begin{align}\label{eq:decouple_critical}
\begin{split}
\P_N\big(M_r\leq t^{-1}\bv(r)\big)
&\leq
\P_N\Big(\bigcap_{i=1}^n (E_i^a)^c\Big)
+\sum_{i=1}^n
\P_N\Big(\inf_{y\in \tB(x_i,\rho)}\eta^i_y<-a\Big)\\
&\leq \big(1-c\exp\{-c \lambda a^2F(\zeta\rho)\}\big)^n
+c'm^2\exp\{-ca^2F(\zeta\rho)\}\\
&\leq \exp\Big\{\!-\!cm^2\exp\big\{-\Cr{exp1}\lambda a^2F(\zeta\rho)\big\}\Big\}+c'm^2\exp\{-\Cr{exp2}a^2F(\zeta\rho)\}.
\end{split}
\end{align}
where we also used the bound $(1-x)\leq e^{-x}$ to get the last inequality.

We now fix $a$ and $\lambda$. 
We pick $a^2=A\log(m)/F(\zeta\rho)$ with $A=3\Cr{exp2}^{-1}\vee \Cr{A}$ and choose $\lambda=\Cl[c]{endlambda}\in(0,1)$ sufficiently
small such that $2-\Cr{exp1}\lambda A\geq 1$. These choices, combined with the last display, give that for all $\zeta\geq f(\Cr{endlambda})$, $r\geq \zeta m$ and $t\geq c(\zeta m)^3$,
\begin{equation}\label{eq:crossingMend}
  \P_N\big(M_r\leq t^{-1}\bv(r)\big)
\leq \exp\{-c m^{2-\Cr{exp1}\lambda A}\}+c'm^{2-\Cr{exp2}A}
\leq  \exp\{-cm\}+c'm^{-1}
\leq c''m^{-1}.
\end{equation}
It remains to express $m$ in terms of $t$. In view of \eqref{eq:apply_crossing_lem}, with the above choice of $\lambda$, there exists $\beta\geq 1$ such that $m^\beta\geq f(\Cr{endlambda})$
for all sufficiently large $m$ depending on $K$. Hence in view of \eqref{eq:crossingMend}, taking $\zeta=m^\beta$ and using the bound $t\geq c(\zeta m)^{3}$ yields that for all $r\geq C(K),C'(K)\leq t\leq c\bv(r)$,
\[
\P_N\big(M_r\leq t^{-1}\bv(r)\big)
\leq ct^{-1/(3(\beta+1))},
\]
which concludes the proof in this case.
The proof when either $t\leq C'(K)$ or $r\leq C(K)$ simply follows from Corollary~\ref{the:mainMra}.
\end{proof}
We finish with the
\begin{proof}[Proof of Lemma~\ref{lem:harm_avg}]
By \eqref{eq:field_decomp}, \eqref{eq:gDecomp} and
\eqref{eq:g2-asymp-unif}, we have that for all $y\in\tB(x,r)$,
\[
\E_N\big[(\eta_y^{\partial B(x,\lambda r)})^2\big]
=
E_y\Big[
g_N\big(y,X_{H_{\partial B(x,\lambda r)}}\big),
H_{\partial B(x,\lambda r)}<\infty
\Big]
\leq \frac{c}{F(\lambda r)}.
\]
We now claim that
\[
\E_N\Big[
\sup_{y\in\tB(x,r)}\eta_y^{\partial B(x,\lambda r)}
\Big]\leq c\Big(\sup_{y\in\tB(x,r)}\E_N\Big[
(\eta_y^{\partial B(x,\lambda r)})^2
\Big]\Big)^{1/2}
\leq \frac{c}{\sqrt{F(\lambda r)}}.
\]
This follows by adapting the proof of
\cite[Lemma~8.2]{rodriguez_2025_anoma} to $\Slab$. The
argument used in
\cite[Lemma~8.4]{rodriguez_2025_anoma} still applies using the Harnack inequality derived in \cite[Corollary 5.2]{rodriguez_2025_on}, which holds on $\Slab$.
The subsequent entropy estimate uses a covering bound
of order $(r'/\rho)^3$ for balls in $\Slab$, in place
of the two-dimensional bound below \cite[(8.16)]{rodriguez_2025_anoma}; this changes only the
constant. 
We may now conclude since \eqref{eq:harm_avg} follows from the last two displays and the BTIS inequality (see~\cite[Theorem 2.1.1]{_2007_rando}). 
%The proof under $\P_N^R$ follows with the same arguement.
\end{proof}

\section{Limiting supercritical density and local uniqueness}\label{sec:Iu}
This section takes the first step towards Theorem~\ref{thm:limit}
by identifying the asymptotic density of the largest cluster in
the killed field at a supercritical level.
Proposition~\ref{pro:convergences} shows that $M_r^{-a}$,
suitably normalised, converges in probability to one.
The density appearing in the normalisation is the same function that
appears in the limiting law of Theorem~\ref{thm:limit}.

\begin{proposition}\label{pro:convergences}
If $1\leq r=r_N\leq \Cr{range}N$ satisfies \eqref{def:AntiCon} and $\lim_N r=\infty$, we have that for all $a\geq c(F(r))^{-\frac12}$ and $\epsilon\in(0,1)$
\begin{equation}\label{eq:pro_cvg}
    \lim_{N}   \P_N^{2r}\bigg(\bigg|\frac{M^{-a}_r}{|B_r|\cdot\big(2\Phi\big(({a^2}/{{g^{2r}_N}})^{1/2} \big)-1\big)}-1\bigg|> \epsilon\bigg)=0 .
\end{equation}
\end{proposition}

The proof of Proposition~\ref{pro:convergences} draws on two results of independent interest.
It relies on various consequences of the isomorphism between the field and random interlacement, which we collect in Lemma~\ref{lem:isom_consequence}. It also requires a local uniqueness estimate for the surviving trajectories in the interlacement set; see Proposition~\ref{lem:LU}.
We combine these results in
Section~\ref{sec:proofWithIu} to prove Proposition~\ref{pro:convergences}.

\smallskip 
Before diving into the isomorphism theorem,
we need the following auxiliary Lemma.
In the sequel, for $t\geq0$, we abbreviate
\begin{equation}\label{def:vartheta}
  \vartheta(t)\defeq 2\Phi(\sqrt{t})-1.
\end{equation}
The notation choice $\vartheta$ reflects the role of
$\vartheta(a^2/g_N^{2r})$ as the asymptotic density of the
largest level $-a$ cluster in $B_r$ under $\P_N^{2r}$,
evident in~\eqref{eq:pro_cvg}.
\begin{lemma}
For all $s,t\geq0$,
\begin{align}
\label{eq:vartheta_bound}
  &\vartheta(t)\asymp (\sqrt{t}\wedge1)\\
  \label{eq:vartheta_scaling}
  &(s\wedge1)\vartheta(t)
  \leq\vartheta(s^2t)
  \leq(s\vee1)\vartheta(t).
\end{align}
\end{lemma}
\begin{proof}
For~\eqref{eq:vartheta_bound}, $\vartheta(t)\asymp \int_0^{\sqrt{t}}e^{-x^2/2}\,dx\asymp1$ when $t\geq 1$ and $\vartheta(t)\asymp \int_0^{\sqrt{t}}1\,dx\asymp\sqrt{t}$ when $t\in[0,1)$, completing the proof of~\eqref{eq:vartheta_bound}.
We now move on to showing~\eqref{eq:vartheta_scaling}.
The function $v\mapsto\vartheta(v^2)$ is increasing
and concave on $[0,\infty)$.
For $0\leq s\leq1$, the upper bound follows by monotonicity.
Using Jensen's inequality and the fact that $\vartheta(0)=0$, we have $\vartheta\big(\big(s\sqrt{t}+(1-s)0\big)^2\big)\geq s\vartheta(t)$, which yields the lower bound. When $s>1$, the lower bound follows from monotonicity and Jensen's inequality yields
\[
\vartheta(t)=\vartheta((s^{-1} s\sqrt{t} +(1-s^{-1})0)^2)\geq s^{-1}\vartheta(s^2t),
\]
which completes the proof.
\end{proof}

\smallskip

We now introduce the objects entering the isomorphism theorem in the form of~\cite[(Isom')]{drewitz_2022_clust}. Let $\bar\P_N^{r}$ denote the law of random interlacements associated with the diffusion with additional killing upon hitting $\partial B_r$, and consider, independently, a field $\varphi$ with law $\P_N^{r}$. Note that these objects can equivalently be defined
on the cable system $\ttB_r$ associated with the weighted graph
$(\overline B_r,\lambda^{r},\kappa^{r})$, where $\lambda_{x,x'}^r=\lambda_{x',x}$ for $x,x'\in \overline{B}_r$ and $\kappa_x=  N^{-2}$ for all $x\in B_r$ and $\kappa_x= \infty$ for $x\in \partial B_r$.

For $u>0$ and $x\in \ttB_r$, we write $\ell_x^u$ for local time at site $x$ of the interlacement process at level $u$. We denote by $\Sigma_x$ the connected component of $x$ in $\{y\in\ttB_r:\varphi_y\neq0\}$.
On an extension $\tilde{\P}_N^{r}$ of $\P_N^{r}\otimes\bar\P_N^{r}$, conditionally on $(|\varphi_x|)_{x\in \ttB_r}$ and $(\ell_x^u)_{x\in \ttB_r,u\geq 0}$, let $\sigma^u$ be constant on each connected component of $\{x\in\ttB_r:\varphi_x^2+2\ell_x^u>0\}$, equal to $1$ on those intersecting $\I^u$, and take independent, uniformly distributed values in $\{-1,1\}$ on the remaining components. Its values on the zero set of $\varphi$ may be chosen arbitrarily. Then \cite[(Isom')]{drewitz_2022_clust} gives, for each $u>0$,
\begin{equation}\label{eq:isomIu}
  \big(\sigma_x^u\sqrt{\varphi_x^2+2\ell_x^u}\big)_{x\in\ttB_r}
  \stackrel{d}{=}
  \big(\varphi_x+\sqrt{2u}\big)_{x\in\ttB_r},
\end{equation}
where the laws on the left- and right-hand sides are $\tilde{\P}_N^{r}$ and $\P_N^{r}$, respectively.

\medskip
We now collect in Lemma~\ref{lem:isom_consequence} a few useful consequences of \eqref{eq:isomIu}. In particular,~\eqref{eq:wired_sgn_conn} below can be viewed as an analogue of the celebrated formula~\cite[Proposition 5.2]{lupu_2016_from} when the level is nonzero. For convenience, we state the results under the measure $\tilde{\P}^r_N$, but Lemma~\ref{lem:isom_consequence} also holds on general graphs where \eqref{eq:isomIu} is valid, see Remark~\ref{rem:isomIu}, (1). To formulate the result neatly, we consider the graph $\ttB_{r,\mathrm{wired}}$, which is obtained from $\ttB_r$ by joining all open endpoints in $I_x$, $x\in \overline{B}_r$ to a single vertex $\partial$, and set $\varphi_\partial=0$.
For $a\in\R$ and $A,B\subseteq\ttB_r$, we write
$A\xleftrightarrow[\mathrm{wired}]{a}B$ if $A$ and $B$ are connected
by a continuous path in
$\{z\in\ttB_{r,\mathrm{wired}}:\varphi_z\neq a\}$.

\begin{lemma}\label{lem:isom_consequence}
For all $x,y\in B_r$ and $u\geq 0$,
\begin{align}
  \label{eq:wired_sgn_conn}
  &\tilde{\E}^{r}_N\big[\sigma_x^{u}\sigma_y^{u}\big]=
  \P_N^{r}\big(x\xleftrightarrow[\text{wired}]{\sqrt{2u}} y\big)
 %\asymp \Big(\frac{g^r_N(x,y)}{\sqrt{g^r_N(x,x)g^r_N(y,y)}} + \frac{2u}{\sqrt{g_N^r(x,x)g_N^r(y,y)}}\Big)\wedge1
 \\
  \label{eq:sign_int_identity}
    &\tilde{\E}^{r}_N[\sigma_x^u]=\tilde{\P}^{r}_N(\Sigma_x\cap \I^u\neq \emptyset)= 2\P_N^{r}(\varphi_x\leq \sqrt{2u})-1\asymp \frac{\sqrt{2u}}{\sqrt{g^r_N(x,x)}}\wedge 1\\
  \label{eq:cov_sigma}
  &\text{Cov}^{\tilde{\P}_N^{r}}(\1_{\{\Sigma_x\cap \I^u\neq \emptyset\}},\1_{\{\Sigma_y\cap \I^u\neq \emptyset\}})
  \leq \text{Cov}^{\tilde{\P}_N^{r}}(\sigma_x^u,\sigma_y^u)\leq  \frac{g^{r}_N(x,y)}{\sqrt{g^{r}_N(x,x)g^{r}_N(y,y)}}.
\end{align}
\end{lemma}

\begin{proof}
  We start with~\eqref{eq:sign_int_identity}. 
  By the definition of $\sigma_x^u$, on the event $\{\Sigma_x\cap \I^u=\emptyset\}$, the expected value of $\sigma_x^u$ is $0$, hence
  \[
  \tilde{\E}^{r}_N[\sigma_x^u]=\tilde{\E}^{r}_N[\sigma_x^u\1\{\Sigma_x\cap \I^u\neq \emptyset\}]+\tilde{\E}^{r}_N[\sigma_x^u\1\{\Sigma_x\cap \I^u=\emptyset\}]=\tilde{\P}^{r}_N(\Sigma_x\cap \I^u\neq \emptyset).
  \]
  Furthermore, by~\eqref{eq:isomIu},
  \begin{equation}\label{eq:sigma_ind}
      \sigma_x^u \text{ under } \tilde{\P}_N^r
  \quad\text{has the same law as}\quad
  2\1\{\varphi_x\geq -\sqrt{2u}\}-1 \text{ under } \P_N^r,
  \end{equation}
  which concludes the proof of~\eqref{eq:sign_int_identity} using \eqref{eq:vartheta_bound}.
  For \eqref{eq:wired_sgn_conn}, using again~\eqref{eq:isomIu}, we have
  \[
  \tilde{\E}^{r}_N[\sigma_x^u \sigma_y^u]=
  \tilde{\P}^{r}_N(\Sigma_x\cap \I^u\neq \emptyset,\Sigma_y\cap \I^u\neq \emptyset)
  + \tilde{\P}^{r}_N\big(\Sigma_x\cap \I^u = \emptyset,\Sigma_y\cap \I^u= \emptyset, y\in \Sigma_x\big),
  \]
  and note that the r.h.s is exactly $\P_N^{r}(x\xleftrightarrow[\text{wired}]{\sqrt{2u}} y)$. 
For~\eqref{eq:cov_sigma}, the last display and~\eqref{eq:sign_int_identity} implies that $\text{Cov}^{\tilde{\P}_N^{r}}(\1_{\{\Sigma_x\cap \I^u\neq \emptyset\}},\1_{\{\Sigma_y\cap \I^u\neq \emptyset\}})\leq \text{Cov}^{\tilde{\P}_N^{r}}(\sigma_x^u,\sigma_y^u)$. 
  Using~\eqref{eq:sigma_ind},
  \[
  \text{Cov}^{\tilde{\P}_N^{r}}\!(\sigma_x^u,\sigma_y^u)\!=\!\text{Cov}^{{\P}_N^{r}}(2\1_{\{\varphi_x\geq -\sqrt{2u}\}}-1,2\1_{\{\varphi_y\geq -\sqrt{2u}\}}-1)
  \!=\!4\text{Cov}^{{\P}_N^{r}}(\1_{\{\varphi_x\geq -\sqrt{2u}\}},\1_{\{\varphi_y\geq -\sqrt{2u}\}}).
  \]
Write $\rho\defeq g_N^{r}(x,y)/\sqrt{g_N^{r}(x,x)g_N^{r}(y,y)}$ and $b_z\defeq\sqrt{2u/g_N^{r}(z,z)}$, $z\in\{x,y\}$. Let $f(x,y,\rho)$ be the density of a bivariate standard normal distribution with correlation $\rho$.
By Price's theorem, for any $s,t\in R$ $\frac{d}{d\rho}\int_{s}^\infty\int_{t}^\infty f(x,y,\rho)\,dx\,dy$=$\int_{s}^\infty\int_{t}^\infty  \frac{d^2}{dxdy} f(x,y,\rho)\,dx\,dy=f(t,s,\rho)$. It then follows from this that
\[
\text{Cov}^{\tilde{\P}_N^{r}}\!(\sigma_x^u,\sigma_y^u)
\!=\!4\int_0^\rho \!f(-b_x,-b_y,t) \,dt
\!=\!\frac{2}{\pi}\!\int_0^\rho\!
\exp\!\left\{-\frac{b_x^2-2t b_xb_y+b_y^2}{2(1-t^2)}\right\}
\frac{dt}{\sqrt{1-t^2}}
\!\leq\!\frac{2}{\pi}\arcsin\rho\leq\rho.
\]
\end{proof}

\medskip
We next turn to Proposition~\ref{lem:LU}, a local uniqueness estimate formulated in terms of \emph{forward-surviving} trajectories, which we now define.
We call a doubly-infinite trajectory on $\ttB_r$ \emph{forward surviving} $(+)$ if its forward part ends through
the outer boundary $\partial B_r$, and \emph{surviving} $(++)$ if both its forward and backward parts end through $\partial B_r$.
See~\cite[(7.4)]{drewitz_2023_criti} and \cite[Section 2]{prevost_2023_perco} for similar considerations. 
The surviving condition precludes short trajectories that 
end ``early'' upon exiting through one of the half-open cables $I_x$, $x\in B_r$ before reaching $\partial B_r$.
Let $W_+^*$ and $W_{++}^*$ denote the corresponding classes of trajectories modulo time-shift. Write $\nu_N^{r}$ for the intensity measure of $\I^u$ under $\P_N^{r}$ and set
\begin{equation}\label{eq:forward_interlacement_measures}
  \nu_{N,+}^{r}\defeq\1_{W_+^*}\nu_N^{r},
  \qquad
  \nu_{N,++}^{r}\defeq\1_{W_{++}^*}\nu_N^{r}.
\end{equation}
For $u>0$, the restrictions to these classes are Poisson point processes with respective intensities $u\nu_{N,+}^{r}$ and $u\nu_{N,++}^{r}$. We denote the unions of their trajectory ranges by $\I_+^u$ and $\I_{++}^u$, respectively, so that $\I_{++}^u\subseteq\I_{+}^u$.

For $x\in B_r$, $1\leq s\leq t$, consider the event
\begin{equation} \label{def:local_uniqueness}
   \LU_{u,s,t}^+(x)\defeq \bigcap_{z,y\in\I^u_+\cap B(x,s)} 
   \big\{z \leftrightarrow y \text{ in } \I^u_+\cap \tB(x,t)\big\}.
\end{equation}

\begin{proposition}[Local Uniqueness]\label{lem:LU}
  For all $u>0$ $\lambda\in(0,\frac12),x\in B_r$ and $r\leq \Cr{range}N$ such that $\lambda r\geq 10$ and $B(x,2\lambda r)\subseteq B_r$, we have
  \begin{equation}\label{eq:LUsingle}
    \bar{\P}_N^{2r}\big( \LU_{u,\lambda r, 2\lambda r}^+(x)^c \big)\leq 
    cuF^{2r}\!(\lambda r)\,\exp\Big\{ \!-\!c'\lambda \big(u F^{2r}(\lambda r)\big)^{\frac13}\Big\}.
  \end{equation}
  Furthermore,
    \begin{equation}\label{eq:LUchain}
    \bar{\P}_N^{2r}\big( \LU_{u,(1-2\lambda)r,r}^+(0)^c \big)\leq 
     c\lambda^{-3} uF^{2r}\!(\lambda r)\,\exp\big\{ \!-\!c'\lambda \big(u F^{2r}(\lambda r)\big)^\frac13 \big\}.
  \end{equation}
\end{proposition}
Results of this kind has been derived on $\Z^d,d\geq 3$ for the full interlacement process at fixed intensity in~\cite[Proposition~1]{rath_2011_on}.
\cite[Theorem~5.1]{drewitz_2023_criti} obtained near-critical
stretched-exponential bounds for the full interlacement process on a
broad family of transient graphs, with exponent $1/3$ on $\Z^3$.
Sharper bounds without the $1/3$ exponent were subsequently obtained
in~\cite[Theorem~1.1]{prevost_2025_first}.  We believe that the $1/3$
exponent in Proposition~\ref{lem:LU} can also be improved, though we do not pursue this here.
Our proof adapts the argument of~\cite[Theorem~5.1]{drewitz_2023_criti}. In our massive setting, the full interlacement process contains short
trajectories that are killed before exploring a substantial region and are therefore difficult to connect to uniformly.
We hence restrict to forward-surviving
trajectories; see Remark~\ref{rem:isomIu},~(3), for more discussions on this.

\smallskip

As a key input, we first show that, for every $K\subseteq B_r$, the intensity of
forward-surviving trajectories hitting $K$ is comparable to that
of all trajectories hitting $K$.
For any $K\subseteq \ttB_r$, let $W_K$ denote the class of doubly-infinite trajectories on $\ttB_r$ hitting $K$ modulo time-shift.
It is worth noting that the result below does not hold for general $K\subseteq \ttB_r$. See Remark~\ref{rem:isomIu}, (2).
\begin{lemma}
  For all $1\leq r\leq \Cr{range}N$ and $K\subseteq B_r$,
\begin{equation}\label{eq:surv_intensity}
 c\,\capacity^{r}(K)\leq \nu_{N,+}^{r}(W_K)\leq  \capacity^{r}(K)
\end{equation}
%\quad \text{and}\quad  \nu_{N,-}^{r}(W_K)\leq c (\frac{r^2}{N^2}\wedge 1)\capacity_N^{r}(K)
\end{lemma}
\begin{proof}
 For every $x\in B_r$, the inclusion $B_r\subseteq B(x,2r)$, translation invariance and \eqref{eq:hitting_lb} give
$P_x(H_{\partial B_r}<\infty)\geq P_0(H_{\partial B_{2r}}<\infty)\geq c.$
 By the definition of $\nu_{N,+}^{r}$, see for instance \cite[(2.11)]{prevost_2023_perco}, we therefore have
 \[
 \nu_{N,+}^{r}(W_K)=\sum_{x\in K}e_{K}^{r}(x)P_x(H_{\partial B_r}<\infty)\geq c\sum_{x\in K}e_{K}^{r}(x)=c\,\capacity^{r}(K).
 \]
 The upper bound follows similarly by the trivial bound $P_x(H_{\partial B_r}<\infty)\leq 1$.

\end{proof}

The following is a counterpart of~\eqref{eq:sign_int_identity} for the surviving interlacement, which is obtained through an $\mathbf{h}$-transformed version of~\eqref{eq:isomIu}, see~\cite[Theorem 5.5]{prevost_2023_perco}.
\begin{lemma}\label{lem:h-isom}
For all $x\in B_r$,
\begin{equation}\label{eq:h-isom}
  \tilde{\P}_N^{r}(\Sigma_x \cap \I^u_{+}\neq \emptyset)\geq
  \tilde{\P}_N^{r}(\Sigma_x \cap \I^u_{++}\neq \emptyset)
  = 2\P_N^{r}\Big(\varphi_x\leq \sqrt{2u}P_x(H_{\partial B_r}<\infty)\Big)-1.
\end{equation}
\end{lemma}
\begin{proof}
   For any $y\in \tSlab$, let $\mathbf{h}(y)=P_y(H_{\partial B_r}<H_{\Delta})$, where here hitting the cemetery state $\Delta$ means that $X$ reaches the open end of one of the killing cables $I_x$ where $x\in B_r$.
  Note that $\mathbf{h}$ is harmonic in $\ttB_r$.
 By an argument analogous to the proof of
\cite[Corollary 5.9]{prevost_2023_perco}, interlacements on the
$\mathbf{h}$-transform of $\ttB_r$, after the corresponding change of space and time,
have the same law as the surviving interlacement process
under $\bar{\P}_N^{r}$.
  It then follows from \cite[Theorem 5.5]{prevost_2023_perco} that
  \[
  \tilde{\P}^{r}_N(\Sigma_x \cap \I^u_{++}\neq \emptyset)=2
  \P_N^{r}\Big(\varphi_x\leq \sqrt{2u}h(x)\Big)-1.
  \]
  This concludes the proof since $h(x)=P_x(H_{\partial B_r}<\infty)$ for $x\in \Slab$.
\end{proof}

We also need a connection estimate for $\I^u_+$, analogous to
\cite[Lemma~4.3]{drewitz_2025_geome}.

\begin{lemma}\label{lem:il_connect_two}
For all $u>0$, $\lambda\in(0,\frac12)$ and $r\leq\Cr{range}N$
such that $\lambda r\geq1$, $x\in\Slab$ such that
$B(x,2\lambda r)\subseteq B_r$, and all 
$U,V\subseteq B(x,\lambda r)$,
\begin{equation}\label{eq:il_connect_two}
  \bar\P_N^r\left(
    U\leftrightarrow V
    \text{ in }\I_+^u\cap\tB(x,2\lambda r)
  \right)
  \geq
  1-\exp\left\{
    -\frac{cu\,\capacity^r(U)\capacity^r(V)}{\lambda r\wedge h_N}
  \right\}.
\end{equation}
\end{lemma}

\begin{proof}
We have seen in the proof of~\eqref{eq:surv_intensity} that
$\inf_{y\in B_r}P_y(H_{\partial B_r}<\infty)\geq c$.
Thus, applying the strong Markov property at time $H_V$ yields that for every $y\in U$,
\begin{align*}
  P_y(H_V<H_{B(x,2\lambda r)^c}|H_{\partial B_r}<\infty)
  &\geq cP_y(H_V<H_{B(x,2\lambda r)^c})\\
  &\geq
  c\,\capacity^{B(x,2\lambda r)^c}(V)
  \inf_{z\in V}g_N^{B(x,2\lambda r)^c}(y,z)
  \geq
  \frac{c\,\capacity^r(V)}{\lambda r\wedge h_N},
\end{align*}
where the second inequality follows from~\eqref{eq:last_exit},
and the last follows from~\eqref{eq:killed_green_smallh} and \eqref{eq:CapVarForm}.

Conditionally on its first entrance point in $U$ being $y\in U$, the forward
part of each trajectory in $\I^u_+$ has law $P_y(\cdot|H_{\partial B_r}<\infty)$.
The last display and Poisson thinning therefore give
\[
  \bar\P_N^r\left(
    U\not\leftrightarrow V
    \text{ in }\I_+^u\cap\tB(x,2\lambda r)
  \right)
  \leq
  \exp\left\{
    -u\nu_{N,+}^r(W_U)
    \frac{c\,\capacity^r(V)}{\lambda r\wedge h_N}
  \right\}.
\]
Applying~\eqref{eq:surv_intensity} yields
\eqref{eq:il_connect_two}.
\end{proof}

\begin{proof}[Proof of Proposition~\ref{lem:LU} (sketch)]
We first explain how~\eqref{eq:LUchain} follows from~\eqref{eq:LUsingle}.
For $n\geq 1$ and $\lambda \in(0,\frac12)$ consider a collection of vertices $\{x_1,\dots,x_n\}\subset B_r$ such that
\begin{align}\label{eq:balls_LU}
\begin{split}
  &n\leq c\lambda^{-3}, \quad B_{(1-2\lambda)r}\subseteq \bigcup_{i=i}^nB(x_i,\tfrac{1}{10}\lambda r),\quad
  \bigcup_{i=1}^nB(x_i,2\lambda r)\subseteq B_r\\
  &B(x_n,\tfrac{1}{10}\lambda r)\subseteq B(x_1,\lambda r)
  \quad \text{and}\quad 
   B(x_i,\tfrac{1}{10}\lambda r)\subseteq B(x_i+1,\lambda r)\text{ for all $1\leq i\leq n-1$}.
\end{split}
\end{align}
By~\eqref{eq:balls_LU}, $\I^u_+\cap B_{(1-2\lambda )r}\subseteq \bigcup_{i=1}^n (\I^u_+\cap B(x_i,\frac{1}{10}\lambda r))$. Hence, on the event $\bigcap_{i=1}^n \big\{\LU^+_{u,\lambda r,2\lambda r}(x_i)\big\}\cap \big\{\I^u_+\cap B(x_i,\tfrac{1}{10}\lambda r)\neq\emptyset\big\}$,
every vertex in $\I^u_+\cap B_{(1-2\lambda )r}$ lies in a single connected component in $\I^u_+\cap \bigcup_{i=1}^n  B(x_i,2\lambda r)$. It then follows from~\eqref{eq:surv_intensity},~\eqref{eq:balls_LU} and~\eqref{eq:LUsingle} that,
\[
\bar{\P}_N^{2r}\big( \LU_{u,(1-2\lambda)r,r}^+(0)^c \big)
\leq c\lambda^{-3}\,\Big(uF^{2r}\!(\lambda r)
  \exp\big\{ \!-\!c'\lambda \big(u F^{2r}(\lambda r)\big)^{\frac13}\big\}
  +\exp\big\{ \!-\!c'u F^{2r}(\lambda r)\big\}
\Big),
\]
which concludes the proof of~\eqref{eq:LUchain}.

We now sketch the proof of~\eqref{eq:LUsingle}, which is obtained by adapting the proof of~\cite[Theorem 5.1]{drewitz_2023_criti}. We follow the same decomposition of $\I^u_+$ and write $\I^u_+=\bigcup_{k=1}^4\I_+^{u/4,k}$, where the four interlacement sets are independent and each has law $\bar\P_N^{2r}$ at level $u/4$. For $k\in\{1,\dots,4\}$, let $Z_k^1,\dots,Z_k^{N_k}$ be the trajectories underlying $\I_+^{u/4,k}$ that hit $B(x,\lambda r)$. In view of~\eqref{eq:surv_intensity}, $N_k$ is still Poisson with parameter $\asymp u\capacity^{2r}(B(x,\lambda r))$.
Starting at its first entrance into $B(x,\lambda r)$, decompose each $Z_k^i$ into its $M_k^i$ successive excursions $Z_k^{i,1},\dots,Z_k^{i,M_k^i}$, each starting upon entering $B(x,\lambda r)$ and stopped at its first exit from $\tB(x,2\lambda r)$. Now let $A^u_r$ denote the event analogous to the one defined below \cite[(5.14)]{drewitz_2023_criti}, which asks for $N_k\leq cuF^{2r}(\lambda r)$ for all $k\in\{1,2,3,4\}$ and $M^k_i\leq cuF^{2r}(\lambda r)$ for all $i\in\{1,\dots ,N_k\}$.

Proceeding as in the argument below \cite[(5.16)]{drewitz_2023_criti}, one can obtain that,
\begin{equation}\label{eq:LUA}
  \bar{\P}_N^{2r}\big(\LU_{u,\lambda r, 2\lambda r}^+(x)^c\cap A^u_r\big)\leq c(uF^{2r}(\lambda r))^4\exp\big\{-c'(uF^{2r}(\lambda r))^{\frac13}\big\}.
\end{equation}
Indeed, the key inputs to showing \eqref{eq:LUA} in \cite{drewitz_2023_criti} are \cite[Lemma 4.3]{drewitz_2025_geome} and~\cite[Lemma 5.3]{drewitz_2023_criti}.
Since the trajectories under consideration are all forward-surviving, ~\eqref{lem:cap_walk_space_2} and \eqref{eq:il_connect_two} provide the analogue of \cite[Lemma 5.3]{drewitz_2023_criti} and \cite[Lemma 4.3]{drewitz_2025_geome}.
It then remains to show that
\begin{equation}\label{eq:LUAc}
  \bar{\P}_N^{2r}\big( (A^u_r)^c\big)\leq cuF^{2r}\!(\lambda r)\times\exp\big\{-c'\lambda uF^{2r}(\lambda r)\big\}.
\end{equation}
Following the geometric domination argument below \cite[(5.13)]{drewitz_2023_criti}, \eqref{eq:LUAc} reduces to showing
\begin{equation}\label{eq:gamblersBound}
\inf_{y\in \partial B(x,{2\lambda r})}P_{y}(H_{B_{2r}^c}<H_{B(x,{\lambda r})})\geq c\lambda. 
\end{equation}
Let $f(\cdot)\defeq P_{\cdot}(H_{B_{3r}^c}<H_{B_{\lambda r}})$.
Since $B_{2r}\subseteq B(x,3r)$, it suffices to show that for all $y\in \partial B_{2\lambda r}$, $f(y)\geq c(\lambda)$.
For all $y\in \partial B_{2\lambda r}$, we can find $(z_1,\dots,z_k)\subset \partial B_{2\lambda r}$ such that $1\leq k\leq c$, $B(z_i,\lambda r/4)\cap B(z_{i+1},\lambda r/4)\neq\emptyset$ for all $1\leq i\leq k-1$ and $y\in B(z_1,\lambda r/4)$ and $\lambda r e_1\in B(z_k,\lambda r/4)$. Furthermore, since $f$ is harmonic on $B(z,\lambda r/2)$ for all $z\in \partial B_{2\lambda r}$, we may apply the Harnack inequality on $\Slab$ (see~\cite[Corollary 5.2]{rodriguez_2025_on}) $k$
times to get, 
\begin{equation}\label{eq:applyHarnack}
  f(y)\geq c f(2\lambda r e_1), \quad \text{for all $y\in \partial B_{2\lambda r}$.}
\end{equation}
Now let $(X^1_t)_{t\geq 0}$ be the first coordinate of $(X_t)_{t\geq 0}$. Hence $X^1$ has the same law as a continuous-time simple random walk on $\Z$ with jump rate $\frac13$ killed with rate $N^{-2}$. Let us denote the law of $X_1$ starting at $z\in \Z$ by $\P^1_z$. It's easy to see that
\begin{equation}
  f(2\lambda r e_1)\geq P^1_{2\lambda r}(H_{\{3 r\}}< H_{\{\lambda r\}})\geq \frac{\lambda r}{(3-\lambda )r}e^{-cr^2/N^2}\geq c\lambda,
\end{equation}
which, combined with~\eqref{eq:applyHarnack}, concludes the proof of~\eqref{eq:gamblersBound}, and hence also the proof of~\eqref{eq:LUAc} and~\eqref{eq:LUsingle}.
\end{proof}

\Needspace{4\baselineskip}
\begin{remark}\label{rem:isomIu}\leavevmode
  \begin{enumerate}[label={\arabic*)}]
    \item  The proof of Lemma~\ref{lem:isom_consequence} relies solely on~\eqref{eq:isomIu}. Hence, in view of \cite[Theorem 1.1 and Lemma 3.4]{drewitz_2022_clust}, Lemma~\ref{lem:isom_consequence} holds on any weighted graph where the corresponding on-diagonal Green's function is uniformly bounded.
  If, moreover, $\kappa\equiv0$ and the interlacement set is almost
  surely connected for every $u>0$ (e.g. on $\Z^d$ with $d\geq 3$), \eqref{eq:wired_sgn_conn} holds without wiring, i.e. we have $\E[\sigma_x^u\sigma_y^u]=\P(\{x\stackrel{\geq\sqrt{2u}}{\longleftrightarrow}y\}\cup\{x\stackrel{\leq \sqrt{2u}}{\longleftrightarrow}y\})$.
  \item 
  For general $K\subseteq\ttB_r$, the exit probabilities appearing
  in the proof of~\eqref{eq:surv_intensity} can be arbitrarily small
  when $K$ is close to the open end of a killing cable $I_x$,
  $x\in B_r$. Thus, the lower bound in~\eqref{eq:surv_intensity}
  need not hold for such sets.
  In particular, we cannot apply~\eqref{eq:surv_intensity} directly
  to $\Sigma_x$ to control the probability of
  $\{\Sigma_x\cap\I_+^u\neq\emptyset\}$.
  This is why, in Lemma~\ref{lem:h-isom}, we instead pass to the
  surviving interlacement process $\I_{++}^u$ and invoke the
  $\mathbf h$-transform version of \eqref{eq:isomIu}.
  \item Without the forward-surviving restriction in Proposition~\ref{lem:LU}, the analogue
  of~\eqref{eq:rw_cap_with_killing} would acquire an additional
  error term of the form $C(r')^2/(sN^2)$, accounting for trajectories
  that are killed before accumulating sufficient capacity.
  This term decays only polynomially in $s$ and is insufficient
  to compensate for the union bound cost in the proof
  of~\eqref{eq:LUsingle}.
  Thus, the proof strategy used here does not yield any nontrivial bound for the failure of local uniqueness for the full interlacement process.
  It is unclear whether any nontrivial bound holds in this setting.
  \end{enumerate}
\end{remark}

\subsection{Proof of Proposition~\ref{pro:convergences}}\label{sec:proofWithIu}
Under $\tilde{\P}_N^{2r}$, for $u>0$ and $a=\sqrt{2u}$, let
\begin{equation}\label{eq:tMra}
\tilde{M}_r^{-a}\!\defeq\!
\sup\big\{|\K|\!:\!\K\text{ is a connected component of $\{x\in \tB_{2r}\!:\!\sigma_x^u\sqrt{\varphi_x^2+2\ell_x^u}\geq 0\}\cap \tB_r$}\big\},
\end{equation}
and consider the random count $N^u\defeq \sum_{x\in B_r}\1{\{\Sigma_x\cap \I^u\neq \emptyset\}}$.

As we now explain, Proposition~\ref{pro:convergences} immediately follows from the following Lemmas.
\begin{lemma}\label{lem:Ncvg}
If $1\leq r=r_N\leq \Cr{range}N$ satisfies \eqref{def:AntiCon} and $\lim_N r=\infty$, we have that for all $u>0$ such that $\sqrt{2u}\geq c(F(r))^{-\frac12}$,
\begin{equation}\label{eq:Ncvg}
  \frac{N^u}{|B_r|\vartheta(2u/{g^{2r}_N})} \stackrel{\tilde{\P}_N^{2r}}{\longrightarrow} 1.
\end{equation}
\end{lemma}

\begin{lemma}\label{lem:MNcvg}
If $1\leq r=r_N\leq \Cr{range}N$ satisfies \eqref{def:AntiCon} and $\lim_N r=\infty$, we have that for all $u>0$ and $a=\sqrt{2u}$ such that $a\geq c(F(r))^{-\frac12}$,
\begin{equation}
  \frac{\tilde{M}_r^{-a}}{N^u}
  \stackrel{\tilde{\P}_N^{2r}}{\longrightarrow} 1.
\end{equation}
\end{lemma}

\begin{proof}[Proof of Proposition~\ref{pro:convergences}]
  By~\eqref{eq:isomIu}, with $u=\frac{a^2}{2}$
$M_r^{-a}$ under $\P_N^{2r}$ has the same law as $\tilde{M}_r^{-a}$ under $\P_N^{2r}$. Proposition~\ref{pro:convergences} follows from this observation and the continuous mapping theorem.
\end{proof}

To show the proof of Lemma~\ref{lem:Ncvg} and~\ref{lem:MNcvg},
we start with the following auxiliary lemma, which
gives equivalent formulations of~\eqref{def:AntiCon}
and of the assumptions of Theorem~\ref{thm:limit_precise}.
\begin{lemma}\label{lem:equivFrh}
For every sequence $r=r_N$ such that $1\leq r_N\leq \Cr{range}N$,
\begin{align}
  \label{claim:dcrossThres}
  \lim_{N\to\infty}\frac{\capacity^{2r}(B_r)}{\capacity(B_r)}=\infty &\quad\text{if and only if}\quad  \lim_N\frac{r}{\Dcross}=\infty \text{ and } \lim_N\frac{N}{r}=\infty;\\
\label{cond:anticonRunbdd}
  \lim_{N\to\infty}\frac{g^r_N\capacity^{2r}(B_r)}{g_N\capacity(B_r)}=\infty &\quad\text{if and only if}\quad \lim_N {r}=\infty , \lim_N\frac{r}{\Dcross}=\infty \text{ and } \lim_N\frac{N}{r}=\infty.
\end{align}
\end{lemma}
The proof of Lemma~\ref{lem:equivFrh} is postponed to Section~\ref{sec:AFbox}.
We are now ready to explain the

\begin{proof}[Proof of Lemma~\ref{lem:Ncvg}]
We start with proving the following claim, which says that the on-diagonal killed Green's function is asymptotically translation invariant.
  \begin{equation}\label{eq:killed_trans_inv}
    \limsup_N \sup_{x\in B_r}  \Big|\frac{g^{2r}_N(x,x)}{g^{2r}_N(0,0)}-1\Big|=0.
  \end{equation}
  For all $x\in B_r$, we have $B(x,r)\subseteq B_{2r}\subseteq B(x,3r)$. It follows from \eqref{eq:gNK} that $g^{r}_N(0,0)\leq g^{2r}_N(x,x)\leq g^{4r}_N(0,0)$. By \eqref{eq:gDecomp} and \eqref{eq:killed_green_smallh},
  \[
  g^{4r}_N(0,0)-g^{2r}_N(0,0)=E^{4r}_0[g_N^{4r}(0,X_{H_{B_{2r}^c}}),H_{B_{2r}^c}<\infty]\leq \sup_{y\in \partial B_{2r}}g_N^{4r}(0,y) \leq \frac{c}{r\wedge h_N}.
  \]
  It follows again from~\eqref{eq:killed_green_smallh} that $g^{2r}_N(0,0)\geq c (1\vee \frac{\log(r)}{h_N})$. It then follows from this and the above display that for all $x\in B_r$,
  \[
  \Big|\frac{g^{2r}_N(x,x)}{g^{2r}_N(0,0)}-1\Big|\leq c\frac{(r\wedge h_N)^{-1}}{\big(1\vee \frac{\log(r)}{h_N}\big)},
  \]
  which tends to $0$ as $N\to\infty$, concluding the proof of~\eqref{eq:killed_trans_inv}.

Next, we claim that
\begin{equation}\label{eq:NuExp}
  \lim_N {\tilde{\E}_N^{2r}[N^u]}/\big({|B_r|\vartheta(2u/{g^{2r}_N(0,0)})}\big)=1 \quad\text{and}\quad
  \tilde{\E}_N^{2r}[N^u]\asymp |B_r|\big(\tfrac{a}{\sqrt{g^{2r}_N\!(0,0)}}\wedge 1\big)
  .
\end{equation}
Indeed, it follows from~\eqref{eq:sign_int_identity} that $\tilde{\E}_N^{2r}[N^u]=\sum_{x\in B_r} \vartheta\big(2u/{g^{2r}_N(x,x)}\big)$.
By~\eqref{eq:killed_trans_inv} and~\eqref{eq:vartheta_scaling},
$\vartheta(2u/g_N^{2r}(x,x))/\vartheta(2u/g_N^{2r}(0,0))\to1$
uniformly in $x\in B_r$ and $a>0$, which proves the first claim
in~\eqref{eq:NuExp}. The second follows from~\eqref{eq:vartheta_bound}.

We now show that $\variance^{\tilde{\P}_N^{2r}}\!(N^u)/\tilde{\E}_N^{2r}[N^u]^2\to0$ as $N\to\infty$, which implies that~\eqref{eq:Ncvg} also holds in the $L^2$ sense. By \eqref{eq:cov_sigma}, \eqref{eq:killed_green_smallh} and a calculation similar to \cite[(3.26) and (3.27)]{rodriguez_2025_on},
\[
 \variance^{\tilde{\P}_N^{2r}}\!(N^u)
 \!=\!\!\!\!\sum_{x,y\in B_r}\!\!\!\text{Cov}^{\tilde{\P}_N^{2r}}\big(\1_{\{\Sigma_x\cap \I^u\neq \emptyset\}},\1_{\{\Sigma_y\cap \I^u\neq \emptyset\}}\big)
 \!\leq\! 
 \frac{c}{g^{2r}_N(0,0)}\!\!\sum_{x,y\in B_r}\!\! g^{2r}_N(x,y)\leq \frac{c|B_r|^2}{g^{2r}_N(0,0)F^{2r}(r)}.
\]
Combining this with~\eqref{eq:NuExp} and the assumption $\sqrt{2u}\geq cF(r)^{-\frac12}$ then yields that $\variance^{\tilde{\P}_N^{2r}}\!(N^u)/\tilde{\E}_N^{2r}[N^u]^2\leq (F(r)/F^{2r}(r))\vee (g^{2r}_N F^{2r}(r))^{-1}$, which tends to $0$ as $N\to\infty$.

\end{proof}

\begin{proof}[Proof of Lemma~\ref{lem:MNcvg}]
  We start with the upper bound, which is simpler. 
Let $\hat{\K}_r$ denote the largest connected component in $\{x\in \tB_{2r}:\Sigma_x\cap \I^u=\emptyset,\sigma_x^u=1\}\cap \tB_r$.
In view of~\eqref{eq:tMra}, $\{\tilde{M}_r^{-a}>N^u\}\subseteq\{ \tilde{M}_r^{-a}=|\hat{\K}_r|\}$ since every vertex whose sign cluster is visited by $\I^u$ would be included in the sum $N^u$. It then follows that for all $\epsilon\in(0,\frac12)$, there exists $C(\epsilon)<\infty$ such that for all $N\geq C(\epsilon)$,
\begin{align}\label{eq:MNub}
\begin{split}
  \tilde{\P}_N^{2r}(\tilde{M}_r^{-a}>N^u)
  &\leq \tilde{\P}_N^{2r}(|\hat{\K}_r|> N^u)\\
  &\leq 
  \tilde{\P}_N^{2r}(|\hat{\K}_r|> (1-\epsilon)|B_r|\vartheta(a^2/{g^{2r}_N})) +\tilde{\P}_N^{2r}(N^u< (1-\epsilon)|B_r|\vartheta(a^2/{g^{2r}_N}))\\
  &\leq c\,{\tilde{\E}_N^{2r}[|\hat{\K}_r|]}\big({|B_r|\vartheta(a^2/{g^{2r}_N})}\big)^{-1}
  + \tfrac{\epsilon}{2}
\end{split}
\end{align}
where we used Markov's inequality and the fact that 
${N^u}/({|B_r|\vartheta(a^2/{g^{2r}_N})}) \stackrel{\tilde{\P}_N^{2r}}{\longrightarrow} 1$ in the last bound.
Now, since $|\hat{\K}_r|$ is stochastically dominated by $M_r^0$ under $\P_N^{2r}$,
it follows from~\eqref{eq:vartheta_bound} and~\eqref{eq:upperMR0} that the term involving the expectation in~\eqref{eq:MNub} is at most
$c\big((F(r)/F^{2r}(r))^{1/2}\vee (g^{2r}_N F^{2r}(r))^{-1/2}\big)$, which is bounded by $\epsilon/2$ for all $N\geq C(\epsilon)$ with $C(\epsilon)$ sufficiently large. This concludes the proof of the upper bound.

For the lower bound, let $\lambda \in(0,\frac14)$ soon to be fixed. Consider the following random counts
\begin{align}
\begin{split}
  &N_1^u\defeq \sum_{x\in B_r} \1{\{\Sigma_x\subseteq\tB_{(1-2\lambda)r},
  \Sigma_x\cap \I^u_+\neq\emptyset\}},
  %\{x\in B_r:\Sigma_x\cap \I^u\cap \tB_{(1-2\lambda)r}\neq\emptyset\}
  \quad\quad  N_2^u\defeq\!\!\!\!\!\!\! \sum_{x\in B_r\setminus B_{(1-4\lambda)r}}
  \!\!\!\!\!
  \1{\{\Sigma_x\cap \I^u \neq\emptyset\}},\\
  &N_3^u\defeq\!\!\!\!\! \sum_{x\in B_{(1-4\lambda)r}}
  \!\!\!\!\!
  \1{\{\Sigma_x\cap \tB_{(1-2\lambda)r}^c\neq\emptyset\}},
  \quad\quad
  N_4^u\defeq \sum_{x\in B_r} \1{\{\Sigma_x \cap \I^u_+=\emptyset, \Sigma_x\cap\I^u\neq\emptyset\}}
\end{split}
\end{align}
It follows from the definition of $N^u$ immediately below \eqref{eq:tMra} that $N^u\leq \sum_{i=1}^4N_i^u$ and hence for all $\epsilon\in(0,1)$,
\begin{equation}
  \tilde{\P}_N^{2r}(\tilde{M}_r^{-a}\leq (1-\epsilon)N^u)
  \leq \tilde{\P}_N^{2r}(\tilde{M}_r^{-a}< N_1^u)
  +  \tilde{\P}_N^{2r}(N_2^u+N_3^u+N_4^u\geq \tfrac\epsilon2 N^u).
\end{equation}
Note that 
$\{\LU_{u,(1-2\lambda)r,r}^+(0)\}\subseteq\{N_1^u\leq \tilde{M}_r^{-a}\}$. Indeed, on the event $\{\LU_{u,(1-2\lambda)r,r}^+(0)\}$, all vertices counted by $N_1^u$ are part of one connected component in $\{x\in \tB_{2r}:\sigma_x^u\sqrt{\varphi_x^2+2\ell_x^u}\geq 0\}\cap \tB_r$. And hence for all $\lambda\in(0,\frac14)$,
\[
\tilde{\P}_N^{2r}(\tilde{M}_r^{-a}< N_1^u) \leq \tilde{\P}_N^{2r}(\LU_{u,(1-2\lambda)r,r}^+(0)^c)\leq c(\lambda)\,\exp\big\{ -c'(\lambda)(u F^{2r}(r))^{1/3}\big\},
\]
which, since $\lim_N uF^{2r}(r)=\infty$, tends to $0$ as $N\to\infty$. It remains to show the second term vanishes in the limit. First note that
\[
\tilde{\P}_N^{2r}(\sum_{i=2}^4N_i^u\geq \tfrac\epsilon2 N^u)
\leq \tilde{\P}_N^{2r}(\sum_{i=2}^4N_i^u\geq \tfrac\epsilon4 |B_r|\vartheta(a^2/{g^{2r}_N}))
+ \tilde{\P}_N^{2r}(N^u\leq \tfrac{1}{2} |B_r|\vartheta(a^2/{g^{2r}_N})).
\]
It follows from~\eqref{eq:sign_int_identity}, \eqref{eq:onearm_killed} and $|B_r\setminus B_{(1-4\lambda)r}|\leq c\lambda |B_r|$ that
\[
\tilde{\E}_N^{2r}[N_2^u]\leq  c\lambda |B_r|\vartheta(a^2/{g^{2r}_N})\quad\text{and}\quad
\tilde{\E}_N^{2r}[N_3^u]\leq c|B_r|({g^{2r}_N F^{2r}(\lambda r)})^{-1/2}.
\]
Furthermore, since $\sup_{x\in B_r}(1-P_x(H_{\partial B_{2r}}<\infty))\leq c(r/N)^2$, we have by~\eqref{eq:h-isom},~\eqref{eq:vartheta_scaling} and~\eqref{eq:sign_int_identity} that $\tilde{\E}_N^{2r}[N_4^u]\leq c(\frac rN)^2 |B_r|\vartheta(a^2/{g^{2r}_N})$. Combining this with the last display and Lemma~\ref{lem:equivFrh} yields, for every $\lambda\in(0,\frac14)$,
\begin{equation}\label{eq:count_less_exp}
  \limsup_N
\frac{\tilde{\E}_N^{2r}[\sum_{i=2}^4N_i^u]}
{|B_r|\vartheta(a^2/g_N^{2r})}\leq c\lambda.
\end{equation}
It then follows from Markov's inequality and~\eqref{eq:Ncvg} that
$\limsup_N\tilde{\P}_N^{2r}(N_2^u+N_3^u+N_4^u\geq \tfrac\epsilon2 N^u)\leq c\epsilon^{-1}\lambda$.
Letting $\lambda\downarrow0$ then completes the proof.
\end{proof}

\section{Proof of Theorems~\ref{thm:antiCon} and~\ref{thm:limit}} \label{sec:final}
In this section, we establish the precise limiting behaviour of $M_r$
in Theorem~\ref{thm:limit_precise}, using Proposition~\ref{pro:convergences} as a key input.
As we explain shortly, Theorem~\ref{thm:limit_precise} directly implies
Theorem~\ref{thm:limit}. 
%The key input isProposition~\ref{pro:convergences}, whose proof relies on theconnectivity properties of the interlacement process and isdeferred to Section~\ref{sec:Iu}.
We finish the section by proving Theorem~\ref{thm:antiCon}
using Theorem~\ref{thm:limit_precise} and
Proposition~\ref{lem:ConLb}.

\begin{theorem}\label{thm:limit_precise}
If $1\leq r=r_N\leq \Cr{range}N$ satisfies \eqref{def:AntiCon} and $\lim_N r=\infty$, then $\E_N[M_r]\sim |B_r|\frac{1}{\pi}\arctan(\sigma_{N,r})$, i.e.~\eqref{eq:EcvgClaim} holds, and
  \begin{equation}\label{eq:limit_precise}
    \overline{M}_r - \pi \frac{2\Phi\big(
    (g^{r}_N)^{-\frac12}  
    (\eta_0^{2r})_+\big)-1}{\arctan(\sigma_{N,r})} \stackrel{\P_N}{\longrightarrow} 0.
  \end{equation}
\end{theorem}
Theorem~\ref{thm:limit} is a direct consequence of
Theorem~\ref{thm:limit_precise}. Indeed, by \eqref{eq:varEta0} and \eqref{eq:ggradCable}, we have that if $r$ satisfies \eqref{def:AntiCon} and $\lim_Nr=\infty$,
\begin{equation}\label{eq:vareta_asymp}
  \variance_N(\eta_0^{2r})\sim g_N(0,re_1),\quad \text{as $N\to\infty$}.
\end{equation}
Since the limit of $\sigma_{N,r}$ exists, $\pi \frac{2\Phi\big(
    (g^{r}_N)^{-\frac12}  
    (\eta_0^{r})_+\big)-1}{\arctan(\sigma_{N,r})}$ 
converges
in distribution to $F_\sigma(Z_+)$, where
\begin{equation}\label{eq:Fsigma}
    F_{\sigma}(x)\defeq
\begin{cases}
\sqrt{2\pi}\,x, &\sigma=0,\\[1ex]
\displaystyle
\pi\frac{2\Phi(\sigma x)-1}{\arctan \sigma}, & 0<\sigma<\infty,\\[1.5ex]
2\mathbf1_{\{x>0\}}, & \sigma=\infty.
\end{cases}
\end{equation}
By \eqref{eq:limit_precise} and the continuous mapping theorem,
\begin{equation}\label{eq:limit}
    \overline{M}_r \stackrel{d}{\longrightarrow}  F_{\sigma}(Z_+).
\end{equation}
Since $M_r/|B_r|=\overline M_r\,\E_N[M_r]/|B_r|$, Theorem~\ref{thm:limit} follows from \eqref{eq:limit} and \eqref{eq:EcvgClaim}, using that $\arctan(\sigma_{N,r})\sim\sigma_{N,r}$ when $\sigma=0$.

\smallskip
Before giving the proof, let us explain the mechanism behind Theorem~\ref{thm:limit} and the value $\frac12$ in Theorem~\ref{thm:antiCon}. 
By the Markov property \eqref{eq:field_decomp}, $\varphi\!=\!\psi^{2r}\!+\eta^{2r}$, where the harmonic part $\eta^{2r}$ is generated by the field on $\partial B_{2r}$ and is independent of $\psi^{2r}$. Two scales compete. On the one hand, the oscillation of $\eta^{2r}$ over $\tB_r$ is of order $(r\!\wedge\! h_N)^{-1/2}$, and Proposition~\ref{prop:shift_h0} shows that a shift of the level by an amount $b\gg (r\!\wedge\! h_N)^{-1/2}$ suffices to compare the local level sets of $\varphi$ with those of the killed field $\psi^{2r}$ at a shifted level. On the other hand, by \eqref{eq:varEta0}, the value $\eta_0^{2r}$ is a centred Gaussian whose standard deviation is at least of order $F(r)^{-1/2}$. Under \eqref{def:AntiCon}, one has $F(r)\!\ll\! r\!\wedge\! h_N$ by \eqref{claim:dcrossThres}, so that $|\eta_0^{2r}|\!\gg\!(r\!\wedge\! h_N)^{-1/2}$ with probability tending to one, and $\eta^{2r}$ acts on $\tB_r$ as an essentially constant shift of the level. With probability tending to $\frac12$ this shift is negative, in which case $M_r(\varphi)\leq M_r(\psi^{2r})$ and the expectation of the latter is of strictly smaller order than $\E_N[M_r]$; this yields \eqref{eq:anticonSmall} and explains why the limiting
law in Theorem~\ref{thm:limit} has an atom of mass $\frac12$ at zero.
Heuristically, with the remaining probability $\frac12$ the shift is positive and pushes the field in $\tB_r$ into a supercritical state, which, by means of Proposition~\ref{pro:convergences}, admits a precise limiting law.
Under \eqref{def:Con}, by contrast, $F(r)\!\asymp\! r\!\wedge\! h_N$, so that the harmonic part is of the same order as its own oscillation over $\tB_r$, and it does not single out a sign.

\smallskip

We now explain how Proposition~\ref{pro:convergences} yields Theorem~\ref{thm:limit_precise}.
\begin{proof}[Proof of Theorem~\ref{thm:limit_precise}]
We abbreviate $\eta_0^{2r}=\eta_0$ and start by showing convergence on the event $\{\eta_0\leq 0\}$.
Using Lemma~\ref{pro:upperMRa} and the Cauchy-Schwarz inequality gives,
\begin{align*}
 \frac{\E_N[M_r\1\{\eta_0\leq 0\}]}{\E_N[M_r]}
 &\leq \frac{1}{{\E_N[M_r]}}\Big(\E_N[M_r^0(\psi^{2r})] +  
 \E_N\big[M_r\1\big\{\eta_0\leq 0\}\1\{M_r^0(\psi^{2r})< M_r\big\}\big]\Big)\\
&\leq c\Big(\frac{g_N F(2r)}{g^{2r}_N F^{2r}(r)}\Big)^{1/2}
+c \P_N\big(\eta_0\leq 0,M_r^0(\psi^{2r})< M_r\big)^{1/2}.
\end{align*}
For $\delta\in(0,1)$, using~\eqref{eq:shift_anticon} and~\eqref{eq:varEta0} and splitting the above probability according to whether
$\eta_0\leq-\delta/\sqrt{F(r)}$ gives
$\limsup_N\P_N(\eta_0\leq0,M_r^0(\psi^{2r})<M_r)
\leq c\delta$. Combining this with~\eqref{cond:anticonRunbdd} and letting $\delta\downarrow0$ then yields for all $\epsilon\in(0,1)$,
\begin{equation}\label{eq:etaneg}
  \lim_N \frac{\E_N[M_r\1\{\eta_0\leq 0\}]}{\E_N[M_r]}=0 \text{ and hence, }
  \lim_N\P_N\Big(|\overline{M}_r|>\epsilon,\eta_0\leq 0\Big)=0 ,
\end{equation}
where we used Markov's inequality to obtain the second claim in \eqref{eq:etaneg}.

It remains to consider what happens on the event $\{\eta_0>0\}$. Let $\Delta_{N,r}\defeq M_r- |B_r|\vartheta(\eta_0^2/g_N^{2r})$. We claim that for all $\epsilon\in(0,1)$, 
\begin{align}\label{eq:etapos}
 \lim_N\frac{
\E_N\big[\big|\Delta_{N,r} \big|\1\{\eta_0>0\}\big]
}{\E_N[M_r]}=0 \quad\text{and}\quad
  \lim_N\P_N\Big(\frac{\big|\Delta_{N,r} \big|}{\E_N[M_r]}>\epsilon,\eta_0>0\Big)=0 .
\end{align}
Here, as for~\eqref{eq:etaneg}, it suffices to show the first claim in \eqref{eq:etapos} since the second follows from the first and Markov's inequality.
It follows from~\eqref{eq:pro_cvg} and~\eqref{eq:vartheta_scaling} that for all $\epsilon\in(0,1), b>c/\sqrt{F(r)}$ and $\gamma\in(0,c(\epsilon))$ with $c(\epsilon)$ small enough, we have that
\begin{equation}\label{eq:cvg_conditional_prob}
  \lim_N\P_N\Big(\Big|\frac{M_r^{-(1\pm\gamma)b}(\psi^{2 r}) }{|B_r|\cdot \vartheta(b^2/g^{2r}_N)}-1\Big|>\epsilon
\Big)=0
\quad \text{and hence}\quad \frac{M_r}{|B_r|\cdot \vartheta(\eta_0^2/g^{2r}_N)} \stackrel{\P_N(\cdot|\eta_0>0)}{\Rightarrow}1,
\end{equation}
where we used \eqref{eq:shift_anticon} to obtain the convergence in probability statement.
Now, by the Cauchy-Schwarz inequality, writing $A_N(\epsilon)\defeq
\big\{\eta_0>0,\
\big|\frac{M_r}{|B_r|\vartheta(\eta_0^2/g_N^{2r})}-1\big|
>\epsilon\big\}$,
\begin{align}\label{eq:exp_cvg_bound}
\begin{split}
&\frac{
\E_N\big[\big|\Delta_{N,r} \big|\1\{\eta_0>0\}\big]
}{\E_N[M_r]} 
\leq
\epsilon\,
\frac{|B_r|\cdot\E_N[\vartheta(\eta_0^2/g_N^{2r})
\1\{\eta_0>0\}]}{\E_N[M_r]}+
\frac{\E_N\big[
\Delta_{N,r}^2
\big]^{1/2}}{\E_N[M_r]}
\P_N(A_N(\epsilon))^{1/2}.
\end{split}
\end{align}
By~\eqref{eq:gDecomp} and~\eqref{eq:g2-asymp-unif},
$g_N\leq g_N^{2r}+cF(r)^{-1}$. Combining this
with~\eqref{eq:thmM0rExp} and~\eqref{eq:vartheta_bound} yields
\begin{align}
\begin{split}
  \frac{|B_r|\vartheta(\eta_0^2/g_N^{2r})}{\E_N[M_r]} \leq  c \sqrt{g_NF(r)}\Big(\big({|\eta_0|}/{\sqrt{g_N^{2r}}}\big)\wedge 1\Big)
  &\leq c\sqrt{1+g^{2r}_NF(r)}\Big(\big({|\eta_0|}/{\sqrt{g_N^{2r}}}\big)\wedge 1\Big)
  \\
  &\leq c\big(1+|\eta_0|\sqrt{F(r)}\big).
\end{split}
\end{align}
The preceding display and~\eqref{eq:varEta0} imply that
$\E_N[(|B_r|\vartheta(\eta_0^2/g_N^{2r})/\E_N[M_r])^2]\leq c$.
Together with Lemma~\ref{pro:upperMRa}
and the inequality $(a-b)^2\leq2a^2+2b^2$, this shows that both
fractions on the r.h.s of~\eqref{eq:exp_cvg_bound}
are uniformly bounded. Since $\P_N(A_N(\epsilon))\to0$
by~\eqref{eq:cvg_conditional_prob}, letting $N\to\infty$
and then $\epsilon\downarrow0$ proves~\eqref{eq:etapos}. 
In view of~\eqref{eq:etaneg}, we have in fact shown that
\begin{equation}\label{eq:EcvgClaim_int}
\lim_N
\frac{|B_r|\E_N[\vartheta(\eta_0^2/g_N^{2r})
\1\{\eta_0>0\}]}
{\E_N[M_r]}=1.
\end{equation}
To conclude, let $X,Y$ be independent standard Gaussians.
By rotational invariance of $(X,Y)$, 
\begin{align*}
\E_N[\vartheta(\eta_0^2/g_N^{2r})\1\{\eta_0>0\}]
&=\P\big(X>0,\ |Y|\leq
\sqrt{\tfrac{\variance_N(\eta_0)}{g_N^{2r}}}\,X\big)
=\frac1\pi\arctan
\Big(\sqrt{\tfrac{\variance_N(\eta_0)}{g_N^{2r}}}\Big),
\end{align*}
since the event describes a sector of angle
$2\arctan(\sqrt{\variance_N(\eta_0)/g_N^{2r}})$ in the plane.
The preceding identity, together with~\eqref{eq:vareta_asymp}
and~\eqref{eq:EcvgClaim_int}, proves~\eqref{eq:EcvgClaim}.
The desired convergence~\eqref{eq:limit_precise} then follows
from~\eqref{eq:EcvgClaim}, \eqref{eq:etapos}
and~\eqref{eq:etaneg}.

\end{proof}

We finish with the
\begin{proof}[Proof of Theorem~\ref{thm:antiCon}]
We start with~\eqref{eq:tailEstlargestVolgena=0}. By Markov's inequality and \eqref{eq:upperMR0} with $K=\emptyset$, one has, for any $t\geq1$ and $1\leq r\leq \Cr{range}N$, $\P_N(M_r\geq t\bv(r))\leq ct^{-1}$. The complementary bound follows from Proposition~\ref{lem:ConLb} and again \eqref{eq:upperMR0} with $K=\emptyset$, completing the proof of~\eqref{eq:tailEstlargestVolgena=0}.

For~\eqref{eq:anticonBig} and ~\eqref{eq:anticonSmall}, note that it suffices to show that for all $\epsilon\in(0,1)$, there exists $\Cl[c]{anticont}<\infty$ such that for all $t\geq \Cr{anticont}$,
\begin{align}
  \label{eq:anticonLowerBound1}
  &\liminf_N\P_N(\overline{M_r}< t^{-1})\geq \tfrac{1}{2}-\epsilon
  \quad\text{and}\quad\\
  \label{eq:anticonLowerBound2}
  &\liminf_N\P_N(t^{-1}\leq \overline{M_r}\leq t)\geq \tfrac{1}{2}-\epsilon 
  .
\end{align}
Indeed, the required upper bound in \eqref{eq:anticonSmall} follows from \eqref{eq:anticonLowerBound2} since \eqref{eq:anticonLowerBound2} implies that for all $t\geq \Cr{anticont}$,
$
\limsup_N \P_N(\overline{M_r}<t^{-1})\leq 1-\liminf_N P_N(t^{-1}\leq \overline{M_r}\leq t)\leq \frac12 +\epsilon.
$
Similarly for \eqref{eq:anticonBig}.

We now proceed to show \eqref{eq:anticonLowerBound1}. By passing to subsequences, it suffices to consider
the cases where $r$ is uniformly bounded or
$\lim_N r=\infty$.
Suppose first that $\sup_N r\leq C$. We will show that
\begin{equation}\label{eq:half_empty}
  \liminf_N\P_N(M_r=0)\geq \tfrac12.
\end{equation}
First note that \eqref{def:AntiCon} and \eqref{claim:dcrossThres} imply that we must have $\lim_N\Dcross=0$ in this case, which, in view of \eqref{def:Dcross} and \eqref{eq:gN0}, is equivalent to $\lim_N g_N=\infty$. Since $|B_r|\leq c$, it follows from \eqref{eq:ggradDisc} that
$
\sup_{x\in B_r}\E_N[(\varphi_x-\varphi_0)^2]\leq 2\sup_{x\in B_r} |g_N(0,0)-g_N(0,x)|\leq c.
$
A union bound and Markov's inequality then yield, for every $\delta>0$,
\[
\P_N\Big(\max_{x\in B_r}|\varphi_x-\varphi_0|
>\delta\sqrt{g_N}\Big)
\leq \frac{c}{\delta^2g_N}.
\]
It then follows from the above display that, 
\begin{align*}
  \P_N(M_r=0)
  &\!\geq\! \P_N(\varphi_x<0 \text{ for all $x\in B_r$})\\
  &\!\geq\!\P_N(\varphi_0<-2\delta\sqrt{g_N} )-\P_N\Big(\max_{x\in B_r}|\varphi_x-\varphi_0|>\delta\sqrt{g_N}\Big)\!\geq\! \P_N(\varphi_0<-2\delta\sqrt{g_N} )-\tfrac{c}{\delta^2g_N}.
\end{align*}
Since $\lim_Ng_N=\infty$,
first letting $N$ tend to infinity and then
letting $\delta$ decrease to zero proves \eqref{eq:half_empty}. This also concludes the proof of \eqref{eq:anticonLowerBound1} when $\sup_Nr<\infty$ since $\bv(r)\geq c$ by Corollary~\ref{the:mainMra}.
For \eqref{eq:anticonLowerBound2},
by \eqref{eq:half_empty} and
$\P_N(M_r=0)\leq\P_N(\varphi_0<0)=1/2$,
we have $\lim_N\P_N(M_r=0)=1/2$.
Since $1/2\leq\E_N[M_r]\leq|B_r|\leq C$,
for all sufficiently large $t$ the event
$\{t^{-1}\leq\overline M_r\leq t\}$ coincides with
$\{M_r>0\}$. This completes the proof in the regime $\sup_N r\leq C$.

When $\lim_N r=\infty$, we are in the regime of Theorem~\ref{thm:limit_precise}.
Since $2\Phi(0)-1=0$ and $\eta^{2r}_0$ is a centered Gaussian, ${2\Phi\big(
(g^{2r}_N)^{-\frac12}  
(\eta_0^{2r})_+\big)-1}=0$ on the event $\{\eta^{2r}_0\leq 0\}$, which happens with probability $\frac12$.
It then follows from~\eqref{eq:limit_precise} that~\eqref{eq:anticonLowerBound1} holds.
For~\eqref{eq:anticonLowerBound2}, passing to a further subsequence along which
$\lim_N\sigma_{N,r}=\sigma\in[0,\infty]$ and using Theorem~\ref{thm:limit} gives,
\[
\liminf_N\P_N(t^{-1}\leq\overline M_r\leq t)
\geq\P(Z\geq 0)-\P(F_\sigma(Z)<t^{-1},Z\geq 0)-\P(F_\sigma(Z)>t,Z\geq 0).
\]
By Markov's inequality, we have that $\P(F_\sigma(Z)>t,Z\geq 0)\leq ct^{-1}$.
Furthermore, by~\eqref{eq:vartheta_bound} and~\eqref{eq:vartheta_scaling},
$F_\sigma(z)\geq c(z\wedge1)$ for $z>0$, uniformly in $\sigma$, which yields $\P(F_\sigma(Z)<t^{-1},Z\geq 0)\leq ct^{-1}$.
This completes the proof of~\eqref{eq:anticonLowerBound2} upon taking $t$ large enough.
\end{proof}

{
\begin{remark}
A somewhat related sensitivity to boundary conditions appears for the two-dimensional metric-graph Gaussian free field in a box with Dirichlet boundary conditions: by \cite{lupu_2018_the,bi_2026_onear}, the one-arm probability to an intermediate scale decays like $(\log N)^{-1/2}$ at every level, whereas the one-arm probability to the boundary of the box exhibits a transition across level $0$.
\end{remark}
}

\noindent\textbf{Acknowledgements.} 
The author wishes to thank her advisor Pierre-Fran\c cois Rodriguez for many helpful discussions and comments on the draft of this manuscript.
This research is supported by the CDT in Mathematics of Random Systems and by the  President's PhD Scholarship scheme at Imperial College, London.

\appendix
\section{Estimates for $\Fbox$, $F^r$}\label{sec:AFbox}

\begin{proof}[Proof of Lemma~\ref{lem:loop_crossing_asymp}]
We start with \eqref{eq:ConCrossing}. By \eqref{eq:capB} and \eqref{eq:cap-killed}, the assumption
implies that $F(r)\geq cK^{-1}(r\wedge h_N)$.
Let $r'=\delta r$, where $\delta\in(0,1)$ will be chosen
below. If $r\leq h_N$, then $F(r')\leq\delta r$, and hence
$F(r')/F(r)\leq cK\delta$. If instead $r\geq h_N$, then,
since $r\leq\Cr{range}N$,
\[
\frac{F(r')}{F(r)}
\leq \frac{cK}{\log^\star(N/(\delta r))}
\leq \frac{cK}{\log^\star(1/(\Cr{range}\delta))}.
\]
Choosing $\delta=\Cr{range}^{-1}\exp\{-\Cr{rrplbK} K/\lambda \}$ with $\Cr{rrplbK}$ sufficiently large
ensures that $F(r')\leq\lambda F(r)$ in both cases.
Choosing $r\geq\Cr{rlbK}$ with $\Cr{rlbK}(K,\lambda)$ large enough thus ensures that $r'\geq1$. This completes the proof of \eqref{eq:ConCrossing}.

To prove \eqref{eq:lambdaKilledCap}, note that by taking $\Cr{lambdaKilledCap}$ small enough, we may assume that $1\leq r'\leq r/4$. It then follows from
\eqref{eq:gDecomp} and \eqref{eq:g2-asymp-unif} that for all $x,y\in B_{r'}$, $0\leq g_N(x,y)-g_N^r(x,y)\leq \frac{c}{F(r)}$. Combining this with \eqref{eq:CapVarForm}
therefore yields
\[
\frac1{\capacity^r(B_{r'})}
\geq \frac1{\capacity(B_{r'})}-\frac{c}{F(r)}\geq \frac{1}{2\capacity(B_{r'})},
\]
where the last inequality follows upon choosing $\lambda \in(0,\Cr{lambdaKilledCap})$ with $\Cr{lambdaKilledCap}$ small enough, since by \eqref{eq:capB}, $\capacity(B_{r'})\leq cF(r')\leq c\lambda F(r)$.
Consequently, by \eqref{eq:capB},
\[
\frac{\capacity^r(B_{r'})}{\capacity(B_r)}
\leq 2\frac{\capacity(B_{r'})}{\capacity(B_r)}
\leq c\frac{F(r')}{F(r)}
\leq c\lambda.
\]
\end{proof}

\begin{proof}[Proof of Lemma~\ref{lem:equivFrh}]
We start with \eqref{claim:dcrossThres}.
By \eqref{eq:capB} and \eqref{eq:cap-killed}, we have $\frac{\capacity^{2r}(B_r)}{\capacity(B_r)}
\asymp \frac{r\wedge h_N}{F(r)}.$
We first show that $\lim_N r/\Dcross=\infty$ and
$\lim_N N/r=\infty$ imply that the capacity ratio diverges.
If $r\leq h_N$, then $F(r)\asymp r\wedge\Dcross$, and hence $\lim_N\frac{r}{r\wedge\Dcross}=\infty.$
If instead $r\geq h_N$, then
$F(r)\leq h_N/\log^\star(N/r)$. In this case,
\[
\frac{\capacity^{2r}(B_r)}{\capacity(B_r)}
\geq c\log^\star(N/r),
\quad\text{and}\quad
\lim_N\log^\star(N/r)=\infty.
\]
This proves the first implication.

For the converse, suppose that the capacity ratio diverges.
If it's not the case that $\lim_N r/\Dcross=\infty$, then, along a
subsequence, $r\leq C\Dcross$ for some $C<\infty$.
On this subsequence, $F(r)\asymp r$, and therefore
$\frac{\capacity^{2r}(B_r)}{\capacity(B_r)}
\asymp \frac{r\wedge h_N}{r}\leq c,$
a contradiction. Similarly, if $\lim_N N/r=\infty$ fails,
then there exists $c>0$ such that $cN<r\leq \Cr{range}N$ along a subsequence. On this subsequence,
$F(r)\asymp r\wedge h_N$, so the capacity ratio is again
bounded. This completes the proof.

We now prove \eqref{cond:anticonRunbdd}. By \eqref{eq:gDecomp}
and \eqref{eq:g2-asymp-unif},
$0\leq g_N-g_N^r\leq c/F(r)$. Since also $g_NF(r)\geq c$,
we obtain
\[
\frac{g_N\capacity(B_r)}{g_N^r\capacity^{2r}(B_r)}=\frac{g_N^r\capacity(B_r)+(g_N-g_N^r)\capacity(B_r)}{g_N^r\capacity^{2r}(B_r)}
\asymp
\frac{\capacity(B_r)}{\capacity^{2r}(B_r)}
+\frac{1}{g_N^r\capacity^{2r}(B_r)}.
\]
By \eqref{claim:dcrossThres}, the first term vanishes
if and only if $\lim_N r/\Dcross=\infty$ and
$\lim_N N/r=\infty$.
For the second term, \eqref{eq:killed_green_smallh} and \eqref{eq:cap-killed}
give $g_N^r\capacity^{2r}(B_r)\asymp r$ when $r\leq h_N$, and
$g_N^r\capacity^{2r}(B_r)\asymp h_N+\log^\star(r/h_N)$ when
$r\geq h_N$. Thus this term vanishes if and only if
$\lim_N r=\infty$. Combining these observations proves
\eqref{cond:anticonRunbdd}.
\end{proof}

\section{Gradient of heat kernel}\label{A:grad}
\begin{proof}[Proof of Lemma \ref{lem:ct-hk-gradient}]
We start with an estimate on the discrete time heat kernel.
Let $q_n(\cdot,\cdot)$ be the transition probability for the random walk on $\Z^d$ such that
\[
  q_1(0,0)=\frac12,
  \quad\text{and}\quad
  q_1(0,\pm e_i)=\frac{1}{4d},\quad \text{ for all } 1\leq i\leq d.
\]
We abbreviate $q_n(y)=q_n(0,y)$ for all $n\geq 0$ and claim that, for all $y\in\mathbb Z^d\setminus\{0\}$, $n\geq 1$ and every coordinate vector $e$,
\begin{equation}
\label{eq:dt-hk-gradient}
  |\nabla_e q_n(y)|
  \leq
  c(d) n^{-(d+1)/2}
  \exp\left\{-c'\!(d)\frac{|y|^2}{n}\right\}.
\end{equation}
For $M\geq 2$ soon to be fixed, we first consider the case when $M\sqrt{n\log(n\vee2)}\leq |y|$. It follows from \cite[2.50]{lawler_2010_rando} that
\[|\nabla_e q_n(y)|\leq q_n(y+e)+q_n(y) \leq ce^{-c'|y|^2/n}.\]
Now we may choose $M=M(d)$ large enough that $e^{-c'\frac{|y|^2}{2n}}\leq cn^{-(d+1)/2}$,
which completes the proof of
\eqref{eq:dt-hk-gradient} in this regime.

We now move on to the case $1\leq |y|<M\sqrt{n\log(n\vee2)}$. Let $\bar q_n(y)\defeq (d/(\pi n))^{d/2}\exp\{-d|y|^2/n\}$. By \cite[Theorem 2.3.8]{lawler_2010_rando}, we have that for every $k\geq 3$, we can find functions $u_j:\Z^d\to\R$ in the form of \cite[(2.25)]{lawler_2010_rando} for all $3\leq j\leq k$ such that 
\[  q_n(y)
  =
  \bar q_n(y)
  +
  \sum_{j=3}^{k}\frac{u_j(y/\sqrt{n})}{n^{(d+j-2)/2}}
  +
  \tilde E_{k,n}(y),\]
where $\sup_{y\in\mathbb Z^d}|\tilde E_{k,n}(y)|
  \leq c(k) n^{-(d+k-1)/2}.$
It then follows from the above display that for all $y\in \Z^d\setminus\{0\}$ and $n\geq 1$
\begin{equation}\label{eq:dt-gradient-expansion}
\nabla_e q_n(y)
  =
  \nabla_e \bar q_n(y)
  +
  \sum_{j=3}^{k}\frac{u_j(\frac{y+e}{\sqrt{n}})-u_j(\frac{y}{\sqrt{n}})}{n^{(d+j-2)/2}}
  + E_{k,n}(y),
\end{equation}
where $\sup_{y\in\mathbb Z^d}|E_{k,n}(y)|
\leq c(k) n^{-(d+k-1)/2}.$ We now control the terms on the r.h.s of \eqref{eq:dt-gradient-expansion}. 
The mean-value theorem and \cite[(2.29)]{lawler_2010_rando} give
\begin{align}\label{eq:grad-mvt}
\begin{split}
	&|\nabla_e\bar q_n(y)|
	\leq c|y|\big(\frac{d}{n}\big)^{d/2+1}e^{-\frac{d|y|^2}{n}}
	\leq c(d)n^{-\frac{d+1}{2}} \frac{|y|}{\sqrt{n}}
	e^{-\frac{d|y|^2}{n}}\leq c(d)n^{-\frac{d+1}{2}} 
	e^{-\frac{d|y|^2}{2n}}\text{ and,}\\
	&\Big|
	u_j\Big(\frac{y+e}{\sqrt{n}}\Big)-u_j\Big(\frac{y}{\sqrt{n}}\Big)
	\Big|
	\leq
	c(j)n^{-\frac12} \Big(\big(\frac{y}{\sqrt{n}}\big)^{j+1}+1\Big)e^{-\frac{c'\!(d)|y|^2}{n}}, \text{ for all $3\leq j\leq k$.}
\end{split}
\end{align}
Thus, the sum in \eqref{eq:dt-gradient-expansion} is bounded by
\begin{equation}\label{eq:u-grad}
	c(k)n^{-\frac{d+1}{2}}n^{-\frac12}\Big(\big(\frac{y}{\sqrt{n}}\big)^{k+1}+1\Big)e^{-\frac{c'\!(d)|y|^2}{n}}\leq c''\!(d,k)n^{-\frac{d+2}{2}}e^{-\frac{c'\!(d)|y|^2}{2n}}.
\end{equation}
It remains only to control $E_{k,n}(\cdot)$ in
\eqref{eq:dt-gradient-expansion}, which is where the assumption $|y|<M\sqrt{n\log(n\vee2)}$ is used. By choosing $k=k(M)$ sufficiently large, we have that
\[\sup_{y\in\mathbb Z^d}|E_{k,n}(y)|
\leq c(k) n^{-(d+k-1)/2}\leq c(M)n^{-(d+1)/2}(n\vee 2)^{-c'M^2}\leq c(M)n^{-(d+1)/2}e^{-c'|y|^2/n},
\]
which, together with \eqref{eq:dt-gradient-expansion}, \eqref{eq:grad-mvt} and \eqref{eq:u-grad}, concludes the proof of \eqref{eq:dt-hk-gradient}.

We are now ready to show \eqref{eq:ct-hk-gradient}. Let $\bar Z_n$ denote the process governed by $q_n(\cdot,\cdot)$ and let $N_t^r$ be a Poisson random variable with parameter $2rt$, independent of $\bar Z$. It's easy to see that $(\bar Z_{N^r_t})_{t\geq0}$ has the same law as $(Y_t)_{t\geq 0}$. Consequently,
\begin{equation}
	|P^{d,r}_0(Y_t=y+e)-P^{d,r}_0(Y_t=y)|\leq\sum_{n=0}^\infty P(N^r_t=n)|\nabla_e q_n(y)|
\end{equation}
Note that it suffices to consider $t\geq |y|$ such that $rt\geq 2$. Indeed, if $|y|\leq t\leq 2/r$, the bound in \eqref{eq:ct-hk-gradient} can be made trivially true by adjusting $\Cr{hkg}$ and $\Cr{hkg_exp}$ depending on $r$.
Let $I_t=\{n\geq1: rt/2\leq n\leq 4rt\}$. It follows from a standard large deviation estimate that $P(N^r_t\notin I_t)\leq ce^{-c'rt}$. Hence it follows from \eqref{eq:dt-hk-gradient} that 
\begin{align*}
\sum_{n=0}^\infty P(N^r_t=n)|\nabla_e q_n(y)|
&\leq c(d)\left(e^{-c'rt}+\sum_{n\in I_t}P(N^r_t=n) (rt)^{-(d+1)/2}
  e^{-c''(d)\frac{|y|^2}{rt}}\right)\\
&\leq c(d,r)\left(t^{-(d+1)/2}e^{-c'rt/2}+ (rt)^{-(d+1)/2}
  e^{-c''(d)\frac{|y|^2}{rt}}\right).
\end{align*}
Furthermore, since $|y|^2/t \le t$, it follows that $e^{-c' r t/2} \le e^{-c(r)|y|^2/t}$. Combining this with the last two displays completes the proof of \eqref{eq:ct-hk-gradient}.
\end{proof}

\bibliography{clean_bib}
\bibliographystyle{abbrv}

\end{document}